\documentclass[12pt]{article}
\newcommand{\nc}{\newcommand}
\nc{\rnc}{\renewcommand}
\rnc{\baselinestretch}{1.1}
\rnc{\arraystretch}{0.9}
\rnc{\thesubsubsection}{\arabic{subsubsection}.}
\usepackage{latexsym}
\usepackage{graphicx}
\newtheorem{lemma}{Lemma}
\newtheorem{theorem}[lemma]{Theorem}
\newtheorem{corollary}[lemma]{Corollary}
\rnc{\thelemma}{\arabic{lemma}.}
\rnc{\labelenumi}{(\arabic{enumi})}
\newcounter{subeq}
\nc{\be}{\begin{equation}}
\nc{\ee}{\end{equation}}
\nc{\bc}{\begin{center}}
\nc{\ec}{\end{center}}
\nc{\bpic}{\begin{picture}}
\nc{\epic}{\end{picture}}
\nc{\ba}[1]{\begin{array}{@{}#1@{}}}
\nc{\ea}{\end{array}}
\nc{\bea}{\begin{eqnarray}}
\nc{\eea}{\end{eqnarray}}
\nc{\er}[1]{(\ref{#1})}
\nc{\itm}[1]{\\\noindent$\bullet$\ \ #1}
\nc{\ds}{\displaystyle}
\nc{\ts}{\textstyle}
\rnc{\ss}{\scriptstyle}
\nc{\sss}{\scriptscriptstyle}
\nc{\ru}[1]{\rule[-#1ex]{0ex}{#1ex}}
\rnc{\vec}[1]{\mbox{\boldmath$#1$}}
\font\tenmsb=msbm10 scaled \magstep1
\font\sevenmsb=msbm7 scaled \magstep1
\font\fivemsb=msbm5 scaled \magstep1
\newfam\msbfam
\textfont\msbfam=\tenmsb
\scriptfont\msbfam=\sevenmsb
\scriptscriptfont\msbfam=\fivemsb
\def\Bbb#1{{\fam\msbfam\relax#1}}
\nc{\p}[2]{\makebox(0,0)[#1]{$#2$}}
\nc{\pp}[2]{\makebox(0,0)[#1]{$\ss#2$}}
\nc{\ppp}[2]{\makebox(0,0)[#1]{$\sss#2$}}
\nc{\text}[6]{\begin{picture}(#1,#2)\put(#3,#4){\p{#5}{\ds#6}}\end{picture}}
\rnc{\a}{\alpha}
\rnc{\b}{\beta}
\rnc{\d}{\delta}
\nc{\D}{\Delta}
\nc{\e}{\eta}
\nc{\g}{\gamma}
\rnc{\l}{\lambda}
\nc{\m}{\mu}
\rnc{\O}{\Omega}
\nc{\s}{\sigma}
\nc{\N}{\Bbb N}
\rnc{\P}{\Bbb P}
\nc{\Z}{\Bbb Z}
\nc{\mi}{\!-\!}
\nc{\pl}{\!+\!}
\rnc{\t}{\!\times\!}
\nc{\plmi}{\!\pm\!}
\nc{\mipl}{\!\mp\!}
\nc{\OP}{\mathrm{OP}}
\nc{\BP}{\mathrm{BP}}
\nc{\EM}{\mathrm{EM}}
\nc{\ASM}{\mathrm{ASM}}
\nc{\VOS}{\mathrm{VOS}}
\nc{\Par}{\mathrm{Par}}
\nc{\OT}{\mathrm{OT}}
\nc{\SYT}{\mathrm{SYT}}
\nc{\SSYT}{\mathrm{SSYT}}
\nc{\SSSYT}{\overline{\rule{0ex}{1.7ex}{\SSYT}}}
\nc{\GOT}{\mathrm{GOT}}
\nc{\NP}{\mathrm{NP}}
\rnc{\L}{\hat{L}}
\nc{\LL}{\check{L}}
\nc{\A}{A}
\rnc{\AA}{\mathcal{A}}

\begin{document}
\pagestyle{empty}
\bc
\vspace*{4mm}
\textbf{\large Osculating Paths and Oscillating Tableaux}\\
\bigskip\bigskip
\textbf{Roger E. Behrend}\\
\bigskip
\textit{School of Mathematics, Cardiff University,\\Cardiff, CF24 4AG, UK}\\
\smallskip
{\footnotesize\tt behrendr@cardiff.ac.uk}\\[20mm]
\begin{abstract}
\noindent The combinatorics of certain osculating lattice paths is studied,
and a relationship with oscillating tableaux is obtained.
More specifically, the paths being considered
have fixed start and end points on respectively the lower and right boundaries of
a rectangle in the square lattice, each path can take only unit steps
rightwards or upwards, and two different paths
are permitted to share lattice points, but not to cross or share lattice edges.
Such paths correspond to configurations of the six-vertex model of statistical mechanics with
appropriate boundary conditions, and they include cases which correspond to alternating sign matrices and
various subclasses thereof.
Referring to points of the rectangle through which no or two paths pass as vacancies
or osculations respectively, the case of primary interest
is tuples of paths with a fixed number~$l$ of vacancies and osculations.
It is then shown that there exist natural bijections which map each such path tuple~$P$ to a pair~$(t,\eta)$,
where~$\eta$ is an oscillating tableau of length~$l$
(i.e., a sequence of~$l\pl1$ partitions, starting with the empty partition, in which the Young diagrams
of successive partitions differ by a single square), and~$t$ is a certain, compatible sequence of~$l$
weakly increasing positive integers.  Furthermore, each vacancy or
osculation of~$P$ corresponds to a partition in~$\eta$ whose Young diagram is obtained from that of its predecessor
by respectively the addition or deletion of a square.
These bijections lead to enumeration formulae for osculating paths involving sums over oscillating tableaux.
\end{abstract}\ec

\hspace*{10.2mm}{\footnotesize\emph{Keywords:} osculating lattice paths, oscillating tableaux, alternating sign matrices\\
\hspace*{10.2mm}\emph{2000 Mathematics Subject Classification:} 05A15}

\newpage
\pagestyle{plain}
\subsubsection{Introduction}
The enumeration of nonintersecting lattice paths and of semistandard Young tab\-leaux are two basic
problems in combinatorics.  These problems are also closely related since there exist straightforward bijections
between certain tableaux and certain tuples of nonintersecting paths.
Furthermore, the problems are now well-understood, one reason being that a fundamental theorem,
often called the Lindstr\"{o}m-Gessel-Viennot theorem (see for
example~\cite[Theorem~1]{GesVie85},~\cite[Corollary~2]{GesVie89} or~\cite[Theorem~2.7.1]{Sta86}),
enables the cardinality of a set of tuples of such nonintersecting paths to be expressed as the determinant of a
matrix of binomial coefficients, thereby significantly elucidating and facilitating the enumeration.

More specifically, the paths in this context have fixed start and end points in the lattice~$\Z^2$, each path can take only unit steps
rightwards or upwards, and different paths cannot share any lattice point.  A (non-skew) semistandard Young tableaux
(see for example~\cite{Ful99},~\cite{Sag88},~\cite{Sag01} or~\cite[Ch.~7]{Sta99}) is
an array of positive integers which increase weakly from left to right along each row and increase strictly from top to bottom down each column,
and where the overall shape of the array corresponds to the Young diagram of a partition.
Apart from their intrinsic combinatorial interest, such tableaux are important in several other areas of mathematics,
including the representation theory of symmetric and general linear groups.
Each row of a tableau read from right to left itself constitutes a partition, and the usual bijections between
tableaux and nonintersecting paths (see for example~\cite[Sec.~6]{GesVie85},~\cite[Sec.~3]{GesVie89} or~\cite[Sec.~7.16]{Sta99})
essentially involve associating each row of a tableau with the path formed by the lower and right
boundary edges of the Young diagram of that row, and translated to a certain position in the lattice.
The condition that different paths cannot intersect then effectively corresponds to the condition that
the entries of a tableau increase strictly down columns.

It is also relevant here to consider standard Young tableaux and oscillating tableaux.
A standard Young tableau is a semistandard Young tableau with distinct entries which simply comprise~$1,2,\ldots,n$ for some~$n$,
while an oscillating tableau of length~$l$  (see for example~\cite{Ber86,Sag88,Sun86,Sun90})
is a sequence of $l\pl1$ partitions which starts with the empty partition,
and in which the Young diagrams of successive partitions differ by a single square.
It can be seen that a standard Young tableau~$\s$ corresponds naturally to an oscillating tableau~$\e$
in which each Young diagram is obtained from its predecessor by the addition of a square.  More precisely, if $\s_{ij}=k$, then
the Young diagram of the~$(k\pl1)$th partition of~$\e$ is obtained from that of
the~$k$th partition by the addition of a square in row~$i$ and column~$j$.
It can also be shown (as will be done for example in Section~18 of this paper)
that a semistandard Young tableau~$\tau$ corresponds naturally to a pair $(t,\e)$
in which~$t$ consists of the entries of~$\tau$ arranged as a weakly increasing sequence, and~$\e$ is an
oscillating tableau in which each Young diagram is obtained from its predecessor by
the addition of a square (i.e.,~$\e$ corresponds to a standard Young tableau).

The primary aim of this paper is to show that these results can essentially be generalized from
nonintersecting paths to osculating paths, and from pairs~$(t,\e)$ in which the oscillating tableau~$\e$
corresponds to a standard Young tableau to more general pairs~$(t,\e)$
in which each Young diagram of~$\e$ can be obtained from its predecessor
by either the addition or deletion of a square.

More specifically, osculating paths are those in which each path can still take only unit steps
rightwards or upwards in $\Z^2$, but for which two different paths are now permitted to share lattice points,
although not to cross or share lattice edges. Such paths correspond to configurations of the six-vertex model of statistical mechanics
(see for example~\cite[Ch.~8]{Bax82}).  The particular case being considered in this paper is that in which
the paths have fixed start and end points on respectively the lower and right boundaries of
a rectangle in~$\Z^2$.  Referring to points of the rectangle through which no or two paths pass as vacancies
or osculations respectively, the case of primary interest
will be path tuples with a fixed number~$l$ of vacancies and osculations.
It will then be found that there exist natural bijections which,
using data associated with the positions of the vacancies and osculations,
map any tuple~$P$ of such osculating paths to a pair~$(t,\eta)$, referred to as a generalized oscillating tableau,
in which~$\eta$ is an oscillating tableau of length~$l$, and~$t$ is a certain, compatible sequence of~$l$
weakly increasing positive integers.  A feature of these bijections is that
each vacancy or osculation of~$P$ corresponds to a partition in~$\eta$ whose
Young diagram is obtained from that of its predecessor by respectively the addition or deletion of a square.
If~$P$ is a tuple of nonintersecting paths, then there is such a bijection for which the
associated generalized oscillating tableau~$(t,\eta)$ corresponds to a semistandard Young tableau,
but the overall correspondence is in fact somewhat different from the usual ones known between nonintersecting paths
and semistandard Young tableaux.

A summary of the bijections between tuples of osculating paths and generalized oscillating tableaux
is given in Section~15.  A particular path tuple, which is shown in Figure~\ref{Ex}, is used as a running example
throughout the paper, while a further example involving sets of path tuples is given in Section~17.  An example of
a tuple of nonintersecting paths is given in Section~18.
The bijections lead to enumeration formulae for osculating paths involving sums over oscillating tableaux, these appearing
in Corollaries~14 and~17.

Much of the motivation for the work reported in this paper was derived from studies of alternating
sign matrices.  An alternating sign matrix, as first defined in~\cite{MilRobRum82,MilRobRum83},
is a square matrix in which each entry is $0$, $1$ or $-1$, each row and column contains at least one nonzero entry,
and along each row and column the nonzero entries alternate in sign, starting and finishing with a~1.
For reviews of alternating sign matrices and related subjects, see for example~\cite{Bre99,BrePro99,Pro01,Rob91a}.
Of particular relevance here is that there exist straightforward bijections between alternating sign matrices, or certain subclasses
thereof, and certain tuples of osculating paths in a rectangle (see for example Section~4 of this paper and references therein).
Relatively simple enumeration formulae are known for such cases,
but all currently-known derivations of these formulae, as given in~\cite{ColPro05,Fis07,Kup96,Kup02,Oka06,Zei96a,Zei96b},
are essentially non-combinatorial in nature.  Furthermore, it is known that
the numbers of $n\t n$ alternating sign matrices, descending plane partitions with no part larger than~$n$~(see for
example~\cite{And79,Kra06,Lal03,MilRobRum82,MilRobRum83}),
and totally symmetric self-complementary plane partitions in a $2n\t2n\t2n$ box~(see for
example~\cite{And94,Ish06a,Ish06b,MilRobRum86,Ste90})
are all equal, and further equalities between the cardinalities of
certain subsets of these three objects have been conjectured or in a few cases proved, but
no combinatorial proofs of these equalities are known.  It is therefore hoped that the bijections between
osculating paths and generalized oscillating tableaux described in this paper may eventually lead to an improved combinatorial
understanding of some of these matters.

Osculating paths have also appeared in a number of recent studies as a special case of
friendly walkers (see for example \cite{Bou06,Ess03,GutVog02,KraGutVie03} and references therein).  However,
all of these cases use a different external configuration from the rectangle being used here.  In
particular, the paths start and end on two parallel lines rotated by $45^\circ$ with respect to the rows
or columns of the square lattice.  A general enumeration formula for such osculating paths has been conjectured in~\cite{Bra97}.

\emph{Notation.} \ Throughout this paper, $\P$ denotes the set of positive integers, $\N$ denotes the set of nonnegative integers,
$[m,n]$ denotes the set $\{m,m\pl1,\ldots,n\}$ for any $m,n\in\Z$, with $[m,n]=\emptyset$ for $n<m$, and $[n]$ denotes the set
$[1,n]$ for any $n\in\Z$.
For a finite set~$T$, $|T|$ denotes the cardinality of~$T$.
For a condition~$C$, $\d_C$ denotes a function which is
$1$ if $C$ is satisfied and $0$ if not, and for numbers $i$ and $j$, $\d_{ij}$ denotes
the usual Kronecker delta, $\d_{ij}=\d_{i=j}$.  For a positive odd integer $n$, the double
factorial is $n!!=n(n\mi2)(n\mi4)\ldots3.1$, while (-1)!! is taken to be 1.

\subsubsection{Osculating Paths}
In this section, the set of tuples of osculating lattice paths in a fixed $a$ by $b$ rectangle, with the paths starting at
points (specified by a subset $\{\b_1,\ldots,\b_r\}$ of $[b]$) along the lower boundary, ending at points
(specified by a subset $\{\a_1,\ldots,\a_r\}$ of $[a]$) along the right boundary, and
taking only unit steps rightwards or upwards, will be defined precisely.

For any $a,b\in\P$, the subset $[a]\t[b]$ of $\Z^2$ will
be regarded diagrammatically as a rectangle of lattice points with rows labeled
$1$ to $a$ from top to bottom, columns labeled $1$ to $b$ from left to right, and $(i,j)$ being the point in row~$i$ and column $j$.
The motivation for using this labeling is that it will provide consistency with the standard labeling of rows and columns
of matrices and Young diagrams, both of which will later be associated with path tuples.

\begin{figure}[h]
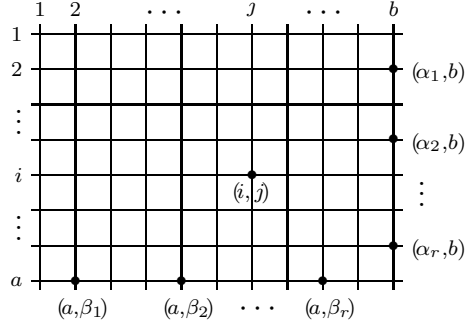
\centering\bpic(12,9.3)
\multiput(0.75,1)(0,1){8}{\line(1,0){10.5}}\multiput(1,0.75)(1,0){11}{\line(0,1){7.5}}
\put(1,8.5){\pp{b}{1}}\put(2,8.5){\pp{b}{2}}\put(7,8.5){\pp{b}{j}}\put(11,8.5){\pp{b}{b}}
\put(0.5,8){\pp{r}{1}}\put(0.5,7){\pp{r}{2}}\put(0.5,4){\pp{r}{i}}\put(0.5,1){\pp{r}{a}}
\put(11.5,6.9){\pp{l}{(\!\a_1,b)}}\put(11.5,4.9){\pp{l}{(\!\a_2,b)}}\put(11.5,1.9){\pp{l}{(\!\a_r,b)}}
\put(2.2,0.5){\pp{t}{(\!a,\b_1\!)}}\put(5.2,0.5){\pp{t}{(\!a,\b_2\!)}}\put(9.2,0.5){\pp{t}{(\!a,\b_r\!)}}
\multiput(4.51,8.6)(4.5,0){2}{\p{}{\cdots}}\multiput(0.4,2.73)(0,3){2}{\p{}{\vdots}}
\put(7,4){\pp{}{\bullet}}\put(6.97,3.8){\pp{t}{(\!i,\,j\!)}}
\put(11.8,3.7){\p{}{\vdots}}\put(7.15,0.2){\p{}{\cdots}}
\put(2,1){\pp{}{\bullet}}\put(5,1){\pp{}{\bullet}}\put(9,1){\pp{}{\bullet}}
\put(11,2){\pp{}{\bullet}}\put(11,5){\pp{}{\bullet}}\put(11,7){\pp{}{\bullet}}\epic
\caption{Labeling of the lattice and boundary points.\label{Lat}}\end{figure}

The general labeling of the lattice, together with the start and end points of paths, is shown diagrammatically in Figure~\ref{Lat}.
\setlength{\unitlength}{4.7mm}

For $\a\in[a]$ and $\b\in[b]$, let $\Pi(a,b,\a,\b)$ be the set of all paths from $(a,\b)$ to $(\a,b)$, in which
each step of any path is $(0,1)$ or $(-1,0)$,
\begin{equation}\ba{r@{}l}\Pi(a,b,\a,\b)\;:=&\;\;
\Bigl\{\Bigl((i_0,j_0)\!=\!(a,\b),(i_1,j_1),\ldots,(i_{L\mi1},j_{L\mi1}),(i_L,j_L)\!=\!(\a,b)\Bigr)\;\Big|\\[2.4mm]
&\hspace{16mm}(i_l,j_l)\mi(i_{l\mi 1},j_{l\mi 1})\in\{(0,1),(-1,0)\}\mbox{ for each }l\in[L]\Bigr\}\,,\ea\end{equation}
where necessarily $L=a\mi\a\pl b\mi\b$.
It follows that $|\Pi(a,b,\a,\b)|=\biggl(\ba{c}a\mi\a\pl b\mi\b\\a\mi\a\ea\biggr)$.

For $\a,\a'\in[a]$ and $\b,\b'\in[b]$, with $\a<\a'$ and $\b<\b'$, paths $P\in\Pi(a,b,\a,\b)$ and $P'\in\Pi(a,b,\a',\b')$
are said to be \emph{osculating} if they
do not cross or share lattice edges, but possibly share lattice points.
More precisely, this means that if $P_l=P'_{l'}=(i,j)$ for some $l$ and $l'$
(which implies that $l=a\mi i\pl j\mi\b$, $l'=a\mi i\pl j\mi\b'$),
then $P_{l\mi1}=(i,j\mi1)$, $P_{l\pl1}=(i\mi1,j)$,
$P'_{l'\mi 1}=(i\pl1,j)$ (if $l'\ne0$) and $P'_{l'\pl 1}=(i,j\pl 1)$ (if $l'\ne a\mi\a'\pl b\mi\b'$).
Any such common point $(i,j)$ will be referred to as an \emph{osculation} of~$P$.

For $r\in[0,\min(a,b)]$, $\a=\{\a_1,\ldots,\a_r\}\subset[a]$ and $\b=\{\b_1,\ldots,\b_r\}\subset[b]$, with
$\a_1\!<\!\ldots\!<\!\a_r$ and $\b_1\!<\!\ldots\!<\!\b_r$,
let $\OP(a,b,\a,\b)$ be the set of $r$-tuples
of pairwise osculating paths in which the $k$-th path is in $\Pi(a,b,\a_k,\b_k)$,
\begin{equation}\label{OP}\ba{r@{}l}\OP(a,b,\a,\b)\;:=&\;\;
\Bigl\{P\!=\!(P_1,\ldots,P_r)\in\Pi(a,b,\a_1,\b_1)\t\ldots\t\Pi(a,b,\a_r,\b_r)\;\Big|\\[2.4mm]
&\hspace{22mm}\mbox{$P_k$ and $P_{k\pl 1}$ are osculating for each $k\in[r\mi1]$}\Bigr\}\,.\ea\end{equation}

Also, for any $a,b\in\P$, let $\BP(a,b)$ be the set of all pairs $(\a,\b)$ of \emph{boundary points},
\begin{equation}\label{BP}\ba{r@{}l}\BP(a,b)\;:=&\;\;
\Bigl\{(\a,\b)\;\Big|\;\a\subset[a],\;\b\subset[b],\;|\a|=|\b|\Bigr\}\,.\ea\end{equation}
It follows that $|\BP(a,b)|=
\!\ds\sum_{r=0}^{\min(a,b)}\!\biggl(\ba{c}a\\r\ea\biggr)\biggl(\ba{c}b\\r\ea\biggr)=\biggl(\ba{c}a\pl b\\a\ea\biggr)$.

Throughout the remainder of this paper, $a$ and $b$ will be used to denote positive integers,
corresponding to the dimensions of a rectangle of lattice points, and $(\a,\b)$ will
denote an element of $\BP(a,b)$.

Now let $\OP(a,b)$ be the set of all tuples of osculating paths in $[a]\t[b]$ with any boundary points,
\begin{equation}\OP(a,b):=\bigcup_{(\a,\b)\in\BP(a,b)}\OP(a,b,\a,\b)\,.\end{equation}

For $P\in\OP(a,b)$, any point $(i,j)\in[a]\t[b]$ through which no path of $P$ passes
will be referred to as a \emph{vacancy} of $P$.
Define $N(P)$ to be the set of all vacancies of $P$,
$X(P)$ to be the set of all osculations of~$P$ and $\chi(P)$ to be the number of osculations of~$P$,
\begin{equation}\label{chi}\chi(P):=|X(P)|\,.\end{equation}
A tuple of \emph{nonintersecting paths} is any $P$ for which $X(P)=\emptyset$.  Nonintersecting paths will
be considered in more detail in Section~18.

Define also a \emph{vacancy-osculation} of $P\in\OP(a,b)$ as either a vacancy or osculation of $P$,
and the \emph{vacancy-osculation set} $Z(P)$ as the set of all
vacancy-osculations of~$P$,
\begin{equation}\label{OPToVOS}Z(P)\;:=\;N(P)\cup X(P)\,.\end{equation}
In other words, $Z(P)$ is the set of points of $[a]\t[b]$ through which either zero or two paths of~$P$ pass.
It will be of particular interest to consider sets of path tuples with~$l$ vacancy-osculations, for fixed $l\in\N$,
\begin{equation}\label{OPl}\ba{l}\OP(a,b,l)\::=\:\Bigl\{P\in\OP(a,b)\;\Big|\;|Z(P)|=l\Bigr\}\\[3mm]
\OP(a,b,\a,\b,l)\::=\:\Bigl\{P\in\OP(a,b,\a,\b)\;\Big|\;|Z(P)|=l\Bigr\}\,,\ea\end{equation}
a primary aim of this paper being to study the properties and cardinality of $\OP(a,b,$ $\a,\b,l)$.

Finally, note that there are trivial bijections, involving reflection or translation, between
certain sets of path tuples.  More precisely, using $\approx$ to denote the existence of a bijection between sets,
\begin{equation}\label{OPRT}\ba{l}\OP(a,b,\a,\b,l)\;\approx\;\OP(b,a,\b,\a,l)\\[3mm]
\qquad\approx\;\OP(\bar{a}\pl a,\bar{b}\pl
b,\{\bar{a}\pl\a_1,\ldots,\bar{a}\pl\a_r\},\{\bar{b}\pl\b_1,\ldots,\bar{b}\pl\b_r\},l\pl\bar{a}\,\bar{b}
\pl\bar{a}\,b\pl a\,\bar{b})\,,\ea\end{equation}
for any $a,b\in\P$, \ $\bar{a},\bar{b}\in\N$ \ and \ $(\a,\b)=(\{\a_1,\ldots,\a_r\},\{\b_1,\ldots,\b_r\})\in\BP(a,b)$.
For the first bijection of~(\ref{OPRT})
each path is reflected in the main diagonal of the lattice, while for the second bijection of~(\ref{OPRT})
each path is translated by $(\bar{a},\bar{b})$.
\setlength{\unitlength}{9.43mm}
\begin{figure}[b]\centering\bpic(0.3,3.6)\put(0,0.15){\pp{l}{4}}\put(0,1.15){\pp{l}{3}}\put(0,2.15){\pp{l}{2}}\put(0,3.15){\pp{l}{1}}
\put(0.45,3.45){\pp{b}{1}}\put(1.45,3.45){\pp{b}{2}}\put(2.45,3.45){\pp{b}{3}}
\put(3.45,3.45){\pp{b}{4}}\put(4.45,3.45){\pp{b}{5}}\put(5.45,3.45){\pp{b}{6}}\epic
\includegraphics[width=50mm]{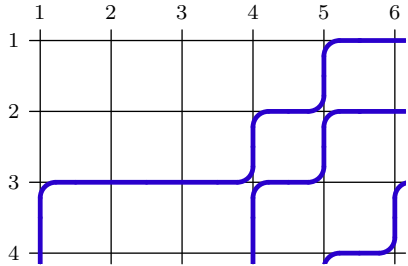}\caption{Example of a tuple of osculating paths.\label{Ex}}\end{figure}

An example of an element of $\OP(4,6,\{1,2,3\},\{1,4,5\},11)$ is
$P=\Bigl(\Bigl((4,1)$, $(3,1)$, $(3,2)$, $(3,3)$, $(3,4)$, $(2,4)$, $(2,5)$, $(1,5)$, $(1,6)\Bigr)$,
$\Bigl((4,4)$, $(3,4)$, $(3,5)$, $(2,5)$, $(2,6)\Bigr),$
$\Bigl((4,5)$, $(4,6)$, $(3,6)\Bigr)\Bigr)$,\ru{1.8}
which is shown diagrammatically in Figure~\ref{Ex}.

For this case, $N(P)=\{(1,1)$, $(1,2)$, $(1,3)$, $(1,4)$, $(2,1)$, $(2,2)$, $(2,3)$, $(4,2)$, $(4,3)\}$,\ru{1.8}
$X(P)=\{(2,5),$ $(3,4)\}$ and $\chi(P)=2$.
This will serve as a running example throughout this paper.

\subsubsection{Edge Matrices}
In this section, it will be seen that each tuple of osculating paths corresponds naturally
to a pair $(H,V)$ of $\{0,1\}$ matrices, which will be referred to as \emph{horizontal} and \emph{vertical edge matrices}.

For $P\in\OP(a,b,\a,\b)$, the correspondence is given simply by the rule that $H_{ij}$ is $0$ or $1$ according to whether or not
$P$ contains a path which passes from $(i,j)$ to $(i,j\pl1)$, and that
$V_{ij}$ is $0$ to $1$ according to whether or not
$P$ contains a path which passes from $(i\pl1,j)$ to $(i,j)$.
Thus $H_{ij}$ is associated with the horizontal lattice edge between $(i,j)$ and $(i,j\pl1)$,
and
$V_{ij}$ is associated with the vertical lattice edge between $(i,j)$ and $(i\pl1,j)$.
It is also convenient to consider boundary edges horizontally between $(i,0)$ and $(i,1)$,
and between $(i,b)$ and $(i,b\pl1)$, for each $i\in[a]$, and vertically
between $(0,j)$ and $(1,j)$, and between $(a,j)$ and $(a\pl1,j)$, for each $j\in[b]$, and
to include in each path $P_k$ the additional points $(a\pl1,\b_k)$ at the start and $(\a_k,b\pl1)$
at the end.
Each point $(i,j)\in[a]\t[b]$ can then be associated with a \emph{vertex configuration} which
involves the four values $H_{i,j\mi1}$, $V_{ij}$, $H_{ij}$ and $V_{i\mi1,j}$, this
being depicted diagrammatically as
\setlength{\unitlength}{6mm}
\begin{equation}\label{VC}
\raisebox{-1.5\unitlength}[1.8\unitlength][1.5\unitlength]{
\bpic(4.2,3)\put(2.2,0.5){\line(0,1){2}}\put(1.2,1.5){\line(1,0){2}}
\multiput(2.2,0.5)(0,2){2}{\ppp{}{\bullet}}\multiput(1.2,1.5)(2,0){2}{\ppp{}{\bullet}}
\put(-0.2,1.4){\ppp{l}{H_{i,j\mi1}}}\put(4.1,1.4){\ppp{r}{H_{ij}}}
\put(2.4,3.08){\ppp{t}{V_{i\mi1,j}}}\put(2.3,-0.15){\ppp{b}{V_{ij}}}\epic}\,.
\end{equation}

It can be seen that for any tuple of osculating paths
there are only six possible path configurations surrounding any lattice point, given diagrammatically as:
\setlength{\unitlength}{1.2mm}
\begin{equation}\raisebox{-5\unitlength}[7.9\unitlength][7\unitlength]{\includegraphics[width=130mm]{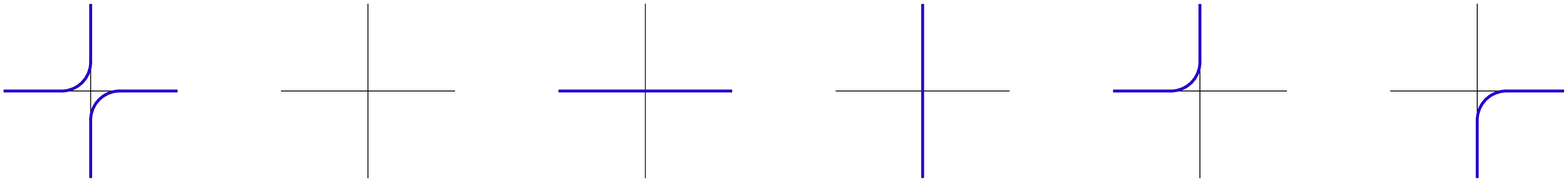}}\end{equation}
Correspondingly, there are six possible vertex configurations:
\begin{equation}\label{6V}\raisebox{-11\unitlength}[9\unitlength][9\unitlength]{
\bpic(100,19)\multiput(5,7)(18,0){6}{\line(0,1){10}}\multiput(0,12)(18,0){6}{\line(1,0){10}}
\multiput(5,7)(18,0){6}{\ppp{}{\bullet}}\multiput(5,17)(18,0){6}{\ppp{}{\bullet}}
\multiput(0,12)(18,0){6}{\ppp{}{\bullet}}\multiput(10,12)(18,0){6}{\ppp{}{\bullet}}
\put(5.1,6){\ppp{t}{1}}\put(5.1,18){\ppp{b}{1}}\put(-0.7,12){\ppp{r}{1}}\put(10.7,12){\ppp{l}{1}}
\put(23.1,6){\ppp{t}{0}}\put(23.1,18){\ppp{b}{0}}\put(17.3,12){\ppp{r}{0}}\put(28.7,12){\ppp{l}{0}}
\put(41.1,6){\ppp{t}{0}}\put(41.1,18){\ppp{b}{0}}\put(35.3,12){\ppp{r}{1}}\put(46.7,12){\ppp{l}{1}}
\put(59.1,6){\ppp{t}{1}}\put(59.1,18){\ppp{b}{1}}\put(53.3,12){\ppp{r}{0}}\put(64.7,12){\ppp{l}{0}}
\put(77.1,6){\ppp{t}{0}}\put(77.1,18){\ppp{b}{1}}\put(71.3,12){\ppp{r}{1}}\put(82.7,12){\ppp{l}{0}}
\put(95.1,6){\ppp{t}{1}}\put(95.1,18){\ppp{b}{0}}\put(89.3,12){\ppp{r}{0}}\put(100.7,12){\ppp{l}{1}}
\put(5,0){\pp{b}{1}}\put(23,0){\pp{b}{2}}\put(41,0){\pp{b}{3}}\put(59,0){\pp{b}{4}}
\put(77,0){\pp{b}{5}}\put(95,0){\pp{b}{6}}
\epic}\end{equation}

The numbers below each vertex configuration will be used to label the six possible types.
Thus, type 1 corresponds to an osculation, and type 2 to a vacancy.

It can also be seen that the six cases of~\er{6V} correspond exactly to the simple but important condition
\begin{equation}\label{AC}H_{i,j\mi1}+V_{ij}=V_{i\mi1,j}+H_{ij}\,,\end{equation}
for each $(i,j)\in[a]\t[b]$.

Accordingly, taking into account all of the previous considerations, sets of edge matrix pairs
for $a,b\in\P$ and $(\a,\b)\in\BP(a,b)$ are defined as
\begin{equation}\label{EM2}
\ba{r@{}l}\EM(a,b):=\Bigl\{(H,V)\,\Big|\,\:&\bullet\mbox{ $H$ and $V$ are matrices with all entries in $\{0,1\}$}\\[1.5mm]
&\bullet\mbox{ $H$ has rows labeled by $[a]$, columns labeled by $[0,b]$}\\[1.5mm]
&\bullet\mbox{ $V$ has rows labeled by $[0,a]$, columns labeled by $[b]$}\\[1.5mm]
&\bullet\mbox{ $H_{i0}=0$ for all $i\in[a]$, \ $V_{0j}=0$ for all $j\in[b]$}\\[1.5mm]
&\bullet\mbox{ $H_{i,j\mi1}+V_{ij}=V_{i\mi1,j}+H_{ij}$ \ for all $(i,j)\in[a]\t[b]$}\,\Bigr\}\ea\end{equation}
and
\begin{equation}\label{EM4}
\ba{l}\EM(a,b,\a,\b)\;:=\\[2mm]
\quad\Bigl\{(H,V)\in\EM(a,b)\,\Big|\,\:
H_{ib}=\d_{i\in\a}\mbox{ for all } i\in[a], \ V_{aj}=\d_{j\in\b}\mbox{ for all }j\in[b]\Bigr\}\,.\ea
\refstepcounter{equation}\label{HV}\addtocounter{equation}{-1}
\end{equation}
It can be seen that the `boundary conditions' on $(H,V)\in\EM(a,b,\a,\b)$
are that the first column of $H$ and first row of $V$ are zero, and that
the last column of $H$ and last row of $V$ are specified by $\a$ and $\b$ respectively.

As already indicated, for any $P=(P_1,\ldots,P_r)\in\OP(a,b,\a,\b)$,
a corresponding $(H,V)\in\EM(a,b,\a,\b)$ is given by
\rnc{\theequation}{\arabic{equation}\alph{subeq}}\setcounter{subeq}{1}
\begin{equation}\ba{rcl}H_{ij}&=&\left\{\ba{ll}1,&\mbox{$(P_k)_l=(i,j)$ and $(P_k)_{l\pl1}=(i,j\pl1)$ for some $k$ and $l$,}\\[1.5mm]
&\mbox{ or $i\in\a$ and $j=b$}\\[2.5mm]
0,&\mbox{otherwise}\ea\right.\ea\end{equation}
\begin{equation}\stepcounter{subeq}\addtocounter{equation}{-1}\ba{rcl}
V_{ij}&=&\left\{\ba{ll}1,&\mbox{$(P_k)_l=(i\pl1,j)$ and $(P_k)_{l\pl1}=(i,j)$ for some $k$ and $l$,}\\[1.5mm]
&\mbox{ or $i=a$ and $j\in\b$}\\[2.5mm]
0,&\mbox{otherwise\,,}\ea\right.\ea\end{equation}
and it can be seen straightforwardly that this mapping is a bijection between
$\OP(a,b,$ $\a,\b)$ and $\EM(a,b,\a,\b)$.\rnc{\theequation}{\arabic{equation}}

It can also be seen that the number of osculations of $P$ can be expressed in terms of
the corresponding $(H,V)$ as
\begin{equation}\label{chiHV}\chi(P)\;=\;\sum_{(i,j)\in[a]\t[b]}\!\!H_{i,j\mi1}\,V_{ij}\;=\;
\sum_{(i,j)\in[a]\t[b]}\!\!V_{i\mi1,j}\,H_{ij}\,,\end{equation}
since $(i,j)$ is an osculation if and only if
$H_{i,j\mi1}=V_{ij}=1$, or equivalently $V_{i\mi1,j}=H_{ij}=1$.

Returning to the running example, the edge matrices are
\begin{equation}\label{ExHV}\left(\ba{c@{\;}c@{\;}c}H_{10}&\ldots&H_{16}\\\vdots&&\vdots\\H_{40}&\ldots&H_{46}\ea\right)\!\!=\!\!
\left(\ba{ccccccc}0&0&0&0&0&1&1\\
0&0&0&0&1&1&1\\
0&1&1&1&1&0&1\\
0&0&0&0&0&1&0\ea\right)\!,\:
\left(\ba{c@{\;}c@{\;}c}V_{01}&\ldots&V_{06}\\\vdots&&\vdots\\V_{41}&\ldots&V_{46}\ea\right)\!\!=\!\!
\left(\ba{cccccc}0&0&0&0&0&0\\
0&0&0&0&1&0\\0&0&0&1&1&0\\
1&0&0&1&0&1\\1&0&0&1&1&0\ea\right)\!.\!\end{equation}

The edge matrix representation of osculating paths
corresponds to the standard representation of configurations of the \emph{six-vertex} or \emph{square ice} lattice model in
statistical mechanics (see for example~\cite[Ch.~8]{Bax82} and references therein).  In this model
$H_{ij}=0$ is represented by a leftward arrow on the corresponding lattice edge,
$H_{ij}=1$ by a rightward arrow, $V_{ij}=0$ by a downward arrow, and
$V_{ij}=1$ by an upward arrow.  The osculating paths then follow the rightward and upward arrows,
and condition~\er{AC} corresponds to \emph{arrow conservation} at each lattice point
(i.e., the numbers of arrows into and out of each point are equal).  One of the main quantities of interest for such statistical
mechanical models is the \emph{partition function}, which is a certain weighted sum over the configurations of the model.
The particular case being considered here is that of configurations of the six-vertex model on an~$a$ by~$b$ rectangle
with fixed boundary conditions in which on the upper boundary all arrows point down,
on the left boundary all arrows point left,
on the lower boundary arrows point up or down according to whether or not their position is in~$\b$, and
on the right boundary arrows point right or left according to whether or not their position is in~$\a$.
Note that the six-vertex model has been extensively studied with a variety of boundary conditions.
See for example~\cite[Ch.~8~\&~9]{Bax82} for the details of studies with periodic (i.e., toroidal) boundary conditions,
and~\cite{BatBaxOroYun95,IzeCokKor92,OwcBax89,YunBat95,ZinJus02} for some studies with other boundary conditions.

\subsubsection{Alternating Sign Matrices}
In this section, it will be seen that each tuple of osculating paths also corresponds naturally
to a $\{-1,0,1\}$ matrix, which will be referred to as an alternating sign matrix.
Although this representation of osculating paths will not be used in obtaining any of the results of
this paper, it is introduced here in order to
present some known formulae for the cardinality of special cases of $\OP(a,b,\a,\b)$
which are usually given in the context of such a representation.
Also, as indicated in Section~1, much of the motivation for the work reported in this paper was derived from the
studies of alternating sign matrices in which these formulae were obtained.
Note, however, that the enumeration results obtained later in this paper
apply to cases of $\OP(a,b,\a,\b,l)$, i.e., for which the path tuples all have $l$ vacancy-osculations,
whereas the enumeration formulae listed in this section apply to certain cases of $\OP(a,b,\a,\b)$, i.e.,
for which there is no restriction on the number of vacancy-osculations.
Note also that the alternating sign matrices defined in other papers comprise a special case of the alternating
sign matrices defined here, and so will be referred to here as `standard alternating sign matrices'.

For any $a,b\in\P$ and $(\a,\b)\in\BP(a,b)$ define the associated set of \emph{alternating sign matrices} as
\begin{equation}\ba{l}\ASM(a,b,\a,\b)\;:=\\[2.5mm]
\quad\quad\ba{r@{}l}\Bigl\{\;A\;\:\Big|\:\;&\bullet\mbox{ $A$ is an $a\t b$ matrix with all entries in $\{-1,0,1\}$}\\[1.5mm]
&\bullet\mbox{ along each row and column of $A$ the nonzero entries, if}\\[0.5mm]
&\mbox{ \ \ there are any, alternate in sign starting with a $1$}\\[1.9mm]
&\bullet\mbox{ $\sum_{j=1}^bA_{ij}=\d_{i\in\a}$ \ for all $i\in[a]$}\\[1.9mm]
&\bullet\mbox{ $\sum_{i=1}^aA_{ij}=\d_{j\in\b}$ \ for all $j\in[b]$}\;\Bigr\}\,.\ea\ea\end{equation}

For any $(H,V)\in\EM(a,b,\a,\b)$, a corresponding $A\in\ASM(a,b,\a,\b)$
is given simply by
\begin{equation}\label{HToASM}A_{ij}\;=\;H_{ij}-H_{i,j\mi1}\end{equation}
or equivalently, due to~\er{AC},
\begin{equation}\label{VToASM}A_{ij}\;=\;V_{ij}-V_{i\mi1,j}\,,\end{equation}
for each $(i,j)\in[a]\t[b]$ (i.e., $A$ is the column difference matrix of $H$ and row difference matrix of $V$).
Note that under this mapping vertex configurations 1--4 of~\er{6V} at $(i,j)$ give $A_{ij}=0$, configuration 5 gives
$A_{ij}=-1$ and configuration 6 gives $A_{ij}=1$.

It can be checked that the mapping~(\ref{HToASM},\ref{VToASM}) is a bijection between
$\EM(a,b,\a,\b)$ and $\ASM(a,b,\a,\b)$, and that the inverse mapping is
\begin{equation}\label{ASMToHV}\ba{c}\ds H_{ij}=\sum_{j'=1}^jA_{ij'}\,,\mbox{ \ for each }(i,j)\in[a]\t[0,b]\\[6mm]
\ds V_{ij}=\sum_{i'=1}^iA_{i'j}\,,\mbox{ \ for each }(i,j)\in[0,a]\t[b]\ea\end{equation}
(i.e., $H$ and $V$ are respectively the partial column and row sum matrices of $A$).

It follows, using~\er{chiHV} and~\er{ASMToHV}, that the number of osculations of any $P\in\OP(a,b)$ can be expressed in terms of
the corresponding $A\in\ASM(a,b)$ as
\begin{equation}\label{chiA}\chi(P)\;=\sum_{\stackrel{\ru{0.8}\{(i,i'\!,j,j')\in[a]^2\times\![b]^2\;|}{
\ss\qquad\qquad i'\ge i,\;j'<j\}}}\!\!\!A_{ij}\,A_{i'j'}\;=
\sum_{\stackrel{\ru{0.8}\{(i,i'\!,j,j')\in[a]^2\times\![b]^2\;|}{\ss\qquad\qquad i'>i,\;j'\le j\}}}\!\!\!A_{ij}\,A_{i'j'}\,.\end{equation}

The alternating sign matrix for the running example is
\begin{equation}A\;=\;\left(\ba{cccccc}
0&0&0&0&1&0\\0&0&0&1&0&0\\
1&0&0&0&-1&1\\0&0&0&0&1&-1\ea\right).\end{equation}

Five previously-studied cases of alternating sign matrices will now be considered, for any $n\in\P$:\\[1.7mm]
\mbox{}$\;\;\bullet\;\;\ASM(n,n,[n],[n])$\\[1mm]
\mbox{}$\;\;\bullet\;\;\ASM(n,n\pl1,[n],[n\pl1]\!\setminus\!\{n\pl1\mi m\})$, \ for $m\in[0,n]$\\[1mm]
\mbox{}$\;\;\bullet\;\;\ASM(n,n\pl m,[n],[n\mi1]\cup\{n\pl m\})$, \ for $m\in\N$\\[1mm]
\mbox{}$\;\;\bullet\;\;\ASM(n,2n\mi1,[n],\{1,3,5,\ldots,2n\mi1\})$\\[1mm]
\mbox{}$\;\;\bullet\;\;\ASM(n,n,\{1,3,5,\ldots,2\lceil\frac{n}{2}\rceil\mi1\},\{1,3,5,\ldots,2\lceil\frac{n}{2}\rceil\mi1\})$\\[1.7mm]
It can be seen that the third case reduces to the first case for $m=0$.

The elements of $\ASM(n,n,[n],[n])$ will be referred to here as
\emph{standard alternating sign matrices}. They are simply $n\t n$
$\{-1,0,1\}$ matrices in which along each row and column the sum of entries is 1, and the nonzero entries
alternate in sign.  They were introduced in~\cite{MilRobRum82,MilRobRum83},
in which an enumeration formula was conjectured which gives
\begin{equation}\label{NASM}|\OP(n,n,[n],[n])|\;=\;\:\prod_{i=0}^{n\mi1}\!\frac{(3i\pl1)!}{(n\pl i)!}\,.\end{equation}
This was eventually proved in~\cite{Zei96a} and, using a different method,~\cite{Kup96}.
The correspondence between standard alternating sign matrices and edge matrices was first
identified in~\cite{RobRum86}, and is also discussed, at least in the statistical mechanical model
version, in~\cite{BouHab95,ElkKupLarPro92,Pro01}. The correspondence between standard alternating sign matrices and
osculating paths is also considered in~\cite[Sec.~5]{BouHab95},~\cite[Sec.~2]{Bra97},~\cite[Sec.~9]{EgeRedRya01} and~\cite[Sec.~IV]{Tam01}.

The six-vertex model boundary conditions for
this case, in which all arrows on the upper and lower boundaries point into the square, and
all arrows on the left and right boundaries point out of the square, are known as \emph{domain wall boundary conditions}
(see for example~\cite{IzeCokKor92}).

It can be seen that the $n\t n$ permutation matrices are included in
$\ASM(n,n,[n],[n])$, and, from~(\ref{chiA}), that if $A$ is a
permutation matrix which corresponds to $P\in\OP(n,n,[n],[n])$, then the number of osculations
$\chi(P)$ is simply the inversion number of $A$.
The inversion number of any $A\in\ASM(a,b)$ could be defined as
$\mathcal{I}(A)=\ru{2}\sum_{\{(i,i'\!,j,j')\in[a]^2\t[b]^2\,|\,i'>i,\;j'<j\}}A_{ij}\,A_{i'j'}$,
this being consistent with a definition of inversion number for standard alternating sign
matrices given in~\cite{MilRobRum83}.  It then follows (see for example~\cite[Theorem~2c]{RobRum86})
that $\mathcal{I}(A)=\chi(P)+\m(A)$, where $P\in\OP(a,b)$ corresponds to $A$, and
$\m(A)$ is the number of $-1$ entries in $A$.

It can be seen that in any standard alternating sign matrix,
there is a single 1 in each of the first and last row and column.
Furthermore, using elementary symmetry considerations, it follows that the number of
$(n\pl1)\t(n\pl1)$ standard alternating sign matrices $A$
with $A_{ij}=1$ is the same for $(i,j)$ taken to be any of the eight cases $(1,m\pl1)$, $(1,n\pl1\mi m)$,
$(m\pl1,1)$, $(n\pl1\mi m,1)$, $(m\pl1,n\pl1)$, $(n\pl1\mi m,n\pl1)$, $(n\pl1,m\pl1)$ or $(n\pl1,n\pl1\mi m)$, for fixed $m\in[0,n]$.
Focussing on the case $(n\pl1,n\pl1\mi m)$, i.e., matrices with the 1 of their last row in column $n\pl1\mi m$, so that
$A_{n\pl1,j}=\d_{j,n\pl1\mi m}$, it follows from~\er{ASMToHV}
that $H_{n\pl1,j}=\d_{j\ge n\pl1\mi m}$, and then from $V_{n\pl1,j}=1$ and~\er{AC} that
$V_{nj}=\d_{j\ne n\pl1\mi m}$, which implies that the submatrices formed by the first $n$ rows of such~$A$ comprise
$\ASM(n,n\pl1,[n],[n\pl1]\!\setminus\!\{n\pl1\mi m\})$.  Obtaining the cardinality of this and the seven
related sets is known as
\emph{refined alternating sign matrix enumeration}, and, following conjectures in~\cite{MilRobRum82,MilRobRum83},
was achieved in~\cite{Zei96b} and more recently in~\cite{ColPro05,Fis07} giving
\begin{equation}|\OP(n,n\pl1,[n],[n\pl1]\!\setminus\!\{n\pl1\mi m\})|\;=\;
\frac{(2n\mi m)!\,(n\pl m)!}{n!\:m!\:(n\mi m)!}\;\prod_{i=1}^{n}\!\frac{(3i\mi2)!}{(n\pl i)!}\,.\end{equation}

The related case of $\ASM(n,n\pl m,[n],[n\mi1]\!\cup\!\{n\pl m\})$ was considered in~\cite{Fis07},
in which an enumeration formula gives
\begin{equation}\ba{l}|\OP(n,n\pl m,[n],[n\mi1]\!\cup\!\{n\pl m\})|\\[2.5mm]
\quad\qquad\qquad\qquad\ds=\;\sum_{i=0}^{n\mi1}|\OP(n\mi1,n,[n\mi1],[n]\!\setminus\!\{n\mi i\})|\:\biggl(\ba{c}m\pl i\\i\ea\biggr)\\[5.5mm]
\quad\qquad\qquad\qquad\ds=\;\frac{1}{(n\mi1)!\,m!}\;\prod_{i=0}^{n\mi2}\!\frac{(3i\pl1)!}{(n\pl i)!}\;\sum_{i=0}^{n\mi1}
\frac{(2n\mi2\mi i)!\,(n\mi1\pl i)!\,(m\pl i)!}{i!^2\,(n\mi1\mi i)!}\,.\ea\end{equation}

A matrix is \emph{horizontally-symmetric} if it is invariant under reflection about the horizontal line through its center,
and \emph{horizontally-and-vertically symmetric} if it is invariant under reflection about both the horizontal and vertical lines through
its center.  It can be seen that for standard alternating sign matrices, such symmetries can only occur for odd-sized matrices.

Considering horizontally-symmetric $(2n\pl1)\t(2n\pl1)$ standard alternating sign matrices $A$, i.e., those for
which $A_{ij}=A_{2n\pl2\mi i,j}$ for each $(i,j)\in[2n\pl1]\t[2n\pl1]$, it can be seen that
the first and last columns always have their $1$'s in the middle row, $A_{i1}=A_{i,2n\pl1}=\d_{i,n\pl1}$, and that
the middle row always consists entirely of alternating $1$'s and $-1$'s, $A_{n\pl1,j}=(-1)^{j\pl1}$.
Proceeding to the corresponding edge matrices, it is then found that along the first column
of vertices $H_{i0}=0$, $V_{i1}=\d_{i\ge n\pl1}$ and $H_{i1}=\d_{i,n\pl1}$, along the last column of vertices
$H_{i,2n}=\d_{i\ne n\pl1}$, $V_{i,2n\pl1}=\d_{i\ge n\pl1}$ and $H_{i,2n\pl1}=1$,
and along the middle row of vertices $V_{nj}=\d_{j\;\mathrm{even}}$, $H_{n\pl1,j}=\d_{j\;\mathrm{odd}}$ and
$V_{n\pl1,j}=\d_{j\;\mathrm{odd}}$.  It follows that each such $A$ is uniquely determined by
the submatrix formed by the first $n$ rows of $A$ with their first and last columns deleted (i.e., $A_{ij}$ with $(i,j)\in[n]\t[2,2n]$),
and that these submatrices comprise $\ASM(n,2n\mi1,[n],\{1,3,\ldots,2n\mi1\})$.
Following conjectures in~\cite{Rob91a,Rob00}, an enumeration formula was proved in~\cite{Kup02} giving
\begin{equation}|\OP(n,2n\mi1,[n],\{1,3,\ldots,2n\mi1\})|\;=\;
\prod_{i=1}^{n}\frac{(6i\mi2)!}{(2n\pl2i)!}\,.\end{equation}

Considering horizontally-and-vertically-symmetric $(2n\pl3)\t(2n\pl3)$ standard alternating sign matrices $A$, i.e., those for
which $A_{ij}=A_{2n\pl4\mi i,j}=A_{i,2n\pl4\mi j}$ for each $(i,j)\in[2n\pl3]\t[2n\pl3]$, it follows,
since both $A$ and $A^t$ are horizontally symmetric, that each such $A$ is uniquely determined by
the submatrix formed by $A_{ij}$ with $(i,j)\in[2,n\pl1]\t[2,n\pl1]$, and that these submatrices comprise
$\ASM(n,n,\{1,3,\ldots,2\lceil\frac{n}{2}\rceil\mi1\},\{1,3,\ldots,2\lceil\frac{n}{2}\rceil\mi1\})$.
Following conjectures in~\cite{Rob91a,Rob00}, an enumeration formula was proved in~\cite{Oka06} giving
\begin{equation}\label{NHVSASM}\ba{l}\ts|\OP(n,n,\{1,3,\ldots,2\lceil\frac{n}{2}\rceil\mi1\},\{1,3,\ldots,2\lceil\frac{n}{2}\rceil\mi1\})|\;=\\[2.7mm]
\hspace{82mm}\ds
\frac{(\lfloor\frac{3n}{2}\rfloor\pl1)!}{3^{\lfloor\frac{n}{2}\rfloor}\:(2n\pl1)!\;\lfloor\frac{n}{2}\rfloor!}\;
\prod_{i=1}^{n}\!\frac{(3i)!}{\;(n\pl i)!}\,.\ea\end{equation}

\subsubsection{Vacancy-Osculation Sets and Matrices}
In this section, the vacancy-osculation set $Z(P)$ associated with each path tuple $P\in\OP(a,b)$ will be
studied further, and it will be shown that, for fixed $a$ and $b$, $P$ is uniquely determined by $Z(P)$.

For any $a,b\in\P$, $(\a,\b)\in\BP(a,b)$ and $l\in\N$, define sets of vacancy-osculation sets as
\begin{equation}\label{VOS}\ba{r@{\:}c@{\:}l@{\qquad}r@{\:}c@{\:}l}\VOS(a,b)&:=&Z(\OP(a,b)),&\VOS(a,b,\a,\b)&:=&Z(\OP(a,b,\a,\b)),\\[2mm]
\VOS(a,b,l)&:=&Z(\OP(a,b,l)),&\VOS(a,b,\a,\b,l)&:=&Z(\OP(a,b,\a,\b,l))\,,\ea\end{equation}
where for each subset $Q$ of $\OP(a,b)$, $Z(Q):=\{Z(P)\mid P\in Q\}$.

It is sometimes convenient to represent or visualize a vacancy-osculation set $S\in\VOS(a,b)$ as
a \emph{vacancy-osculation matrix} $M(S)$. This is an $a\times b$ $\{0,1\}$ matrix defined simply by
\begin{equation}\label{VOSToVOM}M(S)_{ij}:\;=\;\d_{(i,j)\in S}\,,\mbox{ \ for each }(i,j)\in[a]\t[b]\,.\end{equation}
Note that in a vacancy-osculation matrix, the positions of the nonzero entries correspond to the positions
of vertex configurations of types~1 and~2 in~(\ref{6V}), whereas
in an alternating sign matrix they correspond to the positions of configurations of types~5 and~6.

By examining the six possibilities~(\ref{6V}) for each vertex configuration~(\ref{VC}),
it can be seen that if~$S=Z(P)$ is the vacancy-osculation set of~$P\in\OP(a,b)$, and $(H,V)\in\EM(a,b)$ is the edge matrix
pair which corresponds to $P$, then
\begin{equation}\label{HVToVOM}M(S)_{ij}\;=\;\d_{H_{i,j\mi1},V_{ij}}\;=\;\d_{V_{i\mi1,j},H_{ij}}\,.\end{equation}

The vacancy-osculation set and matrix for the running example are
\begin{equation}\label{ExS}S=\{(1,1),(1,2),(1,3),(1,4),(2,1),(2,2),(2,3),(2,5),(3,4),(4,2),(4,3)\}\end{equation}
and
\begin{equation}\label{ExM}M(S)=\left(\ba{cccccc}
1&1&1&1&0&0\\1&1&1&0&1&0\\0&0&0&1&0&0\\0&1&1&0&0&0\ea\right).\end{equation}

\begin{lemma}Each tuple of osculating paths in $\OP(a,b)$ is uniquely determined by its
vacancy-osculation set.\end{lemma}
In other words, (\ref{OPToVOS}) gives an injective, and hence by~(\ref{VOS}) bijective, mapping from
$\OP(a,b)$ to $\VOS(a,b)$.  It immediately follows that
$\OP(a,b,\a,\b)$, $\OP(a,b,l)$ and $\OP(a,b,\a,\b,l)$ are bijectively related to
$\VOS(a,b,\a,\b)$, $\VOS(a,b,l)$ and $\VOS(a,b,$ $\a,\b,l)$ respectively, for any $(\a,\b)\in\BP(a,b)$ and $l\in\N$.

\textit{Proof.} \  The relations of~(\ref{HVToVOM}) can also be written
\begin{equation}\label{MRec}V_{ij}\;=\;\d_{H_{i,j\mi1},M(S)_{ij}}\quad\mbox{and}\quad H_{ij}=\d_{V_{i\mi1,j},M(S)_{ij}}\,,\end{equation}
which can be regarded as recursion relations for $(H,V)$.
Therefore, if $S$ is the vacancy-osculation set of $P\in\OP(a,b)$,
and $(H,V)\in\EM(a,b)$ is the edge matrix pair which corresponds to $P$,
then $(H,V)$, and hence $P$, can be uniquely recovered from $S$ using~(\ref{MRec})
together with the initial conditions $H_{i0}=V_{0j}=0$ for all $i\in[a]$, $j\in[b]$.\hspace{\fill}$\Box$

In fact (\ref{HVToVOM}) or (\ref{MRec}) can also be written
\begin{equation}\label{HVM}H_{i,j\mi1}+V_{ij}\;=\;V_{i,j\mi1}+H_{ij}\;\equiv\;\d_{(i,j)\notin S}\;\equiv\;M(S)_{ij}+1
\pmod2\,,\refstepcounter{equation}\label{VOMToHV}\addtocounter{equation}{-1}\end{equation}
so that, using this with the initial conditions, $(H,V)$ can be expressed explicitly in terms of $S$ as
\rnc{\theequation}{\arabic{equation}\alph{subeq}}\setcounter{subeq}{1}
\begin{equation}\ba{rcll}
H_{ij}&\equiv&\Bigl|\,\Bigl(\{(i\mi k,j\mi k)\mid k\!\in\![0,\min(i,j)\mi1]\}\:\:\cup\\[2mm]
&&\qquad\{(i\mi k\mi1,j\mi k)\mid k\!\in\![0,\min(i\mi1,j)\mi1]\}\Bigr)\!\Bigm\backslash S\:\Bigr|&\pmod2\\[4.5mm]
&\equiv&\ds\d_{i\le j}\;+\sum_{k=0}^{\min(i,j)-1}\!M(S)_{i\mi k,j\mi k}\;\;+\;
\sum_{k=0}^{\min(i\mi1,j)-1}\!\!M(S)_{i\mi k\mi1,j\mi k}&\pmod2\\[6.5mm]
&&\mbox{ \ for each }(i,j)\in[a]\t[0,b]\ea
\end{equation}
\begin{equation}\stepcounter{subeq}\addtocounter{equation}{-1}\ba{rcll}
V_{ij}&\equiv&\Bigl|\,\Bigl(\{(i\mi k,j\mi k)\mid k\!\in\![0,\min(i,j)\mi1]\}\:\:\cup\\[2mm]
&&\qquad\{(i\mi k,j\mi k\mi1)\mid k\!\in\![0,\min(i,j\mi1)\mi1]\}\Bigr)\!\Bigm\backslash S\:\Bigr|&\pmod2\\[4.5mm]
&\equiv&\ds\d_{i\ge j}\;+\sum_{k=0}^{\min(i,j)-1}\!M(S)_{i\mi k,j\mi k}\;\;+\;
\sum_{k=0}^{\min(i,j\mi1)-1}\!\!M(S)_{i\mi k,j\mi k\mi1}&\pmod2\\[6.5mm]
&&\mbox{ \ for each }(i,j)\in[0,a]\t[b]\,,\ea
\end{equation}

\vspace{-2mm}
where each $H_{ij}$ and $V_{ij}$ is taken to be $0$ or $1$.\rnc{\theequation}{\arabic{equation}}

It can be seen that for any $a,b\in\P$, $\emptyset$ and $[a]\t[b]$ are both in $\VOS(a,b)$,
$\emptyset$ corresponding to the path tuple $P\in\OP(a,b,[\min(a,b)],[\min(a,b)])$
in which path $P_i$ passes vertically from $(a,i)$ to $(i,i)$, and then
horizontally from $(i,i)$ to $(i,b)$, for each $i\in[\min(a,b)]$, e.g.,
\[\includegraphics[width=23mm]{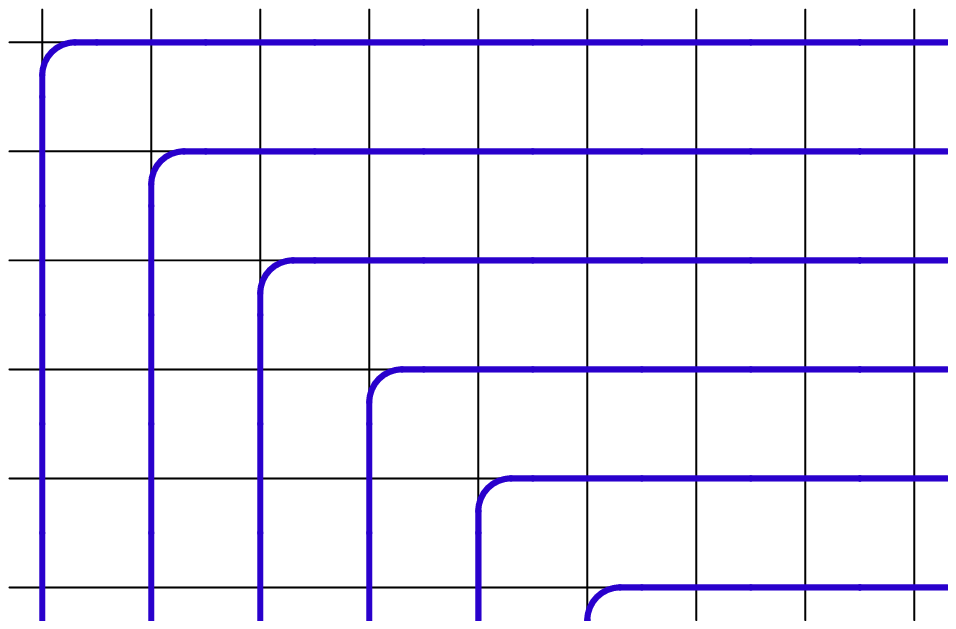}\;\;\raisebox{1mm}{,}\]
and $[a]\t[b]$ corresponding to the single empty path tuple in $\OP(a,b,\emptyset,\emptyset)$, e.g.,
\[\includegraphics[width=23mm]{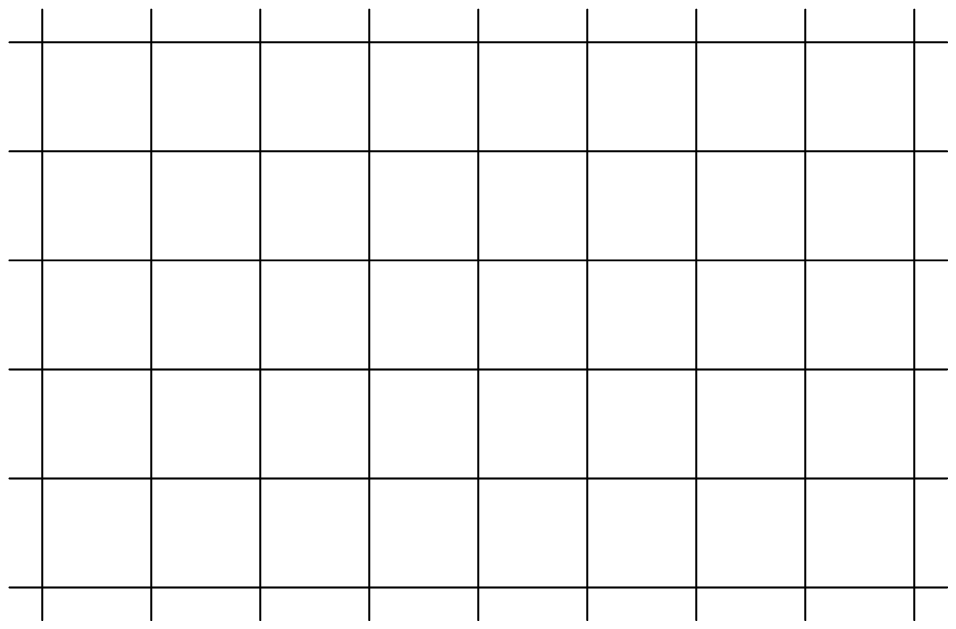}\;\;\raisebox{1mm}{.}\]

Finally, for any $S\in\VOS(a,b)$ corresponding to $P\in\OP(a,b)$,
a vacancy or osculation of $P$ will also be referred to as a vacancy or osculation of $S$,
and the same notation will be used for the sets of vacancies and osculations,
$X(S):=X(P)$ and $N(S):=N(P)$, and for the number of osculations,
$\chi(S):=\chi(P)$.

\subsubsection{Partitions}
In the previous sections, four bijectively-related sets, $\OP(a,b)$,
$\EM(a,b)$, $\ASM(a,b)$ and $\VOS(a,b)$, have been described,
and the aim of forthcoming sections will be to obtain a further, intrinsic
characterization of $\VOS(a,b)$ which facilitates the enumeration of these sets.
The first step in this process will involve a transformation between a boundary point pair $(\a,\b)\in\BP(a,b)$
and a partition, so in this section the relevant notation for partitions is outlined.

A partition $\l=(\l_1,\l_2,\ldots)$ is an infinite sequence of nonnegative integers
in weakly decreasing order, $\l_1\ge\l_2\ge\ldots$, which has only finitely-many nonzero terms.
The nonzero terms are the \emph{parts} of~$\l$, and the number of
parts is the \emph{length} of $\l$, denoted~$\ell(\l)$.  If the sum of parts is $n$, then
$\l$ is said to be a partition of $n$, denoted $|\l|=n$.
The set of all partitions will be denoted as~Par.
When writing a partition as an explicit sequence of terms, some or all of the zero terms will be omitted.
The unique partition of zero will be denoted by~$\emptyset$ and called the empty partition.

For any partition $\l$ define
\begin{equation}Y(\l):=\{(i,j)\in\P^2\mid i\le\ell(\l),\;j\le\l_i\}\,.\end{equation}
The \emph{Young diagram} of $\l$ is then a depiction of $\l$ in which a unit square is centered at each $(i,j)\in Y(\l)$, using matrix-type
labeling of the rows and columns of the lattice.
The \emph{conjugate} of $\l$, denoted~$\l^t$, is the partition whose Young diagram is related to that of $\l$ by
reflection in the main diagonal, $Y(\l^t)=\{(j,i)\mid (i,j)\in Y(\l)\}$.  It follows immediately that
$|\l|=|\l^t|=|Y(\l)|$, $\ell(\l^t)=\l_1$ and $\l^{t\,t}=\l$.

A running example will be the partition $\l=(3,2,2)$, for which $|\l|=7$, $\ell(\l)=3$,
$\l^t=(3,3,1)$, $Y(\l)=\{(1,1),(1,2),(1,3),(2,1),(2,2),(3,1),(3,2)\}$, and the Young diagram is
\setlength{\unitlength}{4mm}
\[\raisebox{-1.27\unitlength}[1.7\unitlength][1.3\unitlength]{
\bpic(3,3)\multiput(0,0)(1,0){3}{\line(0,1){3}}\multiput(0,0)(0,1){2}{\line(1,0){2}}
\multiput(0,2)(0,1){2}{\line(1,0){3}}\put(3,2){\line(0,1){1}}\epic}\,.\]

The \emph{rank} of a partition $\l$ is defined as
\begin{equation}\rho(\l):=|\{i\in[\ell(\l)]\mid\l_i\ge i\}|\,,\end{equation}
and is thus the number of squares on the main diagonal of the Young diagram of $\l$,
$\rho(\l)=|\{(i,j)\in Y(\l)\mid i=j\}|$.

The \emph{Frobenius notation} for a partition $\l$ is
\begin{equation}\label{FN}\ba{rl}\l=(\g_1,\ldots,\g_r\,|\,\d_1,\ldots,\d_r),&\mbox{where }r=\rho(\l),\;\g_i=\l_i\mi i\mbox{ and}\\[2mm]
&\mbox{\quad$\d_i=\l^t_i\mi i$ for each }i\in[r]\,.\ea\end{equation}
Thus, $\g_i$ is the number of squares in row $i$ of the Young diagram to
the right of the main diagonal, and $\d_i$ is the number of squares in column $i$ of the Young diagram below
the main diagonal.  It can be seen that each pair of tuples of nonnegative integers
$(\g_1,\ldots,\g_r)$ and $(\d_1,\ldots,\d_r)$, with $\g_1>\ldots>\g_r$ and
$\d_1>\ldots>\d_r$, corresponds to a unique partition $(\g_1,\ldots,\g_r\,|\,\d_1,\ldots,\d_r)$
with rank~$r$. For the example $(3,2,2)$, the rank is $2$, and the Frobenius notation is $(2,0\,|\,2,1)$.

For any $(i,j)\in\Z^2$, let the \emph{content} of $(i,j)$, or of a unit square centered at $(i,j)$, be $j\mi i$.
Then, for any $d\in\Z$ and any subset $T$ of $\Z^2$, let the $d$-\emph{diagonal} of $T$ be
the set of points of $T$ with content $d$,
\begin{equation}D_d(T):=\{(i,j)\in T\mid j\mi i=d\}\,,\end{equation}
and let the $d$-\emph{rank} of $T$ be the number of points in the $d$-diagonal of $T$,
\begin{equation}\label{Rd}R_d(T)\;:=\;|D_d(T)|\,.\end{equation}
For a partition $\l$, define
\begin{equation}\label{rhod}\rho_d(\l)\;:=\;R_d(Y(\l))\,.\end{equation}
It can be seen that $\rho_0(\l)=\rho(\l)$ and
\begin{equation}\label{Frob}\rho_d(\l)\;=\;\left\{\ba{l}|\{i\in[\rho(\l)]\mid\d_i\ge-d\}|,\quad d\le0\\[2.5mm]
|\{i\in[\rho(\l)]\mid\g_i\ge d\}|,\quad d\ge0\,,\ea\right.\end{equation}
where $(\g_1,\ldots,\g_{\rho(\l)}\,|\,\d_1,\ldots,\d_{\rho(\l)})$ is the Frobenius notation for $\l$.
For the running example,
\begin{equation}\label{RdY322}\rho_d(3,2,2)=\left\{\ba{ll}1,&d\in\{-2,1,2\}\\
2,&d\in\{-1,0\}\\
0,&\mathrm{otherwise\,.}\ea\right.\end{equation}

Partitions $\l$ and $\m$ will be said to \emph{differ by a square},
denoted $\l\!\sim\!\m$, if and only if there exists $i\in\P$ such that
$|\l_k\mi\m_k|=\d_{ki}$ for each $k\in\P$.
Partitions $\l$ and~$\m$ thus differ by a square if and only if there exists $(i,j)\in\P^2$ such that
$Y(\l)$ is the disjoint union of $Y(\m)$ and $\{(i,j)\}$, or $Y(\m)$ is the disjoint union of $Y(\l)$ and $\{(i,j)\}$,
and in such a case the \emph{diagonal difference} between $\l$ and $\m$ is defined as
\begin{equation}\label{DD}\D(\l,\m):=j\mi i\,.\end{equation}
In other words, for $\l\!\sim\!\m$, $\D(\l,\m)$ is the content of the square by which the
Young diagrams of $\l$ and $\m$ differ.
It can be seen that given a partition~$\l$ and a positive integer $i$,
there exists a partition~$\m$ with $\l\!\sim\!\m$ and
$Y(\l)=Y(\m)\cup\{(i,\l_i)\}$ if and only if $\l_i>\l_{i\pl1}$,
and there exists a partition $\m$ with $\l\!\sim\!\m$ and
$Y(\m)=Y(\l)\!\cup\!\{(i,\l_i\pl1)\}$ if and only if $i=1$ or $\l_{i\mi1}>\l_i$.
It follows that for any partition~$\l$, the number of partitions $\m$
with $\l\!\sim\!\m$ is $2\bar{\ell}(\l)\pl1$, where $\bar{\ell}(\l)$ is the number of distinct parts of $\l$.
It can also be seen that for fixed $\l$, each $\m$ with $\l\!\sim\!\m$ is uniquely determined by the diagonal difference~$\D(\l,\m)$.

Define a \emph{change diagonal} of a partition $\l$ to be any integer $d$ for which
there exists a (necessarily unique) partition $\m$ with $\l\!\sim\!\m$ and $\D(\l,\m)=d$.  It can be checked straightforwardly
that
\begin{equation}\label{rhoCond}\ba{l}d\mbox{ is a change diagonal of }\l\mbox{ if and only if }\\[1.5mm]
\qquad\Bigl(\rho_d(\l)-\rho_{d\mi1}(\l),\rho_{d\pl1}(\l)-\rho_d(\l)\Bigr)\;=\;\left\{\ba{ll}
(1,0)\mbox{ \ or \ }(0,1),\quad d<0\\[1mm]
(1,-1)\mbox{ \ or \ }(0,0),\quad d=0\\[1mm]
(0,-1)\mbox{ \ or \ }(-1,0),\quad d>0\,.\ea\right.\ea\end{equation}
Within each of the three cases of $d$ in~(\ref{rhoCond}), the first and second alternatives correspond to the existence
of $\m\!\sim\!\l$ with respectively $|\m|=|\l|\mi1$ and $|\m|=|\l|\pl1$.
For the example $\l=(3,2,2)$, the change diagonals are $-3$, $-1$, $1$, $2$ and $3$, and the corresponding
partitions~$\m\!\sim\!\l$ are
$(3,2,2,1)$, $(3,2,1)$, $(3,3,2)$, $(2,2,2)$ and $(4,2,2)$ respectively.

Finally, for any $a,b\in\P$, define $\Par(a,b)$ to be the
set of all partitions whose Young diagram fits into an $a$ by $b$ rectangle,
\begin{equation}\Par(a,b)\;:=\;\{\l\in\Par\;|\;\ell(\l)\le a\mbox{ and }\l_1\le b\}\,.\end{equation}
By associating a partition in $\Par(a,b)$ with the path formed by the lower and right boundary edges of its Young diagram,
it can be seen that \ru{2.5}$|\Par(a,b)|=\biggl(\ba{c}a\pl b\\a\ea\biggr)$.
For any $\l\in\Par(a,b)$, the $(a,b)$-\emph{complement} of $\l$, denoted $[a]\t[b]\!\setminus\!\l$, is
defined to be the partition
\begin{equation}\label{PC}[a]\t[b]\!\setminus\!\l\;:=\;(b\mi\l_a,b\mi\l_{a\mi1},\ldots,b\mi\l_1)\,,\end{equation}
where, as always, $\l_i$ is $0$ for $i>\ell(\l)$.  It can be seen that
\begin{equation}Y([a]\t[b]\!\setminus\!\l)\;=\;\{(a\pl1\mi i,b\pl1\mi j)\,|\,(i,j)\in[a]\t[b]\!\setminus\!Y(\l)\}\end{equation}
and
\begin{equation}\label{comp}\rho_d([a]\t[b]\!\setminus\!\l)\;=\;R_d([a]\t[b])-\rho_{b\mi a\mi d}(\l)\,.\end{equation}
An example of a complement of $\l=(3,2,2)$ is $[4]\t[6]\!\setminus\!\l=(6,4,4,3)$.

\subsubsection{Correspondence Between Boundary Point Pairs and Partitions}
Returning to the case of osculating paths $\OP(a,b)$,
it will be shown in this section that each boundary point pair
in $\BP(a,b)$ corresponds naturally to a partition in $\Par(a,b)$,
and that the correspondence leads to an important property for
osculating paths.

For any $(\a,\b)=(\{\a_1,\ldots,\a_r\},\{\b_1,\ldots,\b_r\})\in\BP(a,b)$,
the corresponding partition is defined to be
\begin{equation}\label{ToP}
\l_{a,b,\a,\b}\;:=\;
[a]\t[b]\!\setminus\!(b\mi\b_1,\ldots,b\mi\b_r\,|\,a\mi\a_1,\ldots,a\mi\a_r)\,,\end{equation}
where the definitions~(\ref{FN}) of Frobenius notation and~(\ref{PC}) of the $(a,b)$-complement for partitions are being used,
and $\a_1<\ldots<\a_r$, $\b_1<\ldots<\b_r$.

It can be checked straightforwardly that for fixed $a$ and $b$ this mapping is a bijection between
$\BP(a,b)$ and $\Par(a,b)$, and that for $\l\in\Par(a,b)$ the inverse mapping is
given by
\begin{equation}\label{BPInv}\ba{l}(\a_{a,b,\l},\b_{a,b,\l})\;=\;(\{a\mi\d_1,\ldots,a\mi\d_r\},\{b\mi\g_1,\ldots,b\mi\g_r\})\,,\\[2.6mm]
\mbox{\qquad where\quad}[a]\t[b]\!\setminus\!\l\;=\;(\g_1,\ldots,\g_r\,|\,\d_1,\ldots,\d_r)\,.\ea\end{equation}
It can also be seen that
\begin{equation}\label{Parmag}\l_{b,a,\b,\a}\;=\;(\l_{a,b,\a,\b})^t\,,\qquad
|\l_{a,b,\a,\b}|\;=\;ab-(a\pl b\pl1)r+\sum_{i=1}^r(\a_i\pl\b_i)\,,\end{equation}
and that some special cases of~(\ref{ToP}) are
\begin{equation}\label{ParSpec}\ba{c}\l_{a,b,[\min(a,b)],[\min(a,b)]}\;=\;\emptyset\,,\qquad
\l_{a,b,\emptyset,\emptyset}\;=\;(\,\underbrace{b,\ldots,b}_{a}\,)\\[7mm]
\l_{a,b,[a],\b}=(\b_a\mi a,\ldots,\b_2\mi2,\b_1\mi1)\,,\qquad
\l_{a,b,\a,[b]}=(\a_b\mi b,\ldots,\a_2\mi2,\a_1\mi1)\\[5mm]
\l_{a,a,\a,\b}\;=\;(a\mi\bar{\a}_1,\ldots,a\mi\bar{\a}_{\bar{r}}\,|\,
a\mi\bar{\b}_1,\ldots,a\mi\bar{\b}_{\bar{r}})\,,\\[2mm]\qquad\qquad\qquad\qquad\qquad\quad\mbox{where }
\{\bar{\a}_1,\ldots,\bar{\a}_{\bar{r}}\}\!=\![a]\!\setminus\!\a\,,\;\;
\{\bar{\b}_1,\ldots,\bar{\b}_{\bar{r}}\}\!=\![a]\!\setminus\!\b\,.\ea\end{equation}

For $a=4$, $b=6$ and $(\a,\b)=(\{1,2,3\},\{1,4,5\})$, as used for the running example of a tuple of osculating paths,
the corresponding partition is $\:\l_{4,6,\{1,2,3\},\{1,4,5\}}=$ $[4]\t[6]\setminus\!(5,2,1\,|\,3,2,1)=
[4]\t[6]\setminus\!(6,4,4,3)=(3,2,2)$, which was used as the running example of a partition.

For the five special cases of $a$, $b$, $\a$ and $\b$
used in Section~4 to give previously-studied cases of alternating sign
matrices, the corresponding partitions are given in Table~\ref{ASMPar}.
\newpage
\begin{table}[h]\centering$
\ba{c|c}a,\;b,\;\a,\;\b&\l_{a,b,\a,\b}\\[2mm]\hline
\rule{0mm}{7mm}n,\;n,\;[n],\;[n]&\emptyset\\[2mm]
n,\;n\pl1,\;[n],\;[n\pl1]\!\setminus\!\!\{n\pl1\mi m\}&(m)^t\\[2mm]
n,\;n\pl m,\;[n],\;[n\mi1]\!\cup\!\{n\pl m\}&(m)\\[2mm]
n,\;2n\mi1,\;[n],\;\{1,3,\ldots,2n\mi1\}\;&\;\:(n\mi1,n\mi2,\ldots,1)\\[2mm]
\;n,\;n,\;\{1,3,\ldots,2\lceil\frac{n}{2}\rceil\mi1\},\;\{1,3,\ldots,2\lceil\frac{n}{2}\rceil\mi1\}\;&\;\:(n\mi1,n\mi2,\ldots,1)\ea$
\caption{\protect\rule{0ex}{2.5ex}Partitions for alternating sign matrix cases.\label{ASMPar}}
\end{table}

\begin{lemma}For any tuple of osculating paths $P\in\OP(a,b,\a,\b)$ and integer $d$,
\begin{equation}\label{lem2}R_d(N(P))-R_d(X(P))\;=\;\rho_d(\l_{a,b,\a,\b})\,.\end{equation}\end{lemma}
Here, $N(P)$ and $X(P)$ are the sets of vacancies and osculations for~$P$ as defined in Section~2, and $R_d$ and $\rho_d$ are $d$-ranks as
defined in~(\ref{Rd}) and~(\ref{rhod}).  In other words, the difference between the numbers of
vacancies and osculations on any diagonal of the lattice is independent of the path tuple $P$, and equal to the
number of squares on the corresponding diagonal of the Young diagram of $\l_{a,b,\a,\b}$.

Lemma~2 implies that for a path tuple $P\in\OP(a,b,\a,\b)$,
the partition $\l_{a,b,\a,\b}$ can be obtained geometrically
by first marking on the lattice the vacancies $N(P)$,
and marking differently the osculations $X(P)$, then along each diagonal $D_d([a]\t[b])$
removing all of the osculation marks together with an equal number of vacancy marks, and finally along each diagonal
moving, if necessary, the remaining vacancy marks diagonally upwards and leftwards so that they occupy adjacent positions
starting from the upper or left boundary of the lattice.  The occupied points then correspond to $Y(\l_{a,b,\a,\b})$.

\textit{Proof.} \ As usual, let $(\a,\b)=(\{\a_1,\ldots,\a_r\},\{\b_1,\ldots,\b_r\})$
with $\a_1<\ldots<\a_r$ and $\b_1<\ldots<\b_r$.
Now let $J_d(a,b,\a,\b)$ be the number of paths of $P$
which pass through the $d$-diagonal $D_d([a]\t[b])$.
Since path $P_i$ passes continuously by upward or rightward steps from $(a,\b_i)$ to $(\a_i,b)$,
it can be seen by placing a line through the points of $D_d([a]\t[b])$
that $J_d(a,b,\a,\b)$ is the number of
start points $(a,\b_i)$ on or to the left of $(a,a\pl d)$ for $d\le b\mi a$,
and the number of end points $(\a_i,b)$ on or above $(b\mi d,b)$ for $d\ge b\mi a$.
Therefore, as suggested by
the notation, $J_d(a,b,\a,\b)$ is independent of $P\in\OP(a,b,\a,\b)$,
and given explicitly as
\begin{equation}\label{lem21}J_d(a,b,\a,\b)\;=\;
\left\{\ba{ll}|\{i\in[r]\,|\,\b_i\le a\pl d\}|,&d\le b\mi a\\[2.5mm]
|\{i\in[r]\,|\,\a_i\le b\mi d\}|,&d\ge b\mi a\,.\ea\right.\end{equation}
Defining $I(P)$ to be the number of points of $[a]\t[b]$ through which a single path of~$P$ passes
(i.e., points with vertex configurations of types 3--6 in~(\ref{6V})), so that
$[a]\t[b]$ is the disjoint union of $I(P)$, $X(P)$ and $N(P)$,
and using the fact that an osculation is a point through which two paths pass,
it follows that
\begin{equation}\ba{r@{\;}c@{\;}l}\label{lem22}R_d([a]\t[b])&=&R_d(I(P))+R_d(X(P))+R_d(N(P))\\[2.4mm]
J_d(a,b,\a,\b)&=&R_d(I(P))+2R_d(X(P))\,.\ea\end{equation}
Eliminating $J_d(a,b,\a,\b)$ and $R_d(I(P))$ from~(\ref{lem21}) and~(\ref{lem22}) now gives
\begin{equation}\label{lem23}R_d(N(P))\mi R_d(X(P))\;=\;R_d([a]\t[b])-
\left\{\ba{l@{\;\:}l}|\{i\in[r]\,|\,\b_i\le a\pl d\}|,&d\le b\mi a\\[2.5mm]
|\{i\in[r]\,|\,\a_i\le b\mi d\}|,&d\ge b\mi a\,.\ea\right.\end{equation}
Finally, using~(\ref{comp}) with $\l=(b\mi\b_1,\ldots,b\mi\b_r\,|\,a\mi\a_1,\ldots,a\mi\a_r)$,
and then using~(\ref{Frob}) on the second term from the RHS of~(\ref{comp}), it is found that the RHS
of~(\ref{lem23}) is $\rho_d([a]\t[b]\!\setminus\!(b\mi\b_1,\ldots,b\mi\b_r\,|\,a\mi\a_1,\ldots,a\mi\a_r))$ as required.\hspace{\fill}$\Box$

It follows from Lemma~2, by summing~(\ref{lem2}) over all $d\in[1\mi a,b\mi1]$, that
\begin{equation}\label{Pmag}|N(P)|\;=\;|\l_{a,b,\a,\b}|+\chi(P)\,,\quad\mbox{for any }P\in\OP(a,b,\a,\b)\,,\end{equation}
and therefore that
\begin{equation}\label{Smag}|S|\;=\;|\l_{a,b,\a,\b}|+2\chi(S)\,,\quad\mbox{for any }S\in\VOS(a,b,\a,\b)\,.\end{equation}
Note that in (\ref{Pmag}) or (\ref{Smag}), $|\l_{a,b,\a,\b}|$ can be obtained directly from $a$, $b$, $\a$ and $\b$
using~(\ref{Parmag}).

Lemma~2 can be verified for the running example by observing that
\[R_d(N(P))=\left\{\ba{ll}1,&d\in\{-2,2,3\}\\
2,&d\in\{-1,0,1\}\\
0,&\mathrm{otherwise}\ea\right.\qquad\mbox{and}\qquad
R_d(X(P))=\left\{\ba{ll}1,&d\in\{1,3\}\\
0,&\mathrm{otherwise\,,}\ea\right.\]
and that $\rho_d(\l_{4,6,\{1,2,3\},\{1,4,5\}})$ is given by~(\ref{RdY322}).

\subsubsection{Vacancy-Osculation Sets which Differ by a Vacancy-Osculation}
In this section, the relationship between two vacancy-osculation sets which differ by a single element
will be examined, and it will be seen that the corresponding partitions then differ by a square.

Vacancy-osculation sets $S,S'\in\VOS(a,b)$
will be said to \emph{differ by the vacancy-osculation} $(i,j)$ if and only if
the sets differ by the single element $(i,j)$, i.e.,
$S$ is the disjoint union of $S'$ and $\{(i,j)\}$, or $S'$ is the disjoint union of $S$ and $\{(i,j)\}$.
Furthermore, in such a case the \emph{diagonal difference} between $S$ and $S'$ will be defined as
\begin{equation}\label{VOSDelta}\D(S,S'):=j\mi i\,.\end{equation}
The actual position of $(i,j)$ will sometimes be unimportant, so
$S,S'\in\VOS(a,b)$ will be said simply to \emph{differ by a vacancy-osculation},
denoted $S\sim S'$, if there exists some $(i,j)\in[a]\t[b]$ such that $S$ and $S'$ differ by $(i,j)$.
It can be seen that $\sim$ and $\D$ for vacancy-osculation sets have been defined analogously to~$\sim$
and~$\D$ for partitions, and the reason for this will become apparent in Lemma~6.

In terms of the associated vacancy-osculation matrices, $S,S'\in\VOS(a,b)$
differ by a vacancy-osculation if and only if the entries of
the matrices match at all positions except a single
$(i,j)\in[a]\t[b]$, $|M(S)_{kl}\mi M(S')_{kl}|=\d_{(i,j),(k,l)}$ for each $(k,l)\in[a]\t[b]$.

Note that if vacancy-osculation sets $S$ and $S'$ differ by a vacancy-osculation,
it does not necessarily follow that $N(S)$ and $N(S')$, or $X(S)$ and $X(S')$,
individually differ by a single element, it being possible that $N(S)\cap X(S')$ or
$X(S)\cap N(S')$ are nonempty.

Define a \emph{deletion point} of $S\in\VOS(a,b)$ to be any $(i,j)\in S$ for which
there exists $S'\in\VOS(a,b)$ with $S\sim S'$ and $S'=S\setminus\!\{(i,j)\}$,
and define an \emph{addition point} of $S\in\VOS(a,b)$ to be any $(i,j)\in[a]\t[b]\!\setminus\!S$ for which
there exists $S'\in\VOS(a,b)$ with $S\sim S'$ and $S'=S\cup\{(i,j)\}$.
Define a \emph{change point} of $S$ to be any deletion point or addition point of $S$.

\begin{figure}[h]\centering\bpic(150,16)\scriptsize
\put(15,8){\p{}{\left(\ba{c@{\:\:}c@{\:\:}c@{\:\:}c@{\:\:}c@{\:\:}c}1&1&1&1&0&0\\1&1&1&0&1&0\\0&0&0&1&0&0\\0&1&1&0&0&0\ea\right)}}
\put(45,8){\p{}{\left(\ba{c@{\:\:}c@{\:\:}c@{\:\:}c@{\:\:}c@{\:\:}c}1&1&1&1&0&0\\1&1&1&0&1&0\\0&0&0&1&0&0\\0&1&\mathbf{0}&0&0&0\ea\right)}}
\put(75,8){\p{}{\left(\ba{c@{\:\:}c@{\:\:}c@{\:\:}c@{\:\:}c@{\:\:}c}1&1&1&1&0&0\\1&1&\mathbf{0}&0&1&0\\0&0&0&1&0&0\\0&1&1&0&0&0\ea\right)}}
\put(105,8){\p{}{\left(\ba{c@{\:\:}c@{\:\:}c@{\:\:}c@{\:\:}c@{\:\:}c}1&1&1&1&0&0\\1&1&1&0&1&0\\0&0&0&\mathbf{0}&0&0\\0&1&1&0&0&0\ea\right)}}
\put(135,8){\p{}{\left(\ba{c@{\:\:}c@{\:\:}c@{\:\:}c@{\:\:}c@{\:\:}c}1&1&1&1&0&0\\1&1&1&0&1&0\\0&0&0&1&0&0\\0&1&1&0&\mathbf{1}&0\ea\right)}}\epic\\
\includegraphics[width=150mm]{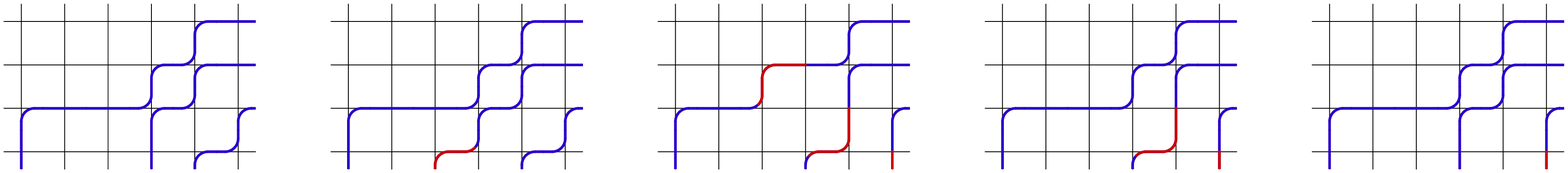}\\
\bpic(150,8)
\put(15,5){\makebox(0,0)[]{{\scriptsize original}}}\put(45,5){\pp{}{(4,3)}}\put(75,5){\pp{}{(2,3)}}\put(105,5){\pp{}{(3,4)}}\put(135,5){\pp{}{(4,5)}}
\epic\\
\bpic(150,17)\scriptsize
\put(15,8){\p{}{\left(\ba{c@{\:\:}c@{\:\:}c@{\:\:}c@{\:\:}c@{\:\:}c}1&1&1&1&0&0\\1&1&1&0&1&0\\0&0&0&1&\mathbf{1}&0\\0&1&1&0&0&0\ea\right)}}
\put(45,8){\p{}{\left(\ba{c@{\:\:}c@{\:\:}c@{\:\:}c@{\:\:}c@{\:\:}c}1&1&1&1&0&0\\1&1&1&0&1&0\\0&0&0&1&0&0\\0&1&1&0&0&\mathbf{1}\ea\right)}}
\put(75,8){\p{}{\left(\ba{c@{\:\:}c@{\:\:}c@{\:\:}c@{\:\:}c@{\:\:}c}1&1&1&\mathbf{0}&0&0\\1&1&1&0&1&0\\0&0&0&1&0&0\\0&1&1&0&0&0\ea\right)}}
\put(105,8){\p{}{\left(\ba{c@{\:\:}c@{\:\:}c@{\:\:}c@{\:\:}c@{\:\:}c}1&1&1&1&0&0\\1&1&1&0&\mathbf{0}&0\\0&0&0&1&0&0\\0&1&1&0&0&0\ea\right)}}
\put(135,8){\p{}{\left(\ba{c@{\:\:}c@{\:\:}c@{\:\:}c@{\:\:}c@{\:\:}c}1&1&1&1&0&0\\1&1&1&0&1&0\\0&0&0&1&0&\mathbf{1}\\0&1&1&0&0&0\ea\right)}}\epic\\
\includegraphics[width=150mm]{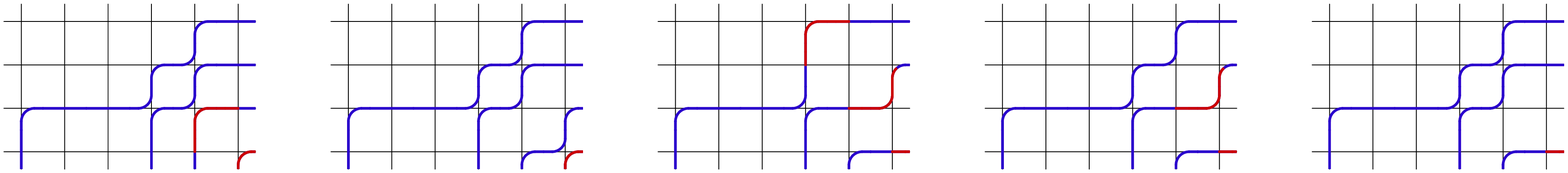}\\
\bpic(150,0)\put(15,2){\pp{}{(3,5)}}\put(45,2){\pp{}{(4,6)}}\put(75,2){\pp{}{(1,4)}}\put(105,2){\pp{}{(2,5)}}\put(135,2){\pp{}{(3,6)}}\epic
\caption{Change points for the running example.\label{ExCP}}
\end{figure}

For the running example, with $S$ given by~(\ref{ExS}), it can be seen in Figure~\ref{ExCP}
and Table~\ref{ExPar} that
$S$ has deletion points $(4,3)$, $(2,3)$, $(3,4)$, $(1,4)$ and $(2,5)$, and addition points $(4,5)$, $(3,5)$, $(4,6)$ and $(3,6)$.
In Figure~\ref{ExCP}, the vacancy-osculation matrix $M(S)$ and diagram of the original path tuple,
as already given in~(\ref{ExM}) and Figure~\ref{Ex}, are shown first, followed by the
vacancy-osculation matrix $M(S')$ and diagram of the corresponding path tuple for each change point $(i,j)$,
where $S'$ is the vacancy-osculation set which differs from $S$ by $(i,j)$.
In each $M(S')$, the $(i,j)$ entry is highlighted, and in each path diagram the sections of paths which differ
from the original are shown in red.
In Table~\ref{ExPar}, the change points are listed together with the
diagonal differences $\D(S,S')=j\mi i$, the $\a'$ and $\b'$ for which $\VOS(4,6,\a',\b')$ contains $S'$, and the partitions $\l_{4,6,\a',\b'}$.
\setlength{\unitlength}{1mm}
\font\tenmsb=msbm10 scaled \magstep1
\font\sevenmsb=msbm7 scaled \magstep1
\font\fivemsb=msbm5 scaled \magstep1
\newfam\msbfam
\textfont\msbfam=\tenmsb
\scriptfont\msbfam=\sevenmsb
\scriptscriptfont\msbfam=\fivemsb

\vspace{4mm}
\begin{table}[h]\centering$\ba{c|c|c|c|c}(i,j)&\;j\mi i\;&\a'&\b'&\;\l_{4,6,\a',\b'}\;\\[2mm]\hline
\rule{0mm}{7mm}
(4,3)&-1&\{1,2,3\}&\{1,3,5\}&(3,2,1)\\[2mm]
\;(2,3),\;(3,4),\;(4,5)&1&\{1,2,3\}&\{1,4,6\}&(3,3,2)\\[2mm]
(3,5),\;(4,6)&2&\{1,2,3,4\}&\{1,4,5,6\}&(2,2,2)\\[2mm]
(1,4),\;(2,5),\;(3,6)&3&\{1,2,4\}&\{1,4,5\}&(4,2,2)\ea$
\caption{\protect\rule{0ex}{2.5ex}Partitions for the change points of the running example.\label{ExPar}}
\end{table}

\begin{figure}[h]
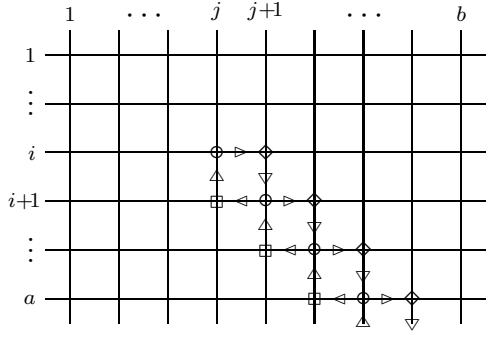
\centering\bpic(9.6,7)
\multiput(0.5,0.5)(0,1){6}{\line(1,0){9}}\multiput(1,0)(1,0){9}{\line(0,1){6}}
\put(0.3,5.5){\pp{r}{1}}\put(0.3,3.5){\pp{r}{i}}\put(0.4,2.5){\pp{r}{i+\!1}}\put(0.3,0.5){\pp{r}{a}}
\put(1,6.2){\pp{b}{1}}\put(4,6.2){\pp{b}{j}}\put(5,6.2){\pp{b}{j+\!1}}\put(9,6.2){\pp{b}{b}}
\multiput(4.5,3.5)(1,-1){4}{\ppp{}{\rhd}}\multiput(5,3)(1,-1){4}{\ppp{}{\bigtriangledown}}
\multiput(4,3)(1,-1){4}{\ppp{}{\bigtriangleup}}\multiput(4.5,2.5)(1,-1){3}{\ppp{}{\lhd}}
\multiput(4,3.5)(1,-1){4}{\p{}{\circ}}\multiput(5,3.5)(1,-1){4}{\p{}{\diamond}}
\multiput(4,2.45)(1,-1){3}{\pp{}{\Box}}
\multiput(2.51,6.3)(4.5,0){2}{\p{}{\cdots}}\multiput(0.2,1.67)(0,3){2}{\p{}{\vdots}}
\epic
\caption{Diagram for Lemmas 3--7.\label{Bfig}}
\end{figure}

A series of lemmas on the properties
of change points and of vacancy-osculation sets which differ by a vacancy-osculation
will now be presented.
In these lemmas, $B(a,b,i,j)$ for any $(i,j)\in\P^2$ will denote the set of points of the diagonal $D_{j\mi i}([a]\t[b])$
weakly below and to the right of $(i,j)$,
\begin{equation}\label{B}B(a,b,i,j):=\Bigl\{(i\pl k,j\pl k)\Bigm|k\in[0,\min(a\mi i,b\mi j)]\Bigr\}\,.\end{equation}
This set is depicted using symbols $\circ$ in Figure~\ref{Bfig}.

\setlength{\unitlength}{6.5mm}

\begin{lemma}Consider a vacancy-osculation set $S\in\mathrm{VOS}(a,b)$ and a point $(i,j)\in[a]\t[b]$.\vspace{-4mm}
\begin{itemize}
\item If $(i,j)\in S$, then $(i,j)$ is a deletion point of $S$ if and only if\\
$S\cap(B(a,b,i,j\pl1)\cup B(a,b,i\pl1,j))=\emptyset$.\vspace{-1mm}
\item If $(i,j)\notin S$, then $(i,j)$ is an addition point of $S$ if and only if\\
$S\cap(B(a,b,i,j\pl1)\cup B(a,b,i\pl1,j))=\emptyset$ and
$H_{ij}=V_{ij}$, where\\ $(H,V)\in\EM(a,b)$ is the edge matrix pair which corresponds to $S$.\end{itemize}
\end{lemma}
The sets $B(a,b,i,j\pl1)$ and $B(a,b,i,j\pl1)$ are depicted in Figure~\ref{Bfig} using the symbols~$\diamond$ and~$\ss\Box$ respectively.
Note that the following conditions are all equivalent.
\begin{equation}\label{equiv1}
\ba{l@{\;\;}l}
\bullet&S\cap(B(a,b,i,j\pl1)\cup B(a,b,i\pl1,j))=\emptyset.\\[1.5mm]
\bullet&\mbox{The vertex configuration at each point of }B(a,b,i,j\pl1)\,\cup\\[1mm]
&B(a,b,i\pl1,j)\mbox{ is of type 3, 4, 5 or 6.}\\[1.5mm]
\bullet&\mbox{A single path of $P$ passes through each point of }B(a,b,i,j\pl1)\,\cup\\
&B(a,b,i\pl1,j),\mbox{ where }P\in\OP(a,b)\mbox{ corresponds to }S.\\[1.5mm]
\bullet&H_{kl}\ne V_{k,l\pl1},\mbox{ for each }(k,l)\in B(a,b,i,j)\!\setminus\!\{(b\pl i\mi j,b)\},\mbox{ and }\\[1mm]
&V_{kl}\ne H_{k\pl1,l},\mbox{ for each }(k,l)\in B(a,b,i,j)\!\setminus\!\{(a,a\pl j\mi i)\}.
\ea
\end{equation}
Also, for any $(k,l)\in[a]\t[b]$, the following are equivalent.
\begin{equation}\label{equiv2}
\ba{l}
\bullet\;\;H_{kl}=V_{kl}.\\[1.5mm]
\bullet\;\;H_{k,l\mi1}=V_{k\mi1,l}.\\[1.5mm]
\bullet\;\;\mbox{The vertex configuration at }(k,l)\mbox{ is of type 1, 2, 5 or 6.}
\ea
\end{equation}
It can be seen using Lemma~3 that the points given in
Figure~\ref{ExCP} and Table~\ref{ExPar} comprise all of the change
points for the vacancy-osculation set~(\ref{ExS}).

\textit{Proof.} \ Let $(H,V)\in\EM(a,b)$ correspond to $S$, and let $S'$ be the set which differs from
$S$ by $(i,j)$,
\[S'\;=\;\left\{\ba{ll}S\setminus\!\{(i,j)\}\,,&(i,j)\in S\\[2mm]
S\cup\{(i,j)\}\,,&(i,j)\notin S\,.\ea\right.\]
It thus needs to be checked whether or not $S'$ is a vacancy-osculation set.
First, define $H'$ and $V'$ by~(\ref{VOMToHV}), using $S'$.
It then follows directly from~(\ref{VOMToHV}) that $(H,V)$ and $(H',V')$ are related by~(\ref{HVD}), which appears in the forthcoming
Lemma~4. The positions of the edges corresponding to
$H_{kl}$ with $(k,l)\in B(a,b,i,j)$, $H_{kl}$ with $(k,l)\in B(a,b,i\pl1,j)$, $V_{kl}$ with $(k,l)\in B(a,b,i,j)$,
and $V_{kl}$ with $(k,l)\in B(a,b,i,j\pl1)$ are shown in Figure~\ref{Bfig} using the symbols~$\rhd$, $\lhd$, $\ss\bigtriangleup$
and $\ss\bigtriangledown$ respectively.  It now follows that $S'\in\VOS(a,b)$ if and only if
$(H',V')\in\EM(a,b)$. The first four conditions on the RHS of~(\ref{EM2}),
which define $\EM(a,b)$, are automatically satisfied due to the definition of $(H',V')$ via~(\ref{VOMToHV}),
so it is only the last condition, which now reads
\begin{equation}\label{ACDD}H'_{k,l\mi1}+V'_{kl}=V'_{k\mi1,l}+H'_{kl}\,,\end{equation}
for each $(k,l)\in[a]\t[b]$, that needs to be checked here.
Since $(H,V)\in\EM(a,b)$,
\begin{equation}\label{ACD}H_{k,l\mi1}+V_{kl}=V_{k\mi1,l}+H_{kl}\end{equation}
is satisfied for all $(k,l)\in[a]\t[b]$.
The condition (\ref{ACDD}) will now be examined for five separate cases:\vspace{-4mm}
\begin{enumerate}
\item \ $(k,l)\in B(a,b,i\pl1,j\pl1)$
\item \ $(k,l)\in B(a,b,i,j\pl1)$
\item \ $(k,l)\in B(a,b,i\pl1,j)$
\item \ $(k,l)=(i,j)$
\item \ All other $(k,l)\in[a]\t[b]$\,.
\end{enumerate}\vspace{-3mm}
In Case (1), (\ref{HVD}) gives $H'_{k,l\mi1}=1\mi H_{k,l\mi1}$,
$V'_{kl}=1\mi V_{kl}$, $V'_{k\mi1,l}=1\mi V_{k\mi1,l}$
and $H'_{kl}=1\mi H_{kl}$.  It thus follows from~(\ref{ACD}) that
in this case~(\ref{ACDD}) is satisfied for all $S\in\VOS(a,b)$.

In Case (2), (\ref{HVD}) gives $H'_{k,l\mi1}=1\mi H_{k,l\mi1}$, $V'_{kl}=1\mi V_{kl}$, $V'_{k\mi1,l}=V_{k\mi1,l}$,
and $H'_{kl}=H_{kl}$.  By examining the six possibilities for $(H_{k,l\mi1},V_{kl},V_{k\mi1,l},H_{kl})$ in~(\ref{ACD}),
it is found that in this case~(\ref{ACDD}) is satisfied if and only if
$(H_{k,l\mi1},V_{kl},V_{k\mi1,l},H_{kl})\ne(1,1,1,1)$ or $(0,0,0,0)$,
or equivalently $(k,l)\notin S$. Thus, Case~(2) leads to the condition
$S\cap B(a,b,i,j\pl1)=\emptyset$.

Case (3) is analogous to Case (2), and leads to the condition $S\cap B(a,b,i\pl1,j)=\emptyset$.

In Case (4), (\ref{HVD}) gives $H'_{i,j\mi1}=H_{i,j\mi1}$, $V'_{ij}=1\mi V_{ij}$, $V'_{i\mi1,j}=V_{i\mi1,j}$ and
$H'_{ij}=1\mi H_{ij}$. It follows that in this case~(\ref{ACDD}) is satisfied if and only if
$(H_{i\!,j\mi1},V_{i\!j},V_{i\mi1,j},H_{i\!j})\ne(1,0,0,1)$ or $(0,1,1,0)$
(i.e., the vertex configuration at $(i,j)$ is not of type~3 or~4 in~(\ref{6V})),
which leads to the condition $H_{ij}=V_{ij}$ if $(i,j)\notin S$.

In Case (5), (\ref{HVD}) gives $H'_{k,l\mi1}=H_{k,l\mi1}$,
$V'_{kl}=V_{kl}$, $V'_{k\mi1,l}=V_{k\mi1,l}$
and $H'_{kl}=H_{kl}$,  so it follows from~(\ref{ACD})
that~(\ref{ACDD}) is satisfied for all $S\in\VOS(a,b)$.

Overall then, the set $S'$ is in $\VOS(a,b)$ if and only if the conditions arising from Cases~2 to~4
are satisfied, which gives the conclusions of Lemma~3.\hspace{\fill}$\Box$

\begin{lemma}Consider vacancy-osculation sets $S,S'\in\mathrm{VOS}(a,b)$, with
corresponding edge matrix pairs $(H,V),(H',V')\in\EM(a,b)$, and a point $(i,j)\in[a]\t[b]$.
Then $S$ and $S'$ differ by the vacancy-osculation $(i,j)$ if and only if
\begin{equation}\label{HVD}\ba{r@{\,\;}c@{\,\;}ll}|H_{kl}\mi H'_{kl}|&=&\d_{(k,l)\in B(a,b,i,j)\cup B(a,b,i\pl1,j)},&
\mbox{for each }(k,l)\in[a]\t[0,b],\\[3mm]
\mbox{and \ \ }|V_{kl}\mi V'_{kl}|&=&\d_{(k,l)\in B(a,b,i,j)\cup B(a,b,i,j\pl1)},&
\mbox{for each }(k,l)\in[0,a]\t[b].\ea\end{equation}
\end{lemma}
The positions of the edges corresponding to
$H_{kl}$ or $H'_{kl}$ with $(k,l)\in B(a,b,i,j)$, $H_{kl}$ or $H'_{kl}$ with $(k,l)\in B(a,b,i\pl1,j)$, $V_{kl}$
or $V'_{kl}$ with $(k,l)\in B(a,b,i,j)$,
and $V_{kl}$ or~$V'_{kl}$ with $(k,l)\in B(a,b,i,j\pl1)$
are shown in Figure~\ref{Bfig} using the symbols~$\rhd$, $\lhd$, $\ss\bigtriangleup$
and $\ss\bigtriangledown$ respectively.

\textit{Proof.} \ If $S$ and $S'$ differ by the vacancy-osculation $(i,j)$, then~(\ref{HVD})
follows directly from~(\ref{VOMToHV}).
Conversely, if~(\ref{HVD}) is satisfied, then, by considering the same five cases for $(k,l)\in[a]\t[b]$ as in the
proof of Lemma~3, it is found, using arguments similar to those used in that proof, that in Case (4) $(k,l)=(i,j)$ is an element of
exactly one of $S$ or $S'$, while in each of the other four cases $(k,l)$ is an element of both or
neither of $S$ and $S'$, so that $S$ and $S'$ differ by $(i,j)$.\hspace{\fill}$\Box$

\begin{lemma}Consider a vacancy-osculation set $S\in\VOS(a,b)$, with corresponding edge matrix pair
$(H,V)\in\EM(a,b)$, and a point $(i,j)\in[a]\t[b]$. Then the following are all equivalent.\vspace{-4mm}
\begin{itemize}
\item $(i,j)$ is a change point of $S$.\vspace{-1mm}
\item $S\cap(B(a,b,i,j\pl1)\cup B(a,b,i\pl1,j))=\emptyset$ and there exists $(k,l)\in B(a,b,i,j)$ such that
$H_{kl}=V_{kl}$.\vspace{-1mm}
\item $S\cap(B(a,b,i,j\pl1)\cup B(a,b,i\pl1,j))=\emptyset$ and $H_{kl}=V_{kl}$ for each $(k,l)\in B(a,b,i,j)$.
\item
\setlength{\unitlength}{9mm}
\begin{equation}\label{lem51}\ba{l}
\bpic(2.5,1.3)\multiput(0,0.5)(0,1){2}{\line(1,0){2}}\multiput(0.5,0)(1,0){2}{\line(0,1){2}}
\multiput(0.5,1)(1,0){2}{\ppp{}{\bullet}}\multiput(1,0.5)(0,1){2}{\ppp{}{\bullet}}
\put(1.05,0.4){\ppp{t}{H_{k\pl1\!,l}}}\put(1,1.61){\ppp{b}{H_{kl}}}
\put(0.43,1){\ppp{r}{V_{\!kl}}}\put(1.6,1){\ppp{l}{V_{\!k\!,l\pl1}}}\epic
\text{1.5}{1.3}{0.6}{1}{}{=}
\bpic(2,1.3)\multiput(0,0.5)(0,1){2}{\line(1,0){2}}\multiput(0.5,0)(1,0){2}{\line(0,1){2}}
\multiput(0.5,1)(1,0){2}{\ppp{}{\bullet}}\multiput(1,0.5)(0,1){2}{\ppp{}{\bullet}}
\put(1,0.36){\pp{t}{0}}\put(1,1.64){\pp{b}{1}}\put(0.4,1){\pp{r}{1}}\put(1.6,1){\pp{l}{0}}\epic
\text{1.6}{1.3}{0.8}{1}{}{\mbox{or}}
\bpic(2,1.3)\multiput(0,0.5)(0,1){2}{\line(1,0){2}}\multiput(0.5,0)(1,0){2}{\line(0,1){2}}
\multiput(0.5,1)(1,0){2}{\ppp{}{\bullet}}\multiput(1,0.5)(0,1){2}{\ppp{}{\bullet}}
\put(1,0.36){\pp{t}{1}}\put(1,1.64){\pp{b}{0}}\put(0.4,1){\pp{r}{0}}\put(1.6,1){\pp{l}{1}}\put(2.4,0.9){\p{}{,}}\epic\\[3mm]
\quad\qquad\qquad\mbox{for each \ }(k,l)\in B(a,b,i,j)\!\setminus\!\{(a,a\pl j\mi i),(b\pl i\mi j,b)\}\,,\ea
\end{equation}
and
\begin{equation}\label{lem52}\ba{c@{}c@{\quad\;\,}c@{\quad\,}c@{\quad\,}c@{}l}
\bpic(3.3,1)\put(0,0.5){\line(1,0){2.2}}\multiput(0.5,0)(1.2,0){2}{\line(0,1){1}}
\multiput(0.5,0)(1.2,0){2}{\ppp{}{\bullet}}\put(1.1,0.5){\ppp{}{\bullet}}
\put(1.12,0.6){\ppp{b}{H_{a\!,a\pl j\mi i}}}
\put(0.95,-0.1){\ppp{tr}{V_{\!a,a\pl j\mi i}}}\put(1.32,-0.1){\ppp{tl}{V_{\!a\!,a\pl j\mi i\pl1}}}\epic
&\text{0}{1}{0}{0.5}{}{=}&
\bpic(2,1)\put(0,0.5){\line(1,0){2}}\multiput(0.5,0)(1,0){2}{\line(0,1){1}}
\multiput(0.5,0)(1,0){2}{\ppp{}{\bullet}}\put(1,0.5){\ppp{}{\bullet}}
\put(1,0.64){\pp{b}{1}}\put(0.5,-0.13){\pp{t}{1}}\put(1.5,-0.13){\pp{t}{0}}\epic
&\text{0}{1}{0}{0.5}{}{\mbox{or}}&
\bpic(1.8,1)\put(0,0.5){\line(1,0){2}}\multiput(0.5,0)(1,0){2}{\line(0,1){1}}
\multiput(0.5,0)(1,0){2}{\ppp{}{\bullet}}\put(1,0.5){\ppp{}{\bullet}}
\put(1,0.64){\pp{b}{0}}\put(0.5,-0.13){\pp{t}{0}}\put(1.5,-0.13){\pp{t}{1}}\epic
&\text{3.7}{1}{0.3}{0.5}{l}{\mbox{, \ \ if \ }j\mi i<b\mi a\,,}\\
\bpic(2.6,1.8)\put(0,0.5){\line(1,0){1}}\put(0.5,0){\line(0,1){1}}
\multiput(0.5,0)(0.5,0.5){2}{\ppp{}{\bullet}}
\put(1.13,0.45){\ppp{l}{H_{ab}}}\put(0.57,-0.13){\ppp{t}{V_{\!ab}}}\epic
&\text{0}{1}{0}{0.5}{}{=}&
\bpic(1,1)\put(0,0.5){\line(1,0){1}}\put(0.5,0){\line(0,1){1}}
\multiput(0.5,0)(0.5,0.5){2}{\ppp{}{\bullet}}
\put(1.13,0.5){\pp{l}{1}}\put(0.5,-0.13){\pp{t}{1}}\epic
&\text{0}{1}{0}{0.5}{}{\mbox{or}}&
\bpic(1,1)\put(0,0.5){\line(1,0){1}}\put(0.5,0){\line(0,1){1}}
\multiput(0.5,0)(0.5,0.5){2}{\ppp{}{\bullet}}
\put(1.13,0.5){\pp{l}{0}}\put(0.5,-0.13){\pp{t}{0}}\epic
&\text{3.7}{1}{0.3}{0.5}{l}{\mbox{, \ \ if \ }j\mi i=b\mi a\,,}\\
\bpic(2.6,2.8)\put(0.5,0){\line(0,1){2}}\multiput(0,0.5)(0,1){2}{\line(1,0){1}}
\multiput(1,0.5)(0,1){2}{\ppp{}{\bullet}}\put(0.5,1){\ppp{}{\bullet}}
\put(0.43,1){\ppp{r}{V_{b,b\pl i\mi j}}}
\put(1.1,1.45){\ppp{l}{H_{b\pl i\mi j,b}}}\put(1.1,0.45){\ppp{l}{H_{b\pl i\mi j\pl1,b}}}
\epic&\text{0}{2}{0}{1}{}{=}&
\bpic(1,2)\put(0.5,0){\line(0,1){2}}\multiput(0,0.5)(0,1){2}{\line(1,0){1}}
\multiput(1,0.5)(0,1){2}{\ppp{}{\bullet}}\put(0.5,1){\ppp{}{\bullet}}
\put(0.41,1){\pp{r}{1}}\put(1.13,1.5){\pp{l}{1}}\put(1.13,0.5){\pp{l}{0}}\epic
&\text{0}{2}{0}{1}{}{\mbox{or}}&
\bpic(1,2)\put(0.5,0){\line(0,1){2}}\multiput(0,0.5)(0,1){2}{\line(1,0){1}}
\multiput(1,0.5)(0,1){2}{\ppp{}{\bullet}}\put(0.5,1){\ppp{}{\bullet}}
\put(0.41,1){\pp{r}{0}}\put(1.13,1.5){\pp{l}{0}}\put(1.13,0.5){\pp{l}{1}}\epic
&\text{3.7}{2}{0.3}{1}{l}{\mbox{, \ \ if \ }j\mi i>b\mi a\,.}
\ea\end{equation}
\end{itemize}
\end{lemma}

\begin{figure}[h]
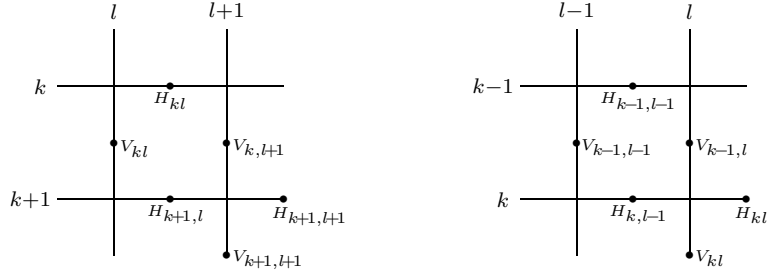
\centering\bpic(3,2.3)
\multiput(0.5,0.5)(0,1){2}{\line(1,0){2}}\multiput(1,0)(1,0){2}{\line(0,1){2}}
\multiput(1,1)(0.5,-0.5){3}{\ppp{}{\bullet}}\multiput(1.5,1.5)(0.5,-0.5){3}{\ppp{}{\bullet}}
\put(0.44,0.5){\pp{r}{k+1}}\put(0.4,1.5){\pp{r}{k}}
\put(1,2.1){\pp{b}{l}}\put(2,2.1){\pp{b}{l+1}}
\put(1.55,0.45){\ppp{t}{H_{k\pl1,l}}}\put(1.5,1.45){\ppp{t}{H_{kl}}}\put(2.4,0.43){\ppp{tl}{H_{k\pl1,l\pl1}}}
\put(1.05,0.95){\ppp{l}{V_{kl}}}\put(2.05,0.95){\ppp{l}{V_{k,l\pl1}}}\put(2.05,0.05){\ppp{tl}{V_{k\pl1,l\pl1}}}
\epic\qquad\qquad
\bpic(3,2.3)
\multiput(0.5,0.5)(0,1){2}{\line(1,0){2}}\multiput(1,0)(1,0){2}{\line(0,1){2}}
\multiput(1,1)(0.5,-0.5){3}{\ppp{}{\bullet}}\multiput(1.5,1.5)(0.5,-0.5){3}{\ppp{}{\bullet}}
\put(0.4,0.5){\pp{r}{k}}\put(0.44,1.5){\pp{r}{k-1}}
\put(1,2.1){\pp{b}{l-1}}\put(2,2.1){\pp{b}{l}}
\put(1.55,0.45){\ppp{t}{H_{k,l\mi1}}}\put(1.55,1.45){\ppp{t}{H_{k\mi1,l\mi1}}}\put(2.4,0.43){\ppp{tl}{H_{kl}}}
\put(1.04,0.95){\ppp{l}{V_{k\mi1,l\mi1}}}\put(2.05,0.95){\ppp{l}{V_{k\mi1,l}}}\put(2.05,0.05){\ppp{tl}{V_{kl}}}
\epic
\caption{Diagrams for the proof of Lemma 5.\label{lemfig}}
\end{figure}

Note that (\ref{lem52}) can be regarded as the boundary case of (\ref{lem51}) for
\begin{equation}\label{boundkl}(k,l)\,=\,(i\pl\min(a\mi i,b\mi j),j\pl\min(a\mi i,b\mi j))\,=\,
\left\{\ba{l}
(a,a\pl j\mi i),\;j\mi i<b\mi a\\[1.5mm]
(a,b),\;j\mi i=b\mi a\\[1.5mm]
(b\pl i\mi j,b),\;j\mi i>b\mi a\,.\ea\right.\end{equation}
\textit{Proof.} \  It will first be shown that the second condition implies the third.  Thus,
let the second condition be satisfied.  If $(k\pl1,l\pl1)\in B(a,b,i,j)$ (i.e., $k\ne a$ and $l\ne b$), then,
using the first diagram of Figure~\ref{lemfig} to visualize the positions of the relevant edge matrix entries,
$V_{k,l\pl1}=1\mi H_{kl}$ (since $(k,l\pl1)\notin S$)
and $H_{k\pl1,l}=1\mi V_{kl}$ (since $(k\pl1,l)\notin S$),
so that $H_{k\pl1,l}=V_{k,l\pl1}$, and therefore (using~(\ref{equiv2})) $H_{k\pl1,l\pl1}=V_{k\pl1,l\pl1}$.
By repeating this process, moving in a downward-rightward direction from $(k,l)$, it follows that $H_{k'l'}=V_{k'l'}$ for
all $(k',l')\in B(a,b,k,l)$.

Similarly, if $(k\mi1,l\mi1)\in B(a,b,i,j)$ (i.e., $(k,l)\ne(i,j)$), then,
using the second diagram of Figure~\ref{lemfig} to visualize the positions of the relevant edge matrix entries,
$H_{k,l\mi1}=V_{k\mi1,l}$ (using~(\ref{equiv2})),
$H_{k\mi1,l\mi1}=1\mi V_{k\mi1,l}$ (since $(k\mi1,l)\notin S$)
and $V_{k\mi1,l\mi1}=1\mi H_{k,l\mi1}$ (since $(k,l\mi1)\notin S$),
and therefore $H_{k\mi1,l\mi1}=V_{k\mi1,l\mi1}$.
By repeating this process, moving in a upward-leftward direction from $(k,l)$, it follows that $H_{k'l'}=V_{k'l'}$ for
all $(k',l')\in B(a,b,i,j)\!\setminus\!B(a,b,k\pl1,l\pl1)$.
\setlength{\unitlength}{15mm}

\vspace{-4mm}
Overall then, it follows that $H_{k'l'}=V_{k'l'}$ for
all $(k',l')\in B(a,b,i,j)$, as required.

Having now shown that the second condition implies the third, the equivalence
between the second and third conditions follows immediately, while the equivalence of these conditions
with the first condition follows from Lemma~3.

Finally, it can be seen that the third and fourth conditions are equivalent, using
the equivalence between the first and fourth conditions of~(\ref{equiv1}).\hspace{\fill}$\Box$

It follows immediately from Lemma~5 that if $(i,j)$ is a change point of $S\in\VOS(a,b)$,
then all elements of $B(a,b,i\pl1,j\pl1)$ are change points of $S$ as well.

It also follows from Lemmas~4 and~5 that if $S,S'\in\VOS(a,b)$ differ by the vacancy-osculation $(i,j)$, then
whichever of the two alternatives occurs on any RHS of~(\ref{lem51}) or the RHS of (\ref{lem52}) for the edge matrix pair corresponding
to~$S$, the other alternative occurs on that RHS for the edge matrix pair corresponding to~$S'$.

\begin{lemma} \ \ \ Consider vacancy-osculation sets $S\in\VOS(a,b,\a,\b)$ and
$S'\in\VOS(a,b,\a',\b')$.  If $S$ and $S'$ differ by a vacancy-osculation, then $\l_{a,b,\a,\b}$ and $\l_{a,b,\a',\b'}$ differ by a square, and
\begin{equation}\label{DDVOSPar}\D(S,S')\;=\;\D(\l_{a,b,\a,\b},\l_{a,b,\a',\b'})\,.\end{equation}
\end{lemma}
For $S'\in\VOS(4,6,\a',\b')$ corresponding to each change point $(i,j)$ of the running example $S$, Lemma~6
states that $(3,2,2)\sim\l_{4,6,\a',\b'}$ and $\D(S,S')=j\mi i=\D((3,2,2),\l_{4,6,\a',\b'})$,
which can be observed in Table~\ref{ExPar}.  It can also be seen in Table~\ref{ExPar} and Figure~\ref{ExCP} that
for cases of Lemma~6, it does not always follow that $|S|-|S'|=|\l_{a,b,\a,\b}|-|\l_{a,b,\a',\b'}|$.

\textit{Proof.} \ Let $S$ and $S'$ differ by the vacancy-osculation $(i,j)$,
so that $\D(S,S')=j\mi i$, and let $(H,V)\in\EM(a,b,\a,\b)$ and $(H',V')\in\EM(a,b,\a',\b')$
correspond to~$S$ and~$S'$ respectively.
Focussing on the edge matrix entries along the right and lower boundaries
of the rectangle, it follows from Lemmas~4 and~5 that~(\ref{HVD}) and~(\ref{lem52}) are satisfied, which gives,
for each $k\in[a]$ and $l\in[b]$:
\[\ba{ll}
\bullet&H_{kb}\,=\,H'_{kb}, \ \ |V_{al}\mi V'_{al}|\,=\,\d_{l\in\{a\pl j\mi i,\,a\pl j\mi i\pl1\}}\mbox{ \ \ and}\\[1.5mm]
&(V_{a,a\pl j\mi i},V_{a,a\pl j\mi i\pl1},V'_{a,a\pl j\mi i},V'_{a,a\pl j\mi i\pl1})=(1,0,0,1)\mbox{ or }(0,1,1,0),\mbox{ \ if }j\mi i<b\mi a.\\[3mm]
\bullet&|H_{kb}\mi H'_{kb}|\,=\,\d_{ak}, \ \ |V_{al}\mi V'_{al}|\,=\,\d_{bl}\mbox{ \ \ and}\\[1.5mm]
&(H_{ab},V_{ab},H'_{ab},V'_{ab})=(1,1,0,0)\mbox{ or }(0,0,1,1),\mbox{ \ if }j\mi i=b\mi a.\\[3mm]
\bullet&|H_{kb}\mi H'_{kb}|\,=\,\d_{k\in\{b\pl i\mi j,\,b\pl i\mi j\pl1\}}, \ \ V_{al}\,=\,V'_{al}\mbox{ \ \ and}\\[1.5mm]
&(H_{b\pl i\mi j,b},H_{b\pl i\mi j\pl1,b},H'_{b\pl i\mi j,b},H'_{b\pl i\mi j\pl1,b})=(1,0,0,1)\mbox{ or }(0,1,1,0),\mbox{ \ if }j\mi i=b\mi a.\ea\]
This gives six separate cases, each of which leads to the conclusions of Lemma~6.
The details of the case $j\mi i<b\mi a$ with $(V_{a,a\pl j\mi i},V_{a,a\pl j\mi i\pl1},V'_{a,a\pl j\mi i},V'_{a,a\pl j\mi i\pl1})
=(1,0,0,1)$ will now be outlined,
the other cases being reasonably similar.  Setting
$\a=\{\a_1,\ldots,\a_r\}$, $\a'=\{\a'_1,\ldots,\a'_{r'}\}$,
$\b=\{\b_1,\ldots,\b_r\}$ and $\b'=\{\b'_1,\ldots,\b'_{r'}\}$, with the elements of each set
labeled in ascending order, it can be seen for the case under consideration that $r=r'$,
$\a=\a'$, and that there exists $m\in[r]$ such that $\b_m=a\pl j\mi i$,
$\b'_m=a\pl j\mi i\pl1$ and $\b'_l=\b_l$, for all $l\in[r]\!\setminus\!\{m\}$.
Now setting $\m=(b\mi\b_1,\ldots,b\mi\b_r\,|\,a\mi\a_1,\ldots,a\mi\a_r)$ and
$\m'=(b\mi\b'_1,\ldots,b\mi\b'_r\,|\,a\mi\a'_1,\ldots,a\mi\a'_r)$, it can be seen
using~(\ref{FN}) that $\m_m=m\pl b\mi\b_m=m\mi a\pl b\pl i\mi j=\m'_m\pl1$ and $\m_k\!=\!\m'_k$, for all $k\in\P\setminus\!\{m\}$.
Therefore, $Y(\m)$ is the disjoint union of $Y(\m')$ and $\{(m,m\mi a\pl b\pl i\mi j)\}$.  Finally, since
$\l_{a,b,\a,\b}=[a]\t[b]\!\setminus\!\m$ and $\l_{a,b,\a',\b'}=[a]\t[b]\!\setminus\!\m'$ from~(\ref{ToP}),
it follows that $Y(\l_{a,b,\a',\b'})$ is the disjoint union of $Y(\l_{a,b,\a,\b})$
and $\{(a\mi m\pl1,a\mi i\pl j\mi m\pl1)\}$, so that
$\l_{a,b,\a,\b}\sim\l_{a,b,\a',\b'}$ and
$\D(\l_{a,b,\a,\b},\l_{a,b,\a',\b'})=j\mi i=\D(S,S')$, as required.\hspace{\fill}$\Box$

\begin{lemma}Consider a vacancy-osculation set $S\in\VOS(a,b,\a,\b)$ and a point $(i,j)\in[a]\t[b]\!\setminus\!S$.
Then $(i,j)$ is an addition point of $S$ if and only if $S\cap(B(a,b,i,j\pl1)\cup B(a,b,i\pl1,j))=\emptyset$ and
$j\mi i$ is a change diagonal of $\l_{a,b,\a,\b}$.
\end{lemma}
\textit{Proof.} \ If $(i,j)$ is an addition point of $S$, then it follows from
Lemma~3 or~5 that $S\cap(B(a,b,i,j\pl1)\cup B(a,b,i\pl1,j))=\emptyset$, and it follows from
Lemma~6 that $j\mi i$ is a change diagonal of $\l_{a,b,\a,\b}$ (since
$(i,j)$ being an addition point of $S$ means that
there exists $S'\in\VOS(a,b,\a',\b')$ with $S\sim S'$ and $S'=S\cup\{(i,j)\}$, for
some $(\a',\b')\in\BP(a,b)$, and Lemma~6 then implies that $\l_{a,b,\a,\b}\sim\l_{a,b,\a',\b'}$ and
$\D(\l_{a,b,\a,\b},\l_{a,b,\a',\b'})=\D(S,S')=j\mi i$).

Conversely, if $j\mi i$ is a change diagonal of $\l_{a,b,\a,\b}$ and
$S\cap(B(a,b,i,j\pl1)\cup B(a,b,i\pl1,j))=\emptyset$, then
it follows from Lemma~5 that $(i,j)$ is an addition point of $S$ provided that
$H_{kl}=V_{kl}$ for some $(k,l)\in B(a,b,i,j)$, where $(H,V)\in\EM(a,b)$ corresponds to $S$.
It will be shown that $H_{kl}=V_{kl}$ is satisfied for $(k,l)$ given by~(\ref{boundkl}).
First, it follows from~(\ref{lem2}) and~(\ref{lem23}) that, for any $d\in\Z$,
\begin{equation}\label{rho1}\rho_d(\l_{a,b,\a,\b})\;=\;R_d([a]\t[b])-
\left\{\ba{l@{\;\:}l}|\{i\in[r]\,|\,\b_i\le a\pl d\}|,&d\le b\mi a\\[2.5mm]
|\{i\in[r]\,|\,\a_i\le b\mi d\}|,&d\ge b\mi a\,.\ea\right.\end{equation}
This can be rewritten as
\begin{equation}\label{rho2}\rho_d(\l_{a,b,\a,\b})\;=\;
\left\{\ba{l}\min(a,b,a\pl d,b\mi d),\;d\in[-a,b]\\[2.5mm]0,\;d\notin[-a,b]\ea\right.\;\;\;-\;\;\;
\left\{\ba{l}\sum_{l=1}^{a\pl d}V_{al},\;d\le b\mi a\\[2.9mm]
\sum_{k=1}^{b\mi d}H_{kb},\;d\ge b\mi a\,,\ea\right.\end{equation}
the first terms on each RHS of~(\ref{rho1}) and~(\ref{rho2}) being equal by elementary
geometry of the rectangle.

Next, it can be shown using~(\ref{rhoCond}) and~(\ref{rho2}) that only the possibilities
\begin{equation}\label{boundHV}\ba{l}
\bullet\;\;(V_{a,a\pl j\mi i},V_{a,a\pl j\mi i\pl1})=(1,0)\mbox{ or }(0,1),\mbox{ \ if \ }j\mi i<b\mi a\\[2.5mm]
\bullet\;\;(H_{ab},V_{ab})=(1,1)\mbox{ or }(0,0),\mbox{ \ if \ }j\mi i=b\mi a\\[2.5mm]
\bullet\;\;(H_{b\pl i\mi j,b},H_{b\pl i\mi j\pl1,b})=(1,0)\mbox{ or }(0,1),\mbox{ \ if \ }j\mi i>b\mi a\,,\ea\end{equation}
can occur for the edge matrix entries at the two edges adjacent to the diagonal $D_{j\mi i}([a]\t[b])$ on
the right or lower boundary.
The details of the case $-a<j\mi i<\min(b\mi a,0)$
will now be outlined, all other cases being similar.
For this case,~(\ref{rho2}) gives $\rho_{j\mi i\mi1}(\l_{a,b,\a,\b})=a\pl j\mi i\mi1\mi\sum_{l=1}^{a\pl j\mi i\mi1}V_{al}$,
$\rho_{j\mi i}(\l_{a,b,\a,\b})=a\pl j\mi i\mi\sum_{l=1}^{a\pl j\mi i}V_{al}$ and
$\rho_{j\mi i\pl1}(\l_{a,b,\a,\b})=a\pl j\mi i\pl1\mi\sum_{l=1}^{a\pl j\mi i\pl1}V_{al}$.
Therefore, $(\rho_{j\mi i}(\l_{a,b,\a,\b})-\rho_{j\mi i\mi1}(\l_{a,b,\a,\b}),\rho_{j\mi i\pl1}(\l_{a,b,\a,\b})-\rho_{j\mi i}(\l_{a,b,\a,\b}))=
(1\mi V_{a,a\pl j\mi i},1\mi V_{a,a\pl j\mi i\pl1})$.  But, using~(\ref{rhoCond}), this must be $(1,0)$ or $(0,1)$, implying that
$(V_{a,a\pl j\mi i},V_{a,a\pl j\mi i\pl1})=(1,0)$ or $(0,1)$, as required.

It can now be seen that, as required, $H_{kl}\,=\,V_{kl}$ for $(k,l)$ given by~(\ref{boundkl}),
this following from~(\ref{boundHV}) and $(a,a\pl j\mi i\pl1)\notin S$ (since $(a,a\pl j\mi i\pl1)\in B(a,b,i,j\pl1)$) if $j\mi i<b\mi a$,
from~(\ref{boundHV}) alone if $j\mi i=b\mi a$, and
from~(\ref{boundHV}) and $(b\pl i\mi j\pl1,b)\notin S$ (since $(b\pl i\mi j\pl1,b)\in B(a,b,i\pl1,j)$) if $j\mi i>b\mi a$.\hspace{\fill}$\Box$

\subsubsection{Vacancy-Osculation Set Progressions}
In this section, it will be shown that the elements of any vacancy-osculation set $S$ can be ordered in such a way that, for each $k\in[|S|]$,
the first $k$ elements themselves constitute a vacancy-osculation set.

For a finite set $T$, let an \emph{ordering} of $T$ be any permutation, represented as a tuple, of all the elements of $T$.
If $T$ is a subset of $\Z^2$, let a \emph{canonical ordering} of $T$ be any
ordering $((i_1,j_1),\ldots,(i_{|T|},j_{|T|}))$ for which
if $k<l$, with $k,l\in[|T|]$, then $i_k<i_l$ or $j_k<j_l$.
Furthermore, let the \emph{lexicographic ordering} of $T$ be the
unique canonical ordering $((i_1,j_1),\ldots,(i_{|T|},j_{|T|}))$ which satisfies
$i_k<i_{k\pl1}$, or $i_k=i_{k\pl1}$ and $j_k<j_{k\pl1}$, for each $k\in[|T|\mi1]$.
(For example, in~(\ref{ExS}) the points of $S$ are written in such an order.)

For a vacancy-osculation set $S\in\VOS(a,b,l)$, define a \emph{progression} for $S$ to be any
ordering $s=(s_1,\ldots,s_l)$ of $S$ which satisfies $\{s_1,\ldots,s_k\}\in\VOS(a,b)$ for each $k\in[0,l]$.
A progression $s$ for $S$ therefore satisfies $\emptyset\!\sim\!\{s_1\}\!\sim\!\{s_1,s_2\}\!\sim\!\ldots\!\sim\!
\{s_1,\ldots,s_l\}$, with~$s_k$ being a deletion point of $\{s_1,\ldots,s_k\}$
and an addition point of $\{s_1,\ldots,s_{k\mi1}\}$ for each $k\in[l]$.

\begin{lemma}Any canonical ordering of a vacancy-osculation set is a progression for that set.\end{lemma}
\textit{Proof.} \ Let $s=(s_1,\ldots,s_l)=((i_1,j_1),\ldots,(i_l,j_l))$ be a canonical ordering of $S\in\VOS(a,b,l)$, and
consider $\{s_1,\ldots,s_m\}$ for any $m\in[l]$.
Since $s$ is canonical, $i_k<i_m$ or $j_k<j_m$ for each $k\in[m\mi1]$, from which it follows,
due to~(\ref{B}), that $\{s_1,\ldots,s_m\}\cap(B(a,b,i_m,j_m\pl1)\cup B(a,b,i_m\pl1,j_m))=\emptyset$.
Therefore, if $\{s_1,\ldots,s_m\}\in\VOS(a,b)$, then, by Lemma~3,~$s_m$ is a deletion point of
$\{s_1,\ldots,s_m\}$, and so $\{s_1,\ldots,s_{m\mi1}\}\in\VOS(a,b)$.  Since
$\{s_1,\ldots,s_l\}=S\in\VOS(a,b)$, it follows by repeated application
of this result starting with $m=l$ that $\{s_1,\ldots,s_m\}\in\VOS(a,b)$
for each $m\in[0,l]$, so that~$s$ is a progression for~$S$.\hspace{\fill}$\Box$

\begin{corollary}There is at least one progression for any vacancy-osculation set.\end{corollary}
\textit{Proof.} \ Since the lexicographic ordering of a vacancy-osculation set is a canonical ordering,
it follows from Lemma~8 that it is a progression for that set.\hspace{\fill}$\Box$

An example of a canonical, but non-lexicographic, ordering of the vacancy-osculation set~(\ref{ExS}) of the running example,
and hence a progression for that set, is
\begin{equation}\label{ExProg}s=((1,1),(1,2),(1,3),(2,1),(2,2),(2,3),(1,4),(3,4),(2,5),(4,2),(4,3))\,.\end{equation}

\setlength{\unitlength}{1mm}
\begin{figure}[h]\centering\bpic(150,16)\scriptsize
\put(13,8){\p{}{\left(\ba{c@{\:\:}c@{\:\:}c@{\:\:}c@{\:\:}c@{\:\:}c}0&0&0&0&0&0\\0&0&0&0&0&0\\0&0&0&0&0&0\\0&0&0&0&0&0\ea\right)}}
\put(37.8,8){\p{}{\left(\ba{c@{\:\:}c@{\:\:}c@{\:\:}c@{\:\:}c@{\:\:}c}\mathbf{1}&0&0&0&0&0\\0&0&0&0&0&0\\0&0&0&0&0&0\\0&0&0&0&0&0\ea\right)}}
\put(62.5,8){\p{}{\left(\ba{c@{\:\:}c@{\:\:}c@{\:\:}c@{\:\:}c@{\:\:}c}1&\mathbf{1}&0&0&0&0\\0&0&0&0&0&0\\0&0&0&0&0&0\\0&0&0&0&0&0\ea\right)}}
\put(87.5,8){\p{}{\left(\ba{c@{\:\:}c@{\:\:}c@{\:\:}c@{\:\:}c@{\:\:}c}1&1&\mathbf{1}&0&0&0\\0&0&0&0&0&0\\0&0&0&0&0&0\\0&0&0&0&0&0\ea\right)}}
\put(112.2,8){\p{}{\left(\ba{c@{\:\:}c@{\:\:}c@{\:\:}c@{\:\:}c@{\:\:}c}1&1&1&0&0&0\\\mathbf{1}&0&0&0&0&0\\0&0&0&0&0&0\\0&0&0&0&0&0\ea\right)}}
\put(137,8){\p{}{\left(\ba{c@{\:\:}c@{\:\:}c@{\:\:}c@{\:\:}c@{\:\:}c}1&1&1&0&0&0\\1&\mathbf{1}&0&0&0&0\\0&0&0&0&0&0\\0&0&0&0&0&0\ea\right)}}\epic\\
\includegraphics[width=150mm]{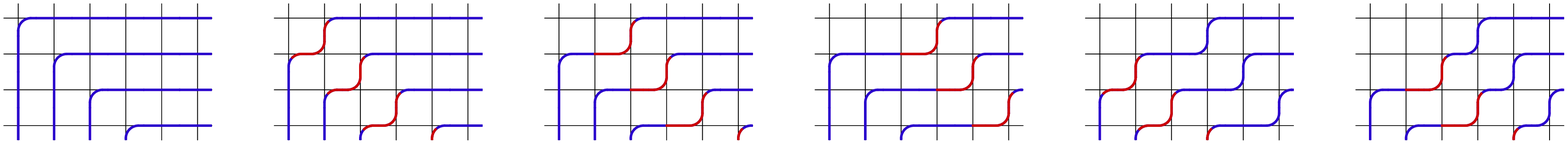}\\
\bpic(150,12.5)
\put(37.8,9){\pp{}{(1,1)}}\put(62.5,9){\pp{}{(1,2)}}\put(87.5,9){\pp{}{(1,3)}}\put(112.2,9){\pp{}{(2,1)}}\put(137,9){\pp{}{(2,2)}}
\epic\\
\bpic(150,13)\scriptsize
\put(13,8){\p{}{\left(\ba{c@{\:\:}c@{\:\:}c@{\:\:}c@{\:\:}c@{\:\:}c}1&1&1&0&0&0\\1&1&\mathbf{1}&0&0&0\\0&0&0&0&0&0\\0&0&0&0&0&0\ea\right)}}
\put(37.8,8){\p{}{\left(\ba{c@{\:\:}c@{\:\:}c@{\:\:}c@{\:\:}c@{\:\:}c}1&1&1&\mathbf{1}&0&0\\1&1&1&0&0&0\\0&0&0&0&0&0\\0&0&0&0&0&0\ea\right)}}
\put(62.5,8){\p{}{\left(\ba{c@{\:\:}c@{\:\:}c@{\:\:}c@{\:\:}c@{\:\:}c}1&1&1&1&0&0\\1&1&1&0&0&0\\0&0&0&\mathbf{1}&0&0\\0&0&0&0&0&0\ea\right)}}
\put(87.5,8){\p{}{\left(\ba{c@{\:\:}c@{\:\:}c@{\:\:}c@{\:\:}c@{\:\:}c}1&1&1&1&0&0\\1&1&1&0&\mathbf{1}&0\\0&0&0&1&0&0\\0&0&0&0&0&0\ea\right)}}
\put(112.2,8){\p{}{\left(\ba{c@{\:\:}c@{\:\:}c@{\:\:}c@{\:\:}c@{\:\:}c}1&1&1&1&0&0\\1&1&1&0&1&0\\0&0&0&1&0&0\\0&\mathbf{1}&0&0&0&0\ea\right)}}
\put(137,8){\p{}{\left(\ba{c@{\:\:}c@{\:\:}c@{\:\:}c@{\:\:}c@{\:\:}c}1&1&1&1&0&0\\1&1&1&0&1&0\\0&0&0&1&0&0\\0&1&\mathbf{1}&0&0&0\ea\right)}}\epic\\
\includegraphics[width=150mm]{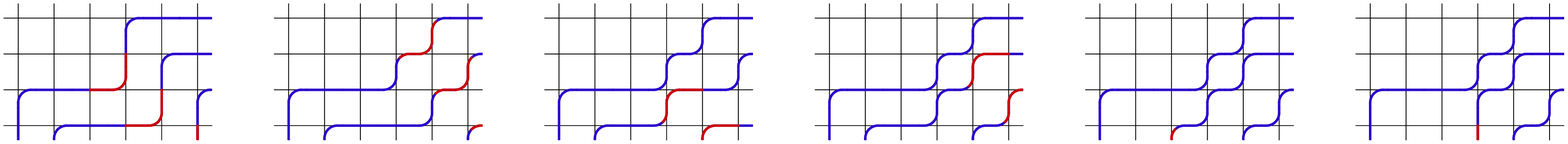}\\
\bpic(150,5)
\put(13,1.5){\pp{}{(2,3)}}\put(37.8,1.5){\pp{}{(1,4)}}\put(62.5,1.5){\pp{}{(3,4)}}
\put(87.5,1.5){\pp{}{(2,5)}}\put(112.2,1.5){\pp{}{(4,2)}}\put(137,1.5){\pp{}{(4,3)}}
\epic
\caption{A progression for the running example.\label{ExProgFig}}
\end{figure}

In Figure~\ref{ExProgFig}, each $s_k$ for this example is shown together with $M(\{s_1,\ldots,s_k\})$
and the diagram of the path tuple corresponding to $\{s_1,\ldots,s_k\}$, for $k\in[0,11]$.
In each matrix, the $s_k$ entry is highlighted, and in each path diagram the sections of paths which differ
from the previous diagram are shown in red.

Note that the converse of Lemma~8 is false.  For example, the ordering obtained from~(\ref{ExProg}) by
interchanging $(1,4)$ and $(2,5)$ is no longer canonical, but is still a progression for the vacancy-osculation set~(\ref{ExS}).

\subsubsection{Ordering Pairs}
In this section, it will be shown that progressions can be systematically
obtained for any vacancy-osculation set by using a function on $\P^2$ and an order on $\Z$ which
satisfy certain properties.

A \emph{total strict order} on a set $T$ is a binary relation $\prec$ on $T$ for which:
\begin{equation}\ba{rl}(1)&\mbox{Exactly one of }t_1\prec t_2, \ t_1=t_2\mbox{ or }t_2\prec t_1\mbox{ holds for each }t_1, t_2\in T.\\[2mm]
(2)&\mbox{If }t_1,t_2,t_3\in T\mbox{ satisfy }t_1\prec t_2\mbox{ and }t_2\prec t_3,\mbox{ then }t_1\prec t_3.\ea\end{equation}
Of primary interest here will be total strict orders on $\Z$.
The standard less-than relation $<$ is one example of a total strict order on~$\Z$, while another example
is that defined, for each $z,z'\in\Z$, by
\begin{equation}\label{ExTSO}z\prec z'\mbox{ \ if and only if \ }|z\mi2|>|z'\mi2|\mbox{ \ or \ }z\mi2=2\mi z'<0\,,\end{equation}
which gives \ $\ldots\prec-1\prec5\prec0\prec4\prec1\prec3\prec2$.

Let an \emph{ordering pair} $(L,\prec)$ be any pair comprised of a function $L:\P^2\rightarrow\P$ and a total strict order $\prec$ on~$\Z$,
which satisfy the property, for $(i,j),(i',j')\in\P^2$, that
\begin{equation}\label{OrdPair}\ba{c}\mbox{If \ }L(i,j)<L(i',j'),\mbox{ \ or \ }L(i,j)=L(i',j')\mbox{ \ and \ }j\mi i\prec j'\mi i',\\[1.6mm]
\mbox{then \ }i<i'\mbox{ \ or \ }j<j'\,.\ea\end{equation}
The function $L$ will be referred to as a \emph{level function}, and $L(i,j)$ will be referred to as the \emph{level} of $(i,j)\in\P^2$.
The full motivation for the use in~(\ref{OrdPair}) of a total strict order on the diagonals of $\P^2$
(i.e., the use of $j\mi i\prec j'\mi i'$ rather than a more general condition for the case $L(i,j)=L(i',j')$)
will become apparent in Section~13.

It follows, for an ordering pair $(L,\prec)$, and $(i,j),(i',j')\in\P^2$, that
\begin{equation}\label{OrdPairEq}\mbox{If \ }L(i,j)=L(i',j')\mbox{ \ and \ }j\mi i=j'\mi i',\mbox{ \ then \ }(i,j)=(i',j')\,.\end{equation}
In other words, each point $(i,j)\in\P^2$ is uniquely determined by its level $L(i,j)$ and content $j\mi i$.
Property~(\ref{OrdPairEq}) can be derived as follows. Let $L(i,j)=L(i',j')$ and $j\mi i=j'\mi i'$, and consider, for example, the
point $(i,j')$. By applying~(\ref{OrdPair}) to $(i,j')$ and $(i,j)$, and to $(i,j')$ and $(i',j')$,
it can be seen that none of the possibilities $L(i,j')<L(i,j)=L(i',j')$, $L(i,j)=L(i',j')<L(i,j')$,
$L(i,j')=L(i,j)=L(i',j')$ and $j'\mi i\prec j\mi i=j'\mi i'$, or $L(i,j')=L(i,j)=L(i',j')$
and $j\mi i=j'\mi i'\prec j'\mi i$, can occur, leaving only
$j'\mi i=j\mi i=j'\mi i'$, and therefore $(i,j)=(i',j')$.

It follows from~(\ref{OrdPair}) and~(\ref{OrdPairEq}) that for an ordering pair $(L,\prec)$
and a vacancy-oscu\-lation set $S\in\VOS(a,b,l)$,
there exists a unique canonical ordering $s=((i_1,j_1),\ldots,$ $(i_l,j_l))$ of~$S$ which satisfies
\begin{equation}\label{IndProg}L(i_k,j_k)<L(i_{k\pl1},j_{k\pl1}),\mbox{ \ or }
L(i_k,j_k)=L(i_{k\pl1},j_{k\pl1})\mbox{ and }j_k\mi i_k\prec j_{k\pl1}\mi i_{k\pl1}\,,\end{equation}
for each $k\in[l\mi1]$.
(Alternatively, a total strict order $\prec$ on $\P^2$ can be defined by $(i,j)\prec(i',j')$ if and only if
$L(i,j)<L(i',j')$, or $L(i,j)=L(i',j')$ and $j\mi i\prec j'\mi i'$, and $s=((i_1,j_1),\ldots,(i_l,j_l))$ is then the unique
ordering of $S$ which satisfies $(i_k,j_k)\prec (i_{k\pl1},j_{k\pl1})$ for each $k\in[l\mi1]$.)
By Lemma~8, this $s$ is a progression for $S$.
Accordingly,~$s$ will be referred to as \emph{the progression for $S$ induced by the ordering pair $(L,\prec)$}.

An example of an ordering pair is $(\L,<)$, where the level function is
\begin{equation}\label{lex}\L(i,j)=i,\mbox{ \ for each }(i,j)\in\P^2\,,\end{equation}
and the total strict order is the standard less-than relation, since it can be seen immediately that these satisfy~(\ref{OrdPair}).
The progression for a vacancy-osculation set induced by $(\L,<)$ is the lexicographic ordering of the set, and
accordingly $(\L,<)$ will be referred to as the \emph{lexicographic ordering pair}.
Similarly, another example of an ordering pair is $(\LL,>)$, where the level function is
\begin{equation}\label{alex}\LL(i,j)=j,\mbox{ \ for each }(i,j)\in\P^2\,,\end{equation}
and the total strict order is the standard greater-than relation. This will be referred to
as the \emph{antilexicographic ordering pair}.

A third example of an ordering pair is $(L,\prec)$, where the level function is
\begin{equation}\label{ExL}L(i,j)=\max(i,j\mi2),\mbox{ \ for each }(i,j)\in\P^2\,,\end{equation}
and the total strict order is given by the example~(\ref{ExTSO}).  In checking that these
satisfy~(\ref{OrdPair}), it may be useful to visualize~(\ref{ExL}) as
\begin{equation}\label{ExLM}
\left(\ba{c@{\;\;}c@{\;}c}\ss L(1,1)&\ss L(1,2)&\ldots\\[-0.5mm]
\ss L(2,1)&\ss L(2,2)\\[-1mm]
\vdots&&\ddots\ea\right)\;=\;
\left(\ba{c@{\;\;}c@{\;\;}c@{\;\;}c@{\;\;}c@{\,}c}
\ss1&\ss1&\ss1&\ss2&\ss3\\[-0.9mm]
\ss2&\ss2&\ss2&\ss2&\ss3\\[-0.9mm]
\ss3&\ss3&\ss3&\ss3&\ss3\\[-2.4mm]
&&&&&\ddots\ea\right).\end{equation}
It can be seen that the progression for the running example~(\ref{ExS}) induced by this
ordering pair is that of~(\ref{ExProg}).

\subsubsection{Oscillating Tableaux}
In this section, oscillating tableaux will be introduced, in preparation for the consideration
in subsequent sections of their relationship with vacancy-osculation sets.

For a partition~$\l$ and a nonnegative integer~$l$, an \emph{oscillating tableau} of \emph{shape}~$\l$ and \emph{length}~$l$
is a sequence
of $l\pl1$ partitions starting with $\emptyset$, ending with $\l$, and in which successive partitions
differ by a square.  Let $\OT(\l,l)$ denote the set of all such oscillating tableaux,
\begin{equation}\label{OT}\ba{l}\OT(\l,l)\,:=\\[2mm]
\hspace{9.9mm}\Bigl\{\e=(\e_0\!=\!\emptyset,\e_1,\ldots,\e_{l\mi1},\e_l\!=\!\l)\:\in\:\Par^{l\pl1}\;\Big|\:\:\e_{k\mi1}\sim\e_k
\mbox{ for each }k\in[l]\Bigr\}\,.\ea\end{equation}
Oscillating tableaux are sometimes called \emph{up-down tableaux}, and were first introduced
in~\cite{Ber86}.  Oscillating tableaux, or certain generalizations thereof, are studied
in the context of Robinson-Schensted-Knuth-type correspondences
in~\cite{Ber86,ChaDul02,Cho03,DelDulFav88,DulSag95,Ges93,GouCha06,Pro91,Rob91b,Rob95,RobTer05,Sun86,Sun90,Ter01},
and in the context of nonintersecting lattice paths in~\cite{BakFor01,Kra96}.
They are also discussed in, for example,~\cite{Sag88},~\cite[Exercise~7.24]{Sta99} and~\cite[Sec.~9]{Sta06}.

Let a \emph{deletion} of an oscillating tableau $\e\in\OT(\l,l)$ be a value $k\in[l]$ for which $|\e_k|=|\e_{k\mi1}|\mi1$.
Denoting the number of deletions in $\e$ as $\mathcal{D}(\e)$,
it can be seen that
\begin{equation}\label{del}\mathcal{D}(\e)=(l\mi|\l|)/2\,,\mbox{
\ for each }\e\in\OT(\l,l)\,,\end{equation}
and that $\OT(\l,l)=\emptyset$ unless $(l\mi|\l|)/2$ is a nonnegative integer.

For an oscillating tableau $\e\in\OT(\l,l)$, define
\begin{equation}\label{OTPro}\O(\e)_k\::=\:\D(\e_{k\mi1},\e_k),\mbox{ \ for each }k\in[l],\end{equation}
and define the \emph{profile} $\O(\e)$ of $\e$ to be the sequence of these diagonal differences,
\begin{equation}\label{OTProSeq}\O(\e)\::=\:(\O(\e)_1,\O(\e)_2,\ldots,\O(\e)_l)\:=
\:(\D(\e_0,\e_1),\D(\e_1,\e_2),\ldots,\D(\e_{l\mi1},\e_l)).\end{equation}
In other words, $\O(\e)$ is the sequence of contents of the squares by which the Young diagrams of successive
partitions of~$\e$ differ.

It can be seen that an oscillating tableau~$\e$ is uniquely determined by its profile~$\O(\e)$,
since $\e_0$ is always~$\emptyset$ and, as indicated in Section~6, for fixed $\m\in\Par$, each $\m'\in\Par$ with $\m\!\sim\!\m'$
is uniquely determined by~$\D(\m,\m')$.

It can also be seen that, for any $\e\in\OT(\l,l)$,
\begin{equation}\label{OTProp}
\mbox{For each }d\in[-\ell(\l)\pl1,\l_1\mi1],\mbox{ there exists  }k\in[l]\mbox{ with }\O(\e)_k=d\,.\end{equation}

It can be deduced straightforwardly that, for any $\e\in\OT(\l,l)$ and $k\in[l]$,
\begin{equation}\label{SDD}\ba{l@{\;\;}l}
\bullet&\mbox{If }\O(\e)_k<0,\mbox{ then there exists }k'\in[k\mi1]\mbox{ with }\O(\e)_{k'}=\O(\e)_k\pl1.\\[3mm]
\bullet&\mbox{If }\O(\e)_k>0,\mbox{ then there exists }k'\in[k\mi1]\mbox{ with }\O(\e)_{k'}=\O(\e)_k\mi1.\\[3mm]
\bullet&\mbox{If }\O(\e)_k<-\ell(\l),\mbox{ then there exists }k'\in[k\pl1,l]\mbox{ with }\O(\e)_{k'}=\O(\e)_k\pl1.\\[3mm]
\bullet&\mbox{If }\O(\e)_k>\l_1,\mbox{ then there exists }k'\in[k\pl1,l]\mbox{ with }\O(\e)_{k'}=\O(\e)_k\mi1.\ea
\end{equation}

For an oscillating tableau $\e\in\OT(\l,l)$ and a total strict order $\prec$ on $\Z$, define the set of \emph{ascents of~$\e$
with respect to~$\prec$} as
\begin{equation}\label{OTAsc}\A(\prec,\e)\;:=\;\{k\in[l\mi1]\mid\O(\e)_k\prec\O(\e)_{k\pl1}\}\,,\end{equation}
and the \emph{number of ascents of~$\e$ with respect to~$\prec$} as
\begin{equation}\label{OTAscN}\AA(\prec,\e)\;:=\;|\A(\prec,\e)|\,.\end{equation}
Also define the distribution function on $\OP(\l,l)$ of the number of ascents with respect to $\prec$ as
\begin{equation}\label{Adist}\mathcal{N}(\l,l,\prec,i):=|\{\e\in\OT(\l,l)\mid\AA(\prec,\e)=i\}|,\mbox{ \ for each }i\in[0,l\mi1]\,.\end{equation}

An example of an oscillating tableau $\e$ of shape $(3,2,2)$ and length $11$ is
\begin{equation}\label{ExOT}\ba{l}\e=\\[1mm]
\;\quad\Bigl(\emptyset,\,(1),\,(2),\,(3),\,(3,1),\,(3,2),\,(3,3),\,(4,3),\,(4,2),\,(3,2),\,(3,2,1),\,(3,2,2)\Bigr).\ea\end{equation}
The Young diagrams for each partition of $\e$, and the terms of the profile of~$\e$, are given in Table~\ref{ExOTTab}.
It can be seen that $8$ and $9$ are the deletions of~$\e$, and that
for the total strict order on $\Z$ given in~(\ref{ExTSO}),
\begin{equation}\label{ExAsc}\A(\prec,\e)=\{1,2,4,5,6,8,10\}.\end{equation}

\setlength{\unitlength}{2.3mm}
\begin{table}[h]\centering\footnotesize$\ba{c||c|c|c|c|c|c|c|c|c|c|c|c}
k&0&1&2&3&4&5&6&7&8&9&10&11\\[1mm]\hline
\rule{0mm}{5mm}\e_k&\ss\emptyset&\ss(1)&\ss(2)&\ss(3)&\ss(3,1)&\ss(3,2)&\ss(3,3)&\ss(4,3)&\ss(4,2)&\ss(3,2)&\ss(3,2,1)&\ss(3,2,2)\\
&\bpic(0.6,4.3)\put(0.3,2.5){\pp{}{\emptyset}}\epic&
\bpic(1,3)\multiput(0,2)(0,1){2}{\line(1,0){1}}\multiput(0,2)(1,0){2}{\line(0,1){1}}\epic&
\bpic(2,3)\multiput(0,2)(0,1){2}{\line(1,0){2}}\multiput(0,2)(1,0){3}{\line(0,1){1}}\epic&
\bpic(3,3)\multiput(0,2)(0,1){2}{\line(1,0){3}}\multiput(0,2)(1,0){4}{\line(0,1){1}}\epic&
\bpic(3,3)\multiput(0,2)(0,1){2}{\line(1,0){3}}\multiput(0,1)(1,0){2}{\line(0,1){2}}\put(0,1){\line(1,0){1}}\multiput(2,2)(1,0){2}{\line(0,1){1}}\epic&
\bpic(3,3)\multiput(0,2)(0,1){2}{\line(1,0){3}}\multiput(0,1)(1,0){3}{\line(0,1){2}}\put(0,1){\line(1,0){2}}\put(3,2){\line(0,1){1}}\epic&
\bpic(3,3)\multiput(0,1)(0,1){3}{\line(1,0){3}}\multiput(0,1)(1,0){4}{\line(0,1){2}}\epic&
\bpic(4,3)\multiput(0,2)(0,1){2}{\line(1,0){4}}\multiput(0,1)(1,0){4}{\line(0,1){2}}\put(0,1){\line(1,0){3}}\put(4,2){\line(0,1){1}}\epic&
\bpic(4,3)\multiput(0,2)(0,1){2}{\line(1,0){4}}\multiput(0,1)(1,0){3}{\line(0,1){2}}\put(0,1){\line(1,0){2}}\multiput(3,2)(1,0){2}{\line(0,1){1}}\epic&
\bpic(3,3)\multiput(0,2)(0,1){2}{\line(1,0){3}}\multiput(0,1)(1,0){3}{\line(0,1){2}}\put(0,1){\line(1,0){2}}\put(3,2){\line(0,1){1}}\epic&
\bpic(3,3)\multiput(0,2)(0,1){2}{\line(1,0){3}}\multiput(0,0)(1,0){2}{\line(0,1){3}}\put(0,1){\line(1,0){2}}\put(2,1){\line(0,1){2}}
\put(0,0){\line(1,0){1}}\put(3,2){\line(0,1){1}}\epic&
\bpic(3,3)\multiput(0,0)(0,1){2}{\line(1,0){2}}\multiput(0,0)(1,0){3}{\line(0,1){3}}\multiput(0,2)(0,1){2}{\line(1,0){3}}\put(3,2){\line(0,1){1}}
\epic\\
\rule{0mm}{5mm}\O(\e)_k&&0&1&2&-1&0&1&3&1&3&-2&-1\\
\ea$\caption{\protect\rule{0ex}{2.5ex}Example of an oscillating tableau.\label{ExOTTab}}\end{table}

The enumeration of oscillating tableaux, and their relationship with standard Young tableaux, will now be discussed.
For a partition~$\l$ and a finite subset $T$ of~$\P$ with $|T|=|\l|$, define the corresponding set of \emph{standard Young tableaux} as
\begin{equation}\ba{r@{}l}\SYT(T,\l)\;:=\;\Bigl\{\s\,\Big|\,\:&
\bullet \ \s\mbox{ is an array of entries }\s_{ij}\in T\mbox{ with }(i,j)\in Y(\l)\\[1.5mm]
&\bullet \ \s_{ij}<\s_{i,j\pl1}\mbox{ \ for each }i\in[\ell(\l)], \ j\in[\l_i\mi1]\\[1.5mm]
&\bullet \ \s_{ij}<\s_{i\pl1,j}\mbox{ \ for each }i\in[\ell(\l)\mi1], \ j\in[\l_{i\pl1}]\Bigr\}\,.\ea\end{equation}
For the case $T=[|\l|]$, define
\begin{equation}\SYT(\l)\::=\:\SYT([|\l|],\l)\,.\end{equation}
A standard Young tableau $\s\in\SYT(T,\l)$ is usually shown by writing each $\s_{ij}$ within the
unit square centered at $(i,j)$ in the Young diagram of $\l$.  These entries thus increase along each row
and down each column, and together comprise $T$. It follows immediately that $|\SYT(T,\l)|$ is independent of $T$,
this number usually being denoted as $f^\l$. A well-known formula for $f^\l$
is the \emph{hook-length formula} (see for example~\cite[Corollary~7.21.6]{Sta99}).
An example of an element of $\SYT((4,2,2))$ is
\setlength{\unitlength}{4.2mm}
\begin{equation}\label{ExSYT}\raisebox{-1.27\unitlength}[1.7\unitlength][1.3\unitlength]{
\bpic(2.7,3)\put(0,1.5){\p{l}{\s\;=}}\epic
\bpic(4,3)\multiput(0,0)(1,0){3}{\line(0,1){3}}\multiput(0,0)(0,1){2}{\line(1,0){2}}
\multiput(0,2)(0,1){2}{\line(1,0){4}}\multiput(3,2)(1,0){2}{\line(0,1){1}}
\put(0.5,0.5){\pp{}{5}}\put(0.5,1.5){\pp{}{4}}\put(0.5,2.5){\pp{}{1}}
\put(1.5,0.5){\pp{}{8}}\put(1.5,1.5){\pp{}{6}}\put(1.5,2.5){\pp{}{2}}
\put(2.5,2.5){\pp{}{3}}\put(3.5,2.5){\pp{}{7}}\epic}\;.\end{equation}

It can be seen that each oscillating tableau $\e$ of shape $\l$ with zero deletions, i.e., $\e\in\OT(\l,|\l|)$, can be bijectively associated
with a standard Young tableau $\s\in\SYT(\l)$, where, for $(i,j)\in Y(\l)$ and $k\in[|\l|]$,
$Y(\e_k)=Y(\e_{k\mi1})\cup\{(i,j)\}$ if and only if $\s_{ij}=k$. For example, the
oscillating tableau $\e\in\OT((4,2,2),8)$ which corresponds to $\s$ of~(\ref{ExSYT}) is shown in Table~\ref{SYTOT}.

\setlength{\unitlength}{2.3mm}
\begin{table}[h]\centering\footnotesize$\ba{c||c|c|c|c|c|c|c|c|c}
k&0&1&2&3&4&5&6&7&8\\[1mm]\hline
\rule{0mm}{5mm}\e_k&\ss\emptyset&\ss(1)&\ss(2)&\ss(3)&\ss(3,1)&\ss(3,1,1)&\ss(3,2,1)&\ss(4,2,1)&\ss(4,2,2)\\
&\bpic(0.6,4.3)\put(0.3,2.5){\pp{}{\emptyset}}\epic&
\bpic(1,3)\multiput(0,2)(0,1){2}{\line(1,0){1}}\multiput(0,2)(1,0){2}{\line(0,1){1}}\epic&
\bpic(2,3)\multiput(0,2)(0,1){2}{\line(1,0){2}}\multiput(0,2)(1,0){3}{\line(0,1){1}}\epic&
\bpic(3,3)\multiput(0,2)(0,1){2}{\line(1,0){3}}\multiput(0,2)(1,0){4}{\line(0,1){1}}\epic&
\bpic(3,3)\multiput(0,2)(0,1){2}{\line(1,0){3}}\multiput(0,1)(1,0){2}{\line(0,1){2}}\put(0,1){\line(1,0){1}}\multiput(2,2)(1,0){2}{\line(0,1){1}}\epic&
\bpic(3,3)\multiput(0,2)(0,1){2}{\line(1,0){3}}\multiput(0,0)(1,0){2}{\line(0,1){3}}\multiput(0,0)(0,1){2}{\line(1,0){1}}
\multiput(2,2)(1,0){2}{\line(0,1){1}}\epic&
\bpic(3,3)\multiput(0,2)(0,1){2}{\line(1,0){3}}\multiput(0,0)(1,0){2}{\line(0,1){3}}\put(0,1){\line(1,0){2}}\put(2,1){\line(0,1){2}}
\put(0,0){\line(1,0){1}}\put(3,2){\line(0,1){1}}\epic&
\bpic(4,3)\multiput(0,2)(0,1){2}{\line(1,0){4}}\multiput(0,0)(1,0){2}{\line(0,1){3}}\put(0,1){\line(1,0){2}}\put(2,1){\line(0,1){2}}
\put(0,0){\line(1,0){1}}\multiput(3,2)(1,0){2}{\line(0,1){1}}\epic&
\bpic(4,3)\multiput(0,0)(1,0){3}{\line(0,1){3}}\multiput(0,0)(0,1){2}{\line(1,0){2}}\multiput(0,2)(0,1){2}{\line(1,0){4}}
\multiput(3,2)(1,0){2}{\line(0,1){1}}\epic\\
\rule{0mm}{5mm}\O(\e)_k&&0&1&2&-1&-2&0&3&-1\\
\ea$\caption{\protect\rule{0ex}{2.5ex}The oscillating tableau $\e$ which corresponds to $\s$ of (\ref{ExSYT}).\label{SYTOT}}\end{table}

Now define the set of \emph{ascents}~$\A(\s)$ of any $\s\in\SYT(\l)$ to
be the set of integers $k\in[|\l|\mi1]$ for which $k\pl1$ appears in $\s$ in the same or a higher row
than does~$k$, and denote the number of ascents of~$\s$ as $\AA(\s):=|\A(\s)|$.
For example, $\A(\s)=\{1,2,5,6\}$ for~$\s$ of~(\ref{ExSYT}).
It can be seen that if $\s\in\SYT(\l)$ corresponds to $\e\in\OT(\l,|\l|)$,
then $\A(\s)=\A(<,\e)$, and therefore $\AA(\s)=\AA(<,\e)$.

In~\cite{Sun86,Sun90}, the general case of oscillating tableaux with $n$ deletions is considered, and it
is shown that $\OT(\l,l)$, with $l=|\l|\pl2n$,
can be bijectively mapped to the set of pairs $(M,\s)$ in which $M$ is a matching on a set $T$
(i.e., a decomposition of $T$ as a disjoint union of 2-element subsets), for any $T\subset[l]$
with $|T|=2n=l-|\l|$, and $\s$ is a standard Young tableau in $\SYT([l]\!\setminus\!T,\l)$.
The number of such sets $T$ is $\biggl(\,\ba{c}l\\|\l|\ea\,\biggr)$
and, for any $T$, the number of matchings~$M$ is $(l\mi|\l|\mi1)!!$
and number of standard Young tableaux~$\s$ is~$f^\l$, which gives
\begin{equation}\label{NOT}|\OT(\l,l)|\;=\;\biggl(\,\ba{c}l\\|\l|\ea\,\biggr)\:(l\mi|\l|\mi1)!!\;f^\l\,,\end{equation}
for $(l\mi|\l|)/2\in\N$.

\subsubsection{Progressions and Oscillating Tableaux}
In this section, it will be shown that any progression for a vacancy-osculation set can be simply and naturally associated with
an oscillating tableau.

\begin{theorem}Consider a progression $s=((i_1,j_1),\ldots,(i_l,j_l))$ for a vacancy-oscu\-lation set
$S\in\VOS(a,b,\a,\b,l)$.  Then there exists an oscillating tableau~$\e(s)$
of shape~$\l_{a,b,\a,\b}$ and length $l$ with the profile
\begin{equation}\label{progToOT}\O(\e(s))\;=\;(j_1\mi i_1,\ldots,j_l\mi i_l)\,.\end{equation}
Furthermore, $\e(s)$ has the same number of deletions as the number of osculations of~$S$,
\begin{equation}\label{DChi}\mathcal{D}(\e(s))=\chi(S)\,,\end{equation}
and if $s$ is a canonical ordering, then, for each $k\in[l]$,
$k$ is a deletion of~$\e(s)$ if and only if $(i_k,j_k)$ is an osculation of $S$.\end{theorem}
As indicated in Section~11, any oscillating tableau is uniquely determined by its
profile, so that (\ref{progToOT}) defines a unique oscillating tableau $\e(\s)$ for the progression $s$.
This will be referred to as \emph{the oscillating tableau associated with~$s$}.

\textit{Proof.} \ For each $k\in[0,l]$, denote $\{(i_1,j_1),\ldots,(i_k,j_k)\}$ as $S^k$, and
let $(\a^k,\b^k)\in\BP(a,b)$ be the boundary point pair for which
$S^k\in\VOS(a,b,\a^k,\b^k)$.  Now define the sequence of partitions
\begin{equation}\nu\:=\:(\nu_0,\ldots,\nu_l)\,:=\:(\l_{a,b,\a^0,\b^0},\ldots,\l_{a,b,\a^l,\b^l})\,.\end{equation}
Since~$s$ is a progression for $S$, it follows that $S^{k\mi1}\sim S^k$, so Lemma~6 implies
that $\l_{a,b,\a^{k\mi1},\b^{k\mi1}}$ $\sim\l_{a,b,\a^k,\b^k}$,
i.e., $\nu_{k\mi1}\sim\nu_k$, for each $k\in[l]$.
The boundary point pair for $S^0=\emptyset$ is
$(\a^0,\b^0)=([\min(a,b)],[\min(a,b)])$, so~(\ref{ParSpec}) gives $\nu_0=\l_{a,b,\a^0,\b^0}=\emptyset$, and the
boundary point pair for $S$ is $(\a^l,\b^l)=(\a,\b)$, so $\nu_l=\l_{a,b,\a^l,\b^l}=\l_{a,b,\a,\b}$.
Therefore, $\emptyset=\nu_0\sim\nu_1\sim\ldots\sim\nu_{l-1}\sim\nu_l=\l_{a,b,\a,\b}$,
so $\nu$ is an oscillating tableau of shape $\l_{a,b,\a,\b}$ and length $l$.
Applying~(\ref{DDVOSPar}) and~(\ref{OTPro}) to $S^{k\mi1}$ and $S^k$ for each $k\in[l]$ gives
$\O(\nu)=(j_1\mi i_1,\ldots,j_l\mi i_l)$, which matches the RHS of~(\ref{progToOT}), implying,
since any oscillating tableau is uniquely determined by its profile,
that $\e(s)=\nu$.  The result~(\ref{DChi}) now
follows by applying~(\ref{Smag}) and~(\ref{del}) to $S$ and~$\e(s)$.

Finally, consider a fixed $k\in[l]$.  Then applying~(\ref{Smag}) to $S^k$ and $S^{k\mi1}$ gives
$k=|\e(s)_k|+2\chi(S^k)=1\pl|\e(s)_{k\mi1}|+2\chi(S^{k\mi1})$,
so that~$k$ is a deletion of $\e(s)$, i.e., $|\e(s)_k|=|\e(s)_{k\mi1}|\mi1$,
if and only if $\chi(S^k)=\chi(S^{k\mi1})\pl1$.
Now it follows from Lemma~4 that
a point of $S^{k\mi1}$ is an osculation of $S^{k\mi1}$ and a vacancy of $S^k$, or a vacancy
of $S^{k\mi1}$ and an osculation of $S^k$, if and only if the point is in $S^{k\mi1}\cap B(a,b,i_k,j_k)$,
i.e., $(X(S^{k\mi1})\cap N(S^k))\cup(N(S^{k\mi1})\cap X(S^k))=S^{k\mi1}\cap B(a,b,i_k,j_k)$.
However, if~$s$ is a canonical ordering, then $S^{k\mi1}\cap B(a,b,i_k,j_k)=\emptyset$.
It can now be seen that for this case the following statements are all equivalent:
$k$ is a deletion of $\e(s)$, $\chi(S^k)=\chi(S^{k\mi1})\pl1$, $X(S^k)$ is the disjoint union of $X(S^{k\mi1})$ and $\{s_k\}$,
$s_k\in X(S^k)$, $s_k\in X(S^k)\cap X(S^{k\pl1})\cap\ldots\cap X(S^l)$, $s_k\in X(S)$
(where $s_k=(i_k,j_k)$).\hspace{\fill}$\Box$

For the example of a progression~$s=\{(i_1,j_1),\ldots,(i_{11},j_{11})\}$ of~(\ref{ExProg}) and Figure~\ref{ExProgFig},
$(j_1\mi i_1,\ldots,j_{11}\mi i_{11})=(0,1,2,-1,0,1,3,1,3,-2,-1)$.  This matches the profile
for the example of an oscillating tableau of~(\ref{ExOT}) and Table~\ref{ExOTTab}, which is therefore the oscillating
tableau associated with $s$. Also, by Theorem~10, since $s$ is a canonical ordering, the deletions $8$ and $9$ of
the oscillating tableau correspond to the osculations $(i_8,j_8)=(3,4)$ and $(i_9,j_9)=(2,5)$ of the vacancy-osculation
set.

\subsubsection{Generalized Oscillating Tableaux}
It was seen in the previous section that a vacancy-osculation set can be associated, via
a progression, with an oscillating tableau.
However a progression is not, in general, uniquely determined by its associated oscillating tableau.  In this section, a generalization
of oscillating tableaux will be introduced, and it will be shown that, via the progression induced by
an ordering pair, a vacancy-osculation set can then be uniquely associated with such a generalized oscillating tableau.

For any subset $T$ of $\P$, total strict order $\prec$ on $\Z$, partition~$\l$ and nonnegative integer~$l$,
define the corresponding set of \emph{generalized oscillating tableaux} as
\begin{equation}\label{GOT}\ba{r@{\;}c@{\;}l}\GOT(T,\prec,\l,l)&:=&\Bigl\{(t,\e)=((t_1,\ldots,t_l),\,\e)\in T^l\!\times\OT(\l,l)
\;\Big|\;t_k<t_{k\pl1},\\[2.3mm]
&&\qquad\mbox{or }t_k=t_{k\pl1}\mbox{ and }k\in\A(\prec,\e),\mbox{ \ for each }k\in[l\mi1]\Bigr\}\\[3.4mm]
&=&\Bigl\{(t,\e)\in T^l\!\times\OT(\l,l)\;\Big|\;t_k<t_{k\pl1},\mbox{ or }t_k=t_{k\pl1}\mbox{ and }\\[2mm]
&&\quad\qquad\qquad\qquad\O(\e)_k\prec\O(\e)_{k\pl1},
\mbox{ \ for each }k\in[l\mi1]\Bigr\}\,,\ea\end{equation}
where in the second line, the definition~(\ref{OTAsc}) of $\A(\prec,\e)$ has simply been used.

It follows that, for finite $T$,
\begin{equation}\label{NGOT}
\ba{r@{}l}|\GOT(T,\prec,\l,l)|\;\;&=\,\ds\sum_{\e\in\OT(\l,l)}\left(\;\ba{c}|T|+\AA(\prec,\e)\\[1.3mm]l\ea\;\right)\\[7mm]
&=\;\ds\sum_{i=0}^{l\mi1}\;\mathcal{N}(\l,l,\prec,i)\,
\left(\,\ba{c}|T|+i\\[1.3mm]l\ea\,\right),\ea\end{equation}
where definitions~(\ref{OTAscN}) and~(\ref{Adist}) have been used.

Using the example (\ref{ExOT}) of an oscillating tableau $\e$, and the example~(\ref{ExTSO}) of a total strict order~$\prec$ on~$\Z$,
for which $\A(\prec,\e)$ is given by~(\ref{ExAsc}), it can be seen that an
example of a generalized oscillating tableau in $\GOT([4],\prec,(3,2,2),11)$ is
\begin{equation}\label{ExGOT}
\ba{l}(t,\e)\,=\,\Bigl((1,1,1,2,2,2,2,3,3,4,4),\\[1.2mm]
\hspace{10.2mm}(\emptyset,(1),(2),(3),(3,1),(3,2),(3,3),(4,3),(4,2),(3,2),(3,2,1),(3,2,2))\Bigr).\ea
\end{equation}

\begin{theorem}Consider $a,b\in\P$, $(\a,\b)\in\BP(a,b)$, $l\in\N$, and
an ordering pair $(L,\prec)$.  For each vacancy-osculation set $S\in\VOS(a,b,\a,\b,l)$,
let $s$ be the progression for $S$ induced by $(L,\prec)$ and
$\e(s)$ be the oscillating tableau associated with $s$, and define
\begin{equation}\label{the11}\Theta_{(L,\prec)}(S)\;:=\;\Bigl((L(s_1),\ldots,L(s_l)),\,\e(s)\Bigr)\,.\end{equation}
Then $\Theta_{(L,\prec)}$ is an injective mapping from $\VOS(a,b,\a,\b,l)$ to the set of generalized oscillating tableaux
$\GOT(L([a]\t[b]),\prec,\l_{a,b,\a,\b},l)$, where $L([a]\t[b]):=\{L(i,j)\:|\:(i,j)\in[a]\t[b]\}$.\end{theorem}
\textit{Proof.} \ The definition of $L([a]\t[b])$ immediately gives $(L(s_1),\ldots,L(s_l))\in L([a]\t[b])^l$, while
Theorem~10 gives $\e(s)\in\OT(\l_{a,b,\a,\b},l)$.  Also,~(\ref{IndProg}),~(\ref{OTAsc})
and~(\ref{progToOT}) give
$L(s_k)<L(s_{k\pl1})$, or $L(s_k)=L(s_{k\pl1})$ and $k\in\A(\prec,\e(s))$, for each $k\in[l\mi1]$.
Therefore, by~(\ref{GOT}), $\Theta_{(L,\prec)}(S)\in\GOT(L([a]\t[b]),\prec,\l_{a,b,\a,\b},l)$
for each $S\in\VOS(a,b,\a,\b,l)$.

Now consider $S,S'\in\VOS(a,b,\a,\b,l)$ for which $\Theta_{(L,\prec)}(S)=\Theta_{(L,\prec)}(S')$.
It follows that $L(s_k)=L(s'_k)$ for each $k\in[l]$, where $s$ and $s'$ are the progressions
for $S$ and $S'$ induced by $(L,\prec)$.
The associated oscillating tableaux are also equal, $\e(s)=\e(s')$, so setting
$s_k=(i_k,j_k)$ and $s'_k=(i'_k,j'_k)$,~(\ref{progToOT}) gives
$j_k\mi i_k=j'_k\mi i'_k$ for each $k\in[l]$.
Therefore, by~(\ref{OrdPairEq}), $(i_k,j_k)=(i'_k,j'_k)$ for each $k\in[l]$, giving $s=s'$ and $S=S'$,
and implying that $\Theta_{(L,\prec)}$ is injective.\hspace{\fill}$\Box$

For an ordering pair $(L,\prec)$,
$\Theta_{(L,\prec)}\!:\VOS(a,b,\a,\b,l)\!\rightarrow\!\GOT(L([a]\t[b]),\prec,\l_{a,b,\a,\b},l)$
as defined by~(\ref{the11}) will be referred to as \emph{the mapping induced by $(L,\prec)$},
and, for any vacancy-osculation set $S\in\VOS(a,b,\a,\b,l)$, $\Theta_{(L,\prec)}(S)$ will be referred to as
\emph{the generalized oscillating tableau associated with~$S$ via $(L,\prec)$}.

For the vacancy-osculation set~$S$ of~(\ref{ExS}) and the ordering pair $(L,\prec)$ given by~(\ref{ExTSO})
and~(\ref{ExL}), it has already been seen that
the progression for $S$ induced by $(L,\prec)$ is~$s$ of~(\ref{ExProg}), and that
the oscillating tableau associated with $s$ is that of~(\ref{ExOT}). It can now be seen that
the generalized oscillating tableau associated with~$S$ via $(L,\prec)$ is that of~(\ref{ExGOT}).

\subsubsection{Bijections Between Vacancy-Osculation Sets and Generalized\\Oscillating Tableaux}
It was seen in the previous section that any ordering pair $(L,\prec)$ leads to an injective mapping from
any set $\VOS(a,b,\a,\b,l)$ of vacancy-osculation sets to
the set $\GOT(L([a]\t[b]),\prec,\l_{a,b,\a,\b},l)$ of generalized oscillating tableaux, and it will be seen in
this section that there exist certain ordering pairs for which this mapping is always bijective.

\begin{lemma}Consider $a,b\in\P$, $(\a,\b)\in\BP(a,b)$, $l\in\N$,
an ordering pair $(L,\prec)$, and $(t,\e)\in\GOT(L([a]\t[b]),\prec,\l_{a,b,\a,\b},l)$.
If, for each $k\in[l]$, there
exists $s_k=(i_k,j_k)\in[a]\t[b]$ such that $L(s_k)=t_k$ and $j_k\mi i_k=\O(\e)_k$,
then $s_k$ is unique, the set $S=\{s_1,\ldots,s_l\}$ is in
$\VOS(a,b,\a,\b,l)$, and $(t,\e)$ is the generalized oscillating tableau associated with
$S$ via $(L,\prec)$.
\end{lemma}
\textit{Proof.} \ For each $k\in[l]$, let there exist
$s_k=(i_k,j_k)\in[a]\t[b]$ with $L(s_k)=t_k$ and $j_k\mi i_k=\O(\e)_k$,
Then $s_k$ is unique due to~(\ref{OrdPairEq}).
Also, since $(t,\e)\in\GOT(L([a]\t[b]),\prec,$ $\l_{a,b,\a,\b},l)$, it follows that $t_k<t_{k'}$, or $t_k=t_{k'}$ and
$\O(\e)_k\prec\O(\e)_{k'}$,  for each $k,k'\in[l]$ with $k<k'$,
and hence that $L(s_k)<L(s_{k'})$, or $L(s_k)=L(s_{k'})$ and $j_k\mi i_k\prec j_{k'}\mi i_{k'}$.
Therefore, by~(\ref{OrdPair}), $i_k<i_{k'}$ or $j_k<j_{k'}$, from which it follows,
due to~(\ref{B}), that for each $m\in[0,l\mi1]$,
\begin{equation}\label{lem11}\{s_1,\ldots,s_m\}\cap(B(a,b,i_{m\pl1},j_{m\pl1}\pl1)\cup B(a,b,i_{m\pl1}\pl1,j_{m\pl1}))=
\emptyset\,.\end{equation}
Now consider any $m\in[0,l\mi1]$, and assume that $\{s_1,\ldots,s_m\}\in\VOS(a,b,\a^m,\b^m)$ and
$\l_{a,b,\a^m,\b^m}=\e_m$ for some $(\a^m,\b^m)\in\BP(a,b)$.  It will be shown that, in this case,
$\{s_1,\ldots,s_{m\pl1}\}\in\VOS(a,b,\a^{m\pl1},\b^{m\pl1})$ and
$\l_{a,b,\a^{m\pl1},\b^{m\pl1}}=\e_{m\pl1}$ for some $(\a^{m\pl1},\b^{m\pl1})\in\BP(a,b)$.
Since $\e_m\sim\e_{m\pl1}$ and $\O(\e)_{m\pl1}=j_{m\pl1}\mi i_{m\pl1}$, $j_{m\pl1}\mi i_{m\pl1}$ is a change diagonal
of $\l_{a,b,\a^m,\b^m}$.  Using this and~(\ref{lem11}),
it follows by Lemma~7 that $s_{m\pl1}$ is an addition point of $\{s_1,\ldots,s_m\}$.
Therefore, $\{s_1,\ldots,s_{m\pl1}\}\in\VOS(a,b,\a^{m\pl1},\b^{m\pl1})$ for
some $(\a^{m\pl1},\b^{m\pl1})\in\BP(a,b)$, and Lemma~6 implies that $\l_{a,b,\a^m,\b^m}\sim\l_{a,b,\a^{m\pl1},\b^{m\pl1}}$
and $j_{m\pl1}\mi i_{m\pl1}=\D(\l_{a,b,\a^m,\b^m},\l_{a,b,\a^{m\pl1},\b^{m\pl1}})$, so that
(since any partition $\m$ with $\l_{a,b,\a^m,\b^m}\sim\m$ is uniquely determined by $\D(\l_{a,b,\a^m,\b^m},\m)$)
$\l_{a,b,\a^{m\pl1},\b^{m\pl1}}=\e_{m\pl1}$.

Since $\emptyset\in\VOS(a,b,[\min(a,b)],[\min(a,b)])$ and, from~(\ref{ParSpec}), $\l_{a,b,[\min(a,b)],[\min(a,b)]}=\emptyset$ $=\e_0$,
it follows by repeated application of the result from the previous paragraph starting with $m=0$ that, for each $m\in[0,l]$,
$\{s_1,\ldots,s_m\}\in\VOS(a,b,\a^m,\b^m,m)$ for some $(\a^m,\b^m)\in\BP(a,b)$ with
$\l_{a,b,\a^m,\b^m}=\e_m$. Taking $m=l$ gives $S=\{s_1,\ldots,s_l\}\in\VOS(a,b,\a,\b,l)$.
It can also be seen that $(s_1,\ldots,s_l)$ is the progression for $S$ induced by
$(L,\prec)$, that $\e$ is the oscillating tableau associated with $s$, and that therefore $(t,\e)$
is the generalized oscillating tableau associated with $S$ via $(L,\prec)$.\hspace{\fill}$\Box$

For any $q,i,j\in\Z$, define
\begin{equation}\label{Lq}L_q(i,j)\::=\:\left\{\ba{l}\max(i\pl q,j),\ q\le0\\[1.5mm]
\max(i,j\mi q),\ q\ge0\,,\ea\right.\end{equation}
or, more compactly, $L_q(i,j)=\max(\min(i,i\pl q),\min(j,j\mi q)$.  Usually,
$L_q$ will be regarded as a level function from $\P^2$ to $\P$ for some fixed $q$,
and in this case it can be visualized as
\begin{equation}\left(\ba{c@{\;\;}c@{\;}c}\ss L(1,1)&\ss L(1,2)&\ldots\\[-0.5mm]
\ss L(2,1)&\ss L(2,2)\\[-1mm]
\vdots&&\ddots\ea\right)\;=\;\left\{\ba{l}
\left(\ba{r@{\;}c}\!
{\sss q}\Biggl\{\ba{c@{\;\;}c@{\;\;}c}
\ss1&\ss2&\ss3\\[-2mm]
\vdots&\vdots&\vdots\\[-1.7mm]
\ss1&\ss2&\ss3\ea\\[-1mm]
\ba{c@{\;\;}c@{\;\;}c}
\ss2&\ss2&\ss3\\[-1mm]
\ss3&\ss3&\ss3\ea\\[-3mm]
&\ddots
\ea\right),\quad q\le0\\[12mm]
\left(\ba{c@{\;\,}c@{\;\,}c@{\;}c}
{\ss1}\ldots{\ss1}&\ss2&\ss3\\[-0.9mm]
{\ss2}\ldots{\ss2}&\ss2&\ss3\\[-0.9mm]
\ss\underbrace{{\ss3}\ldots{\ss3}}_q&\ss3&\ss3\\[-5.5mm]
&&&\ddots\ea\right),\quad q\ge0\,.
\ea\right.\end{equation}
It will also be seen that for the case of osculating paths in an $a$ by $b$ rectangle, $q$ is taken to be~$b\mi a$.

It follows from~(\ref{OrdPair}) that for $q\in\Z$ and a $q$-dependent total strict
order $\prec_q$ on $\Z$, $(L_q,\prec_q)$ is an ordering pair if and only if
$\prec_q$ satisfies the property that, for $z,z'\in\Z$,
\begin{equation}\label{TSOQ1}\mbox{If \ }z<z'\le q\mbox{ \ or \ }z>z'\ge q,\mbox{ \ then \ }z\prec_q z'\,,\end{equation}
i.e., \ $\ldots\prec_q\:q\mi2\:\prec_q\:q\mi1\:\prec_q\:q$ \ and \ $\ldots\prec_q\:q\pl2\:\prec_q\:q\pl1\:\prec_q\:q$.
Any total strict order $\prec_q$ on $\Z$ which satisfies~(\ref{TSOQ1}) will be referred to as a \emph{$q$-order}.

A specific example of a $q$-order is that defined, for any $z,z'\in\Z$, by
\begin{equation}\label{TSOQ2}z\prec_q z'\mbox{ \ if and only if \ }|z\mi q|>|z'\mi q|\mbox{ \ or \ }z\mi q=q\mi z'<0\,,\end{equation}
i.e., \ $\ldots\prec_q\:q\mi2\:\prec_q\:q\pl2\:\prec_q\:q\mi1\:\prec_q\:q\pl1\:\prec_q\:q$.

Any ordering pair $(L_q,\prec_q)$ in which the level function
$L_q$ is given by~(\ref{Lq}) (and~$\prec_q$ is therefore a $q$-order) will be referred
to as a \emph{complete ordering pair}.
It can be seen that the example of an ordering pair $(L,\prec)$ given by~(\ref{ExTSO}) and~(\ref{ExL})
is the complete ordering pair $(L_2,\prec_2)$ with $\prec_2$ given by~(\ref{TSOQ2}).

Now define, for any $q,u,z\in\Z$,
\begin{equation}\label{Gq}
G_q(u,z)\;:=\;\left\{\ba{ll}
(u\mi z,u)\,,&q\le\min(0,z)\\[2mm]
(u\mi q,u\mi q\pl z)\,,&z\le q\le0\\[2mm]
(u\pl q\mi z,u\pl q)\,,&0\le q\le z\\[2mm]
(u,u\pl z)\,,&q\ge\max(0,z)\,.\ea\right.
\end{equation}

It can be checked straightforwardly that for $q,i,j,u,z\in\Z$,
\begin{equation}\label{LqGq}(L_q(i,j),j\mi i)\,=\,(u,z)\mbox{ \ if and only if \ }(i,j)\,=\,G_q(u,z)\,.\end{equation}

Finally, in preparation for Theorem~13 define,
for any $q\in\Z$, $T\subset\P$, total strict order $\prec$ on~$\Z$, $\l\in\Par$, $l\in\N$ and $(t,\e)\in\GOT(T,\prec,\l,l)$,
\begin{equation}\label{Phi}
\Phi_q(t,\e)\;:=\;\{G_q(t_1,\O(\e)_1),\ldots,G_q(t_l,\O(\e)_l)\}\,.
\end{equation}

\begin{theorem}Consider $a,b\in\P$, $(\a,\b)\in\BP(a,b)$, and $l\in\N$.  Then the
mapping $\Theta_{(L_{b\mi a},\prec_{b\mi a})}$ induced by any complete ordering pair $(L_{b\mi a},\prec_{b\mi a})$
is a bijection from $\VOS(a,b,\a,\b,l)$ to $\GOT([\min(a,b)],\prec_{b\mi a},\l_{a,b,\a,\b},l)$, and the inverse
mapping is $\Phi_{b\mi a}$.
\end{theorem}
\textit{Proof.} \ The mapping $\Theta_{(L_{b\mi a},\prec_{b\mi a})}$ induced by a complete ordering pair
$(L_{b\mi a},\prec_{b\mi a})$ is given by~(\ref{the11}).  By Theorem~11, it
is an injective mapping from $\VOS(a,b,\a,\b,l)$ to $\GOT([\min(a,b)],\prec_{b\mi a},\l_{a,b,\a,\b},l)$,
with $L_{b\mi a}([a]\t[b])\!=\![\min(a,b)]$ following from~(\ref{Lq}).  Thus it only
remains to be shown that it is a surjective mapping and that its inverse is~$\Phi_{b\mi a}$.
This will now be done explicitly for the case $a\le b$, i.e., $\min(a,b)=a$.  The case $b\le a$
follows from the case $a\le b$ by reflection of the lattice in the main diagonal using
the first bijection of~(\ref{OPRT}).

Consider any $(t,\e)\in\GOT([a],\prec_{b\mi a},\l_{a,b,\a,\b},l)$, and abbreviate $\O(\e)_k$ as
$\O_k$ for each $k\in[l]$.  By~(\ref{GOT}), $1\le t_1\le\ldots\le t_l\le a$, and for $n,n'\in[l]$,
\begin{equation}\label{the131}
\mbox{If \ }n<n',\mbox{ \ and \ }\O_{n'}\prec_{b\mi a}\O_n\mbox{ \ or \ }\O_n=\O_{n'},
\mbox{ \ then \ }t_n<t_{n'}\,.\end{equation}
Now, for each $k\in[l]$, define $s_k=(i_k,j_k)$ by
\begin{equation}\label{the132}s_k=G_{b\mi a}(t_k,\O_k)\;=\;\left\{\ba{l}
(t_k,t_k\pl\O_k)\,,\quad\O_k\le b\mi a\\[2mm]
(t_k\mi\O_k\pl b\mi a,t_k\pl b\mi a)\,,\quad\O_k\ge b\mi a\,,\ea\right.
\end{equation}
the explicit form here following immediately from the definition~(\ref{Gq})
of $G_q$, with $q=b\mi a\ge0$.
It will first be shown that $s_k\in [a]\t[b]$ for each $k\in[l]$.
It can be seen straightforwardly that $1\le t_k\le a\le b$ implies that $i_k\ge1$ for
$\O_k\le b\mi a$, $j_k\ge1$ for $\O_k\ge0$, $i_k\le a$ and $j_k\le b$.  Therefore it only remains
to be checked that $i_k\ge1$ for $\O_k>b\mi a$, and $j_k\ge1$ for $\O_k<0$.
If $\O_k>b\mi a$, then, by applying the second property of~(\ref{SDD}) $\O_k\mi b\pl a$ times
in succession, it follows that there exist integers $1\le n_0<n_1<\ldots<n_{\O_k\mi b\pl a}=k$
such that $\O_{n_m}=b\mi a\pl m$ for each $m\in[0,\O_k\mi b\pl a]$.  Since, by~(\ref{TSOQ1}),
$\O_k\prec_{b\mi a}\O_k\mi1\prec_{b\mi a}\ldots\prec_{b\mi a}b\mi a$, it follows using~(\ref{the131}) that
$t_{n_m}\ge m\pl1$ for each $m\in[0,\O_k\mi b\pl a]$ and so, taking $m=\O_k\mi b\pl a$,
that $t_k\ge\O_k\mi b\pl a\pl1$, i.e., by~(\ref{the132}), that $i_k\ge1$ as required.
Similarly, if $\O_k<0$, then, by applying the first property of~(\ref{SDD}) $-\O_k$ times
in succession, it follows that there exist integers $1\le n_0<n_1<\ldots<n_{-\O_k}=k$
such that $\O_{n_m}=-m$ for each $m\in[0,-\O_k]$.  Since, by~(\ref{TSOQ1}),
$\O_k\prec_{b\mi a}\O_k\pl1\prec_{b\mi a}\ldots\prec_{b\mi a}0$, it follows using~(\ref{the131}) that
$t_{n_m}\ge m\pl1$ for each $m\in[0,-\O_k]$ and so, taking $m=-\O_k$,
that $t_k\ge-\O_k\pl1$, i.e., by~(\ref{the132}), that $j_k\ge1$ as required.

Having now verified that $s_k\in [a]\times[b]$, and observing from~(\ref{LqGq}) and~(\ref{the132})
that $L_{b\mi a}(s_k)=t_k$ and $j_k\mi i_k=\O_k$ for each $k\in[l]$, it follows by Lemma~12
that $\{s_1,\ldots,s_l\}\in\VOS(a,b,\a,\b,l)$ and that
$\Theta_{(L_{b\mi a},\prec_{b\mi a})}(\{s_1,\ldots,s_l\})=(t,\e)$, so that
$\Theta_{(L_{b\mi a},\prec_{b\mi a})}$ is surjective, with the inverse mapping being~$\Phi_{b\mi a}$.\hspace{\fill}$\Box$

Returning to the running example, if it were not already known that
$(t,\e)$ of~(\ref{ExGOT}) is the generalized oscillating tableau associated with the
vacancy-osculation matrix~$S$ of~(\ref{ExS}) via the complete ordering pair $(L,\prec)$
of~(\ref{ExTSO}) and~(\ref{ExL}), then $S$ could be obtained
by simply evaluating the profile for $\e$, as
done in Table~\ref{ExOTTab}, and then applying~$\Phi_2$ to give
$S=\{G_2(t_1,\O(\e)_1),\ldots,G_2(t_{11},\O(\e)_{11})\}$.

\begin{corollary}Consider $a,b\in\P$, $(\a,\b)\in\BP(a,b)$, $l\in\N$, and a $(b\mi a)$-order~$\prec_{b\mi a}$.
Then
\begin{equation}\label{NOP1}
|\OP(a,b,\a,\b,l)|\;=\,\sum_{\e\in\OT(\l_{a,b,\a,\b},l)}\left(\;\ba{c}\min(a,b)+\AA(\prec_{b\mi a},\e)\\[1.3mm]l\ea\;\right).\end{equation}
\end{corollary}
\textit{Proof.} \ This result follows immediately from Theorem~13,~(\ref{NGOT}),
and the bijection between $\OP(a,b,\a,\b,l)$ and $\VOS(a,b,\a,\b,l)$.\hspace{\fill}$\Box$

It was seen in Theorem~13 that there is a bijection between each set of vacancy-osculation sets $\VOS(a,b,\a,\b,l)$
and a set of generalized oscillating tableaux, and it will now be seen in Theorem~15 that there is a bijection between each
set of generalized oscillating tableaux $\GOT([n],\prec_q,\l,l)$ with $n\ge L_q(\ell(\l),\l_1)$ and a set of
vacancy-osculation sets.

\begin{theorem}Consider $n\in\P$, $q\in\Z$, a $q$-order~$\prec_q$, $\l\in\Par$, and $l\in\N$.
Then for $n<L_q(\ell(\l),\l_1)$, $\GOT([n],\prec_q,\l,l)$ is empty, while for
$n\ge L_q(\ell(\l),\l_1)$, the mapping $\Phi_q$ is a bijection from
$\GOT([n],\prec_q,\l,l)$ to $\VOS(a(n,q),b(n,q),\a_{a(n,q),b(n,q),\l},$ $\b_{a(n,q),b(n,q),\l},l)$,
where $a(n,q)=\max(n,n\mi q)$ and $b(n,q)=\max(n,n\pl q)$, and the inverse mapping is~$\Theta_{(L_q,\prec_q)}$.
\end{theorem}
Here, $L_q$ is defined in~(\ref{Lq}), and
$(\a_{a(n,q),b(n,q),\l},\b_{a(n,q),b(n,q),\l})$ is the boundary point pair in $\BP(a(n,q),b(n,q))$ which corresponds to~$\l$,
as given explicitly by~(\ref{BPInv}).  It can be checked straightforwardly that $\l\in\Par(a(n,q),b(n,q))$, for $n\ge L_q(\ell(\l),\l_1)$,
so that $(\a_{a(n,q),b(n,q),\l},\b_{a(n,q),b(n,q),\l})$ is well-defined.
Note also that $(a(n,q),b(n,q))$ is simply $(n\mi q,n)$ for $q\le0$, or $(n,n\pl q)$ for
$q\ge0$.  It can be seen that if~$q$ and~$\l$ are regarded as fixed and $n$ as variable,
then $a(n,q)$, $b(n,q)$, $\a_{a(n,q),b(n,q),\l}$ and $\b_{a(n,q),b(n,q),\l}$
represent a sequence of rectangles and boundary points with a fixed difference~$q$ between the length and width of the rectangle,
and a fixed corresponding partition~$\l$.  This occurs, for example, in the first three cases of
Table~\ref{ASMPar} if $m$ is regarded as fixed and $n$ as variable.

\textit{Proof.} \  For $n\ge L_q(\ell(\l),\l_1)$, the result follows immediately
from Theorem~13. For $n<L_q(\ell(\l),\l_1)$, the case $\l\ne\emptyset$ and $0\le q\le\l_1\mi1$ will now be
considered in detail, all other cases being similar.  Assume that there exists $(t,\e)\in\GOT([n],\prec_q,\l,l)$.
It can be seen, using~(\ref{OTProp}) and the first two properties of~(\ref{SDD}),
that there exist integers $1\le k_0<k_1<\ldots<k_{\l_1\mi q\mi1}\le l$
and $1\le m_0<m_1<\ldots<m_{\ell(\l)\mi1}\le l$ such that
$\O(\e)_{k_i}=q\pl i$, for each $i\in[0,\l_1\mi q\mi1]$, and
$\O(\e)_{m_i}=-i$, for each $i\in[0,\ell(\l)\mi1]$.
Since, by~(\ref{TSOQ1}), $\l_1\mi1\prec_q\ldots\prec_q q\pl1\prec_q q$, it follows that
$t_{k_i}\ge i\pl1$ for each $i\in[0,\l_1\mi q\mi1]$, and so that $t_{k_{\l_1\mi q\mi1}}\ge\l_1\mi q$.
Also, $-\ell(\l)\pl1\prec_q\ldots\prec_q\mi1\prec_q0$, from which it follows that
$t_{m_i}\ge i\pl1$ for each $i\in[0,\ell(\l)\mi1]$, and so that $t_{m_{\ell(\l)\mi1}}\ge\ell(\l)$.
But $t_{k_{\l_1\mi q\mi1}}\le n$ and $t_{m_{\ell(\l)\mi1}}\le n$ must both hold, giving
$n\ge\max(\ell(\l),\l_1\mi q)$, which contradicts $n<L_q(\ell(\l),\l_1)=\max(\ell(\l),\l_1\mi q)$
(since $q\ge0$), and implies that $\GOT([n],\prec_q,\l,l)$ is empty.\hspace{\fill}$\Box$

Theorem~15 and~(\ref{NGOT}) give
\begin{equation}\label{NOP2}\ba{l}
\ds\sum_{i=0}^{l\mi1}\;\mathcal{N}(\l,l,\prec_q,i)\left(\,\ba{c}n+i\\[1.3mm]l\ea\,\right)\,=\\[6mm]
\hspace{12.8mm}\left\{\ba{l}0,\;\;n<L_q(\ell(\l),\l_1)\\[2mm]
|\OP(a(n,q),b(n,q),\a_{a(n,q),b(n,q),\l},\b_{a(n,q),b(n,q),\l},l)|,\;n\ge L_q(\ell(\l),\l_1).\ea\right.\ea
\end{equation}
Since the RHS of~(\ref{NOP2}) is independent of the $q$-order $\prec_q$, it
follows that the distribution of the number of ascents $\AA$ on any fixed $\OP(\l,l)$ is
the same with respect to any $q$-order, i.e., if~$\prec_q$
and~$\prec_q'$ both satisfy~(\ref{TSOQ1}) for $q\in\Z$, then
\begin{equation}\mathcal{N}(\l,l,\prec_q,i)\;=\;\mathcal{N}(\l,l,\prec_q',i)\,,\end{equation}
for any $\l\in\Par$, $l\in\N$ and $i\in[0,l\mi1]$.

It also follows from~(\ref{NOP2}), the first equation of~(\ref{Parmag})
and the first bijection of~(\ref{OPRT}) that
\begin{equation}\mathcal{N}(\l,l,\prec_q,i)\;=\;\mathcal{N}(\l^t,l,\prec_{-q},i)\,,\end{equation}
for any $\l\in\Par$, $l\in\N$, $q\in\Z$, $q$-order $\prec_q$, $-q$-order $\prec_{-q}$ \ and $i\in[0,l\mi1]$.

\subsubsection{Summary}
The demonstration that there exist bijections
between tuples of osculating paths and generalized oscillating tableaux
constitutes one of the main components of this paper. However, since the details of these bijections
appear across many of the previous sections, they will be restated in this section as a summary.

Starting with mappings from path tuples to generalized oscillating tableaux,
consider $a,b\in\P$, $l\in\N$, and a boundary pair $(\a,\b)\in\BP(a,b)$, with $\BP(a,b)$ defined
in~(\ref{BP}).  Then, for any total strict order $\prec_{b\mi a}$ on $\Z$ which
satisfies~(\ref{TSOQ1}) for $q=b\mi a$, there is a bijection from the set $\OP(a,b,\a,\b,l)$ of path tuples,
as defined in~(\ref{OPl}) using~(\ref{OP}) and~(\ref{OPToVOS}),
to the set $\GOT([\min(a,b)],\prec_{b\mi a},\l_{a,b,\a,\b},l)$ of generalized oscillating tableau,
as defined in~(\ref{GOT}) using~(\ref{DD}),~(\ref{ToP}),~(\ref{OT}) and~(\ref{OTPro}).
For any path tuple $P\in\OP(a,b,\a,\b,l)$,
the corresponding generalized oscillating tableau $(t,\e)\in\GOT([\min(a,b)],\prec_{b\mi a},\l_{a,b,\a,\b},l)$
can be obtained as follows.\vspace{-2mm}
\begin{itemize}
\item Find the set $S$ of $l$ points of $[a]\t[b]$ which are vacancies or osculations of~$P$, i.e., points through
which pass zero or two paths of~$P$ respectively.\vspace{-2mm}
\item Order the points of $S$ as the unique tuple $s=((i_1,j_1),\ldots,(i_l,j_l))$
which satisfies~(\ref{IndProg}), with $L$ taken to be $L_{b\mi a}$, as defined in~(\ref{Lq}),
and $\prec$ taken to be~$\prec_{b\mi a}$.\vspace{-2mm}
\item Obtain the unique oscillating tableau $\e(s)$ whose profile, as defined in~(\ref{OTProSeq}), is given by~(\ref{progToOT}).
This step can be performed by simply
starting with $\e(s)_0=\emptyset$, and then obtaining the Young diagram of each successive partition $\e(s)_k$ from
that of $\e(s)_{k\mi1}$ by adding or deleting (only one or the other being possible) a square with content $j_k\mi i_k$,
for each $k\in[l]$.\vspace{-2mm}
\item The generalized oscillating tableau $(t,\e)$ is then $((L_{b\mi a}(i_1,j_1),\ldots,L_{b\mi a}(i_l,j_l)),$ $\e(s))$.
\end{itemize}
\vspace{-2mm}
In this mapping, for each $k\in[l]$, $k$ is a deletion of~$\e(s)$, i.e., a value for which $|\e(s)_k|=|\e(s)_{k\mi1}|\mi1$,
if and only if $(i_k,j_k)$ is an osculation of $P$.

For the example of Figure~\ref{Ex},~$S$,~$s$,~$\e(s)$ and $(t,\e)$ are given in~(\ref{ExS}),~(\ref{ExProg}),~(\ref{ExOT})
and~(\ref{ExGOT}) respectively.

Proceeding to mappings from generalized oscillating tableaux to path tuples, consider
$n\in\P$, $q\in\Z$, a total strict order~$\prec_q$ on $\Z$ which satisfies~(\ref{TSOQ1}), $l\in\N$, and a partition $\l$,
with $n\ge L_q(\ell(\l),\l_1)$, where $L_q$ is defined in~(\ref{Lq}).
Then there is a bijection from the
set $\GOT([n],\prec_q,\l,l)$ of generalized oscillating tableaux to the set
$\OP(a(n,q),b(n,q),\a_{a(n,q),b(n,q),\l},\b_{a(n,q),b(n,q),\l},l)$ of path tuples,
where $a(n,q)=\max(n,n\mi q)$, $b(n,q)=\max(n,n\pl q)$, and $(\a_{a(n,q),b(n,q),\l},\b_{a(n,q),b(n,q),\l})$ is defined
using~(\ref{BPInv}).  For any generalized oscillating tableau $(t,\e)\in\GOT([n],\prec_q,\l,l)$, the corresponding
path tuple $P\in\OP(a(n,q),b(n,q),\a_{a(n,q),b(n,q),\l},\b_{a(n,q),b(n,q),\l},l)$ can be found as follows.\vspace{-2mm}
\vspace{-2mm}
\begin{itemize}
\item Obtain the profile $\O(\e)$, as defined in~(\ref{OTProSeq})
using~(\ref{DD}) and~(\ref{OTPro}), of the oscillating tableau $\e$.\vspace{-2mm}
\item Obtain the vacancy-osculation set $S=\Phi_q(t,\e)$, where $\Phi_q$ is defined in~(\ref{Phi})
using~(\ref{Gq}).\vspace{-2mm}
\item Use (\ref{VOMToHV}) to obtain the edge matrix pair $(H,V)$ corresponding to $S$.\vspace{-2mm}
\item The path tuple is then the unique $P\in\OP(a(n,q),b(n,q))$ which satisfies~(\ref{HV}).
\end{itemize}
\vspace{-2mm}
For the example of~(\ref{ExGOT}), $\O(\e)$, $S$, $(H,V)$ and $P$ are given in Table~\ref{ExOTTab},~(\ref{ExS}),~(\ref{ExHV})
and Figure~\ref{Ex} respectively.

For two sets related bijectively by one of mappings described here, the other mapping is the inverse mapping, where the same total strict
order on $\Z$ is used throughout.

\subsubsection{Lexicographic and Antilexicographic Ordering Pairs}
In this section, the lexicographic and antilexicographic
ordering pairs $(\L,<)$ and $(\LL,>)$, with $\L$ and $\LL$ defined in~(\ref{lex}) and~(\ref{alex}), will be considered
further, and it will be seen that in certain cases these lead to alternative bijections to those described in the
previous section.

By Theorem~11, for any $a,b\in\P$, $(\a,\b)\in\BP(a,b)$, and $l\in\N$,
the mappings induced by $(\L,<)$ or $(\LL,>)$ are injective from $\VOS(a,b,\a,\b,l)$ to $\GOT([a],<,\l_{a,b,\a,\b},l)$
or $\GOT([b],>,\l_{a,b,\a,\b},l)$ respectively.
These mappings are not always bijective,
but it will be seen in Theorem~16 that $\Theta_{(\L,<)}$ is bijective if $\a=[a]$ and
that $\Theta_{(\LL,>)}$ is bijective if $\b=[b]$.
Note that $\a=[a]$ implies that $a\le b$, and that $\b=[b]$ implies that $b\le a$.
The partitions for these cases are given in~(\ref{ParSpec}), and it can be seen that
\begin{equation}\label{Ext}\ba{l@{\;\;}l}\bullet&\a=[a]\mbox{ \ if and only if \ }(\l_{a,b,\a,\b})_1\le b\mi a\\[2.4mm]
\bullet&\b=[b]\mbox{ \ if and only if \ }\ell(\l_{a,b,\a,\b})\le a\mi b\,.\ea\end{equation}

\begin{theorem}Consider $a,b\in\P$, $(\a,\b)\in\BP(a,b)$, and $l\in\N$. If $\a=[a]$, then the
mapping $\Theta_{(\L,<)}$ induced by the lexicographic ordering pair $(\L,<)$
is a bijection from $\VOS(a,b,[a],\b,l)$ to $\GOT([a],<,\l_{a,b,[a],\b},l)$, and
if $\b=[b]$, then the
mapping $\Theta_{(\LL,>)}$ induced by the antilexicographic ordering pair $(\LL,>)$
is a bijection from $\VOS(a,b,\a,[b],l)$ to $\GOT([b],>,\l_{a,b,\a,[b]},l)$.
The inverse mappings for these cases are respectively
\begin{equation}
\hat{\Psi}(t,\e)\;=\;\{(t_1,t_1\pl\O(\e)_1),\ldots,(t_l,t_l\pl\O(\e)_l)\}\,,
\end{equation}
for each $(t,\e)\in\GOT([a],<,\l_{a,b,[a],\b},l)$,
and
\begin{equation}
\check{\Psi}(t,\e)\;=\;\{(t_1\mi\O(\e)_1,t_1),\ldots,(t_l\mi\O(\e)_l,t_l)\}\,,
\end{equation}
for each $(t,\e)\in\GOT([b],>,\l_{a,b,\a,[b]},l)$.
\end{theorem}
\textit{Proof.} \ The proof is similar to that of Theorem~13. Due to Theorem~11, it only needs to
be shown that $\Theta_{(\L,<)}$ and $\Theta_{(\LL,>)}$ are surjective, with inverses $\hat{\Psi}$
and $\check{\Psi}$.
This will now be done explicitly for the case $\a=[a]$.  The case $\b=[b]$
follows from the case $\a=[a]$ by reflection of the lattice in the main diagonal using the first bijection
of~(\ref{OPRT}).

Consider any $(t,\e)\in\GOT([a],<,\l_{a,b,[a],\b},l)$, and again abbreviate $\O(\e)_k$ as
$\O_k$ for each $k\in[l]$.  By~(\ref{GOT}), $1\le t_1\le\ldots\le t_l\le a$, and for $n,n'\in[l]$,
\begin{equation}\label{the151}
\mbox{If \ }n<n'\mbox{ \ and \ }\O_n\ge\O_{n'},\mbox{ \ then \ }t_n<t_{n'}\,.\end{equation}
For each $k\in[l]$, define $s_k=(i_k,j_k)$ by
\begin{equation}\label{the152}s_k=(t_k,t_k\pl\O_k)\,.
\end{equation}
It will now be shown that $s_k\in [a]\t[b]$ for each $k\in[l]$.
It can be seen immediately that $i_k\in[a]$, $j_k\ge1$ for $\O_k\ge0$,
and $j_k\le b$ for $\O_k\le b\mi a$, so it only remains
to be checked that $j_k\ge1$ for $\O_k<0$, and $j_k\le b$ for $\O_k>b\mi a$.
If $\O_k<0$, then, by applying the first property of~(\ref{SDD}) $-\O_k$ times
in succession, it follows that there exist integers $1\le n_0<n_1<\ldots<n_{-\O_k}=k$
such that $\O_{n_m}=-m$ for each $m\in[0,-\O_k]$.  It then follows using~(\ref{the151}) that
$t_{n_m}\ge m\pl1$ for each $m\in[0,-\O_k]$ and so, taking $m=-\O_k$,
that $t_k\ge-\O_k\pl1$, i.e., by~(\ref{the152}), that $j_k\ge1$ as required.
Similarly, if $\O_k>b\mi a$, then abbreviating $(\l_{a,b,[a],\b})_1=\b_a\mi a$ as $\l_1$, noting that
$\l_1\le b\mi a$, and applying the fourth property of~(\ref{SDD}) $\O_k\mi\l_1$ times
in succession, it follows that there exist integers $k=n_{\O_k\mi\l_1}<\ldots<n_1<n_0\le l$
such that $\O_{n_m}=m\pl\l_1$ for each $m\in[0,\O_k\mi\l_1]$.  It then follows using~(\ref{the151}) and $t_l\le a$ that
$t_{n_m}\le a\mi m$ for each $m\in[0,\O_k\mi\l_1]$ and so, taking $m=\O_k\mi\l_1$,
that $t_k\le a\mi\O_k\pl\l_1$, i.e., by~(\ref{the152}), that $j_k\le a\pl\l_1\le b$ as required.

Having now verified that $s_k\in [a]\times[b]$, and observing from~(\ref{lex}) and~(\ref{the152})
that $\L(s_k)=t_k$ and $j_k\mi i_k=\O_k$ for each $k\in[l]$, it follows by Lemma~12
that $\{s_1,\ldots,s_l\}\in\VOS(a,b,[a],\b,l)$ and that
$\Theta_{(\L,<)}(\{s_1,\ldots,s_l\})=(t,\e)$, so that
$\Theta_{(\L,<)}$ is surjective, with the inverse mapping~$\hat{\Psi}$.\hspace{\fill}$\Box$

\begin{corollary}Consider $a,b\in\P$, $(\a,\b)\in\BP(a,b)$, and $l\in\N$. If $\a=[a]$, then
\begin{equation}\label{NOP3}
|\OP(a,b,[a],\b,l)|\;=\,\sum_{\e\in\OT(\l_{a,b,[a],\b},l)}\left(\,\ba{c}a+\AA(<,\e)\\[1.3mm]l\ea\,\right),\end{equation}
and if $\b=[b]$, then
\begin{equation}\label{NOP4}
|\OP(a,b,\a,[b],l)|\;=\,\sum_{\e\in\OT(\l_{a,b,\a,[b]},l)}\left(\,\ba{c}b+\AA(>,\e)\\[1.3mm]l\ea\,\right).\end{equation}
\end{corollary}
\textit{Proof.} \ This follows immediately from Theorem~16,~(\ref{NGOT}),
and the bijection between $\OP(a,b,\a,\b,l)$ and $\VOS(a,b,\a,\b,l)$.\hspace{\fill}$\Box$

It follows from~(\ref{NOP3}) that $|\OP(a,b,[a],\b,l)|$ is independent of $b$.
This can also be seen by examining the paths directly
since, for any $P\in\OP(a,b,[a],\b,l)$ and $k\in[a]$, path $P_k$ must end with a segment from $(k,\b_a)$ to $(k,b)$.
Similarly, $|\OP(a,b,\a,[b],l)|$ is independent of $a$, as can be seen from~(\ref{NOP4}) or from the
fact that, for any $P\in\OP(a,b,\a,[b],l)$ and $k\in[b]$, path $P_k$ starts with a segment from $(a,k)$ to $(\a_b,k)$.

It also follows from~(\ref{NOP1}),~(\ref{Ext}),~(\ref{NOP3}) and~(\ref{NOP4}) that
the distribution of $\AA$ on $\OP(\l,l)$ is the same with respect to $<$ or
any $q$-order~$\prec_q$ with $q\ge\l_1$, and the same with respect to $>$ or
any $q$-order~$\prec_q$ with $q\le-\ell(\l)$,
\begin{equation}\label{NEq}\ba{l}\mathcal{N}(\l,l,<,i)\;=\;\mathcal{N}(\l,l,\prec_q,i)\mbox{ \ for any \ }q\ge\l_1\\[2.5mm]
\mathcal{N}(\l,l,>,i)\;=\;\mathcal{N}(\l,l,\prec_q,i)\mbox{ \ for any \ }q\le-\ell(\l),\ea\end{equation}
with any $\l\in\Par$, $l\in\N$ and $i\in[0,l\mi1]$.

\subsubsection{Further Example}
In this section, a further example, that of enumerating $\OP(n,n,[n],[n],6)$ for any $n\in\P$, will be considered.

As seen in Section~4, the path tuples of $\OP(n,n,[n],[n])$ correspond to $n\times n$ standard alternating sign matrices,
and, as seen in Table~\ref{ASMPar}, $\l_{n,n,[n],[n]}=\emptyset$.
By~(\ref{Pmag}), any path tuple of $\OP(n,n,[n],[n])$ has equal numbers of vacancies and osculations,
so $\OP(n,n,[n],[n],6)$ is the set of all path tuples of $\OP(n,n,[n],[n])$ which have~3 osculations and~3 vacancies.

In this case, (\ref{NOP1}), (\ref{NOP3}) or (\ref{NOP4}) can be used, giving
\begin{equation}\label{ExNOP}\ba{l}\ds|\OP(n,n,[n],[n],6)|\;=\,\sum_{\e\in\OT(\emptyset,6)}\!\left(\ba{c}n+\AA(\prec_0,\e)\\[1.3mm]6\ea\right)\\[9mm]
\hspace{37.5mm}\ds=\,\sum_{\e\in\OT(\emptyset,6)}\!\left(\ba{c}n+\AA(<,\e)\\[1.3mm]6\ea\right)
\;=\,\sum_{\e\in\OT(\emptyset,6)}\!\left(\ba{c}n+\AA(>,\e)\\[1.3mm]6\ea\right),\ea\end{equation}
where $\prec_0$ is any $0$-order.

By~(\ref{NOT}), $|\OT(\emptyset,6)|=5!!=15$.  These~15 oscillating tableaux are shown in Table~\ref{ASMEx}.
In this table, the Young diagrams of each partition of $\e$,
the terms of the profile of~$\e$, with $\O(\e)_k$ abbreviated
as~$\O_k$ for each $k\in[6]$, and the sets of ascents $\A(\prec_0,\e)$, $\A(<,\e)$ and $\A(>,\e)$ are shown for each $\e\in\OT(\emptyset,6)$,
where $\prec_0$ is taken to be the specific $0$-order of~(\ref{TSOQ2}),
\begin{equation}z\prec_0 z'\mbox{ \ if and only if \ }|z|>|z'|\mbox{ \ or \ }z=-z'<0\,,\end{equation}
for $z,z'\in\Z$.

\vspace{8mm}
\setlength{\unitlength}{2.3mm}
\nc{\pe}{\raisebox{-0.04\unitlength}[2\unitlength]{\bpic(0.4,1)\put(0.2,0.5){\p{}{\emptyset}}\epic}}
\nc{\ps}{\raisebox{-0.04\unitlength}[2\unitlength]{\bpic(1,1)\multiput(0,0)(1,0){2}{\line(0,1){1}}\multiput(0,0)(0,1){2}{\line(1,0){1}}\epic}}
\nc{\pd}{\raisebox{-0.04\unitlength}[2\unitlength]{\bpic(2,1)\multiput(0,0)(1,0){3}{\line(0,1){1}}\multiput(0,0)(0,1){2}{\line(1,0){2}}\epic}}
\nc{\pt}{\raisebox{-0.04\unitlength}[2\unitlength]{\bpic(3,1)\multiput(0,0)(1,0){4}{\line(0,1){1}}\multiput(0,0)(0,1){2}{\line(1,0){3}}\epic}}
\nc{\peu}{\raisebox{-1\unitlength}[2\unitlength]{\bpic(0.4,2)\put(0.2,1.5){\p{}{\emptyset}}\epic}}
\nc{\psu}{\raisebox{-1\unitlength}[2\unitlength]{\bpic(1,2)\multiput(0,1)(1,0){2}{\line(0,1){1}}\multiput(0,1)(0,1){2}{\line(1,0){1}}\epic}}
\nc{\pdu}{\raisebox{-1\unitlength}[2\unitlength]{\bpic(2,2)\multiput(0,1)(1,0){3}{\line(0,1){1}}\multiput(0,1)(0,1){2}{\line(1,0){2}}\epic}}
\nc{\ptu}{\raisebox{-1\unitlength}[2\unitlength]{\bpic(3,2)\multiput(0,1)(1,0){4}{\line(0,1){1}}\multiput(0,1)(0,1){2}{\line(1,0){3}}\epic}}
\nc{\pdv}{\raisebox{-1\unitlength}[2\unitlength]{\bpic(1,2)\multiput(0,0)(0,1){3}{\line(1,0){1}}\multiput(0,0)(1,0){2}{\line(0,1){2}}\epic}}
\nc{\px}{\raisebox{-1\unitlength}[2\unitlength]{\bpic(2,2)\multiput(0,0)(1,0){2}{\line(0,1){2}}\multiput(0,1)(0,1){2}{\line(1,0){2}}
\put(0,0){\line(1,0){1}}\put(2,1){\line(0,1){1}}\epic}}
\nc{\peuu}{\raisebox{-2\unitlength}[2\unitlength]{\bpic(0.4,3)\put(0.2,2.5){\p{}{\emptyset}}\epic}}
\nc{\psuu}{\raisebox{-2\unitlength}[2\unitlength]{\bpic(1,3)\multiput(0,2)(1,0){2}{\line(0,1){1}}\multiput(0,2)(0,1){2}{\line(1,0){1}}\epic}}
\nc{\pdvu}{\raisebox{-2\unitlength}[2\unitlength]{\bpic(1,3)\multiput(0,1)(0,1){3}{\line(1,0){1}}\multiput(0,1)(1,0){2}{\line(0,1){2}}\epic}}
\nc{\ptv}{\raisebox{-2\unitlength}[2\unitlength]{\bpic(1,3)\multiput(0,0)(0,1){4}{\line(1,0){1}}\multiput(0,0)(1,0){2}{\line(0,1){3}}\epic}}
\begin{table}[h]\centering\footnotesize$\ba{|c@{\;\;\,}c@{\;\;\;}c@{\;\;\;}c@{\;\;\;}c@{\;\;\;}c@{\;\;\,}c|crrrrc|c|c|c|}
\hline\rule{0mm}{4ex}
\e_0&\e_1&\e_2&\e_3&\e_4&\e_5&\e_6&\,\O_1\!\!\!&\O_2\!\!\!&\O_3\!\!\!&\O_4\!\!\!&\O_5\!\!\!&\O_6\!\!&\A(\prec_0,\e)&\A(<,\e)&\A(>,\e)\\[3mm]\hline
\pe&\ps&\pe&\ps&\pe&\ps&\pe&0&0&0&0&0&0&\emptyset&\emptyset&\emptyset\\[1.2mm]\hline
\pe&\ps&\pd&\ps&\pe&\ps&\pe&0&1&1&0&0&0&\{3\}&\{1\}&\{3\}\\[1.4mm]\hline
\peu&\psu&\pdv&\psu&\peu&\psu&\peu&0&\mi1&\mi1&0&0&0&\{3\}&\{3\}&\{1\}\\[3.4mm]\hline
\pe&\ps&\pe&\ps&\pd&\ps&\pe&0&0&0&1&1&0&\{5\}&\{3\}&\{5\}\\[1.4mm]\hline
\peu&\psu&\peu&\psu&\pdv&\psu&\peu&0&0&0&\mi1&\mi1&0&\{5\}&\{5\}&\{3\}\\[3.4mm]\hline
\pe&\ps&\pd&\ps&\pd&\ps&\pe&0&1&1&1&1&0&\{5\}&\{1\}&\{5\}\\[1.4mm]\hline
\peu&\psu&\pdv&\psu&\pdv&\psu&\peu&0&\mi1&\mi1&\mi1&\mi1&0&\{5\}&\{5\}&\{1\}\\[3.4mm]\hline
\peu&\psu&\pdu&\psu&\pdv&\psu&\peu&0&1&1&\mi1&\mi1&0&\{5\}&\{1,5\}&\{3\}\\[3.4mm]\hline
\peu&\psu&\pdv&\psu&\pdu&\psu&\peu&0&\mi1&\mi1&1&1&0&\{3,5\}&\{3\}&\{1,5\}\\[3.4mm]\hline
\peu&\psu&\pdu&\px&\pdu&\psu&\peu&0&1&\mi1&\mi1&1&0&\{4,5\}&\{1,4\}&\{2,5\}\\[3.4mm]\hline
\peu&\psu&\pdv&\px&\pdv&\psu&\peu&0&\mi1&1&1&\mi1&0&\{2,5\}&\{2,5\}&\{1,4\}\\[3.4mm]\hline
\peu&\psu&\pdu&\px&\pdv&\psu&\peu&0&1&\mi1&1&\mi1&0&\{3,5\}&\{1,3,5\}&\{2,4\}\\[3.4mm]\hline
\peu&\psu&\pdv&\px&\pdu&\psu&\peu&0&\mi1&1&\mi1&1&0&\{2,4,5\}&\{2,4\}&\{1,3,5\}\\[3.4mm]\hline
\pe&\ps&\pd&\pt&\pd&\ps&\pe&0&1&2&2&1&0&\{4,5\}&\{1,2\}&\{4,5\}\\[1.4mm]\hline
\peuu&\psuu&\pdvu&\ptv&\pdvu&\psuu&\peuu&0&\mi1&\mi2&\mi2&\mi1&0&\{4,5\}&\{4,5\}&\{1,2\}\\[5.6mm]\hline
\ea$\caption{\protect\rule{0ex}{2.5ex}The oscillating tableaux of shape $\emptyset$ and length $6$.\label{ASMEx}}\end{table}

It can be seen from Table~\ref{ASMEx} that
\begin{equation}\label{ExN}\mathcal{N}(\emptyset,6,\prec_0,i)=\mathcal{N}(\emptyset,6,<,i)=\mathcal{N}(\emptyset,6,>,i)\;=\;
\left\{\ba{ll}1,&i=0\\[1.5mm]7,&i=1\\[1.5mm]6,&i=2\\[1.5mm]1,&i=3\\[1.5mm]0,&\mbox{otherwise\,,}\ea\right.\end{equation}
so using this in~(\ref{ExNOP}) gives
\begin{equation}|\OP(n,n,[n],[n],6)|\;=\;
\left(\,\ba{c}n\\[1mm]6\ea\,\right)
+7\left(\,\ba{c}n\pl1\\[1mm]6\ea\,\right)
+6\left(\,\ba{c}n\pl2\\[1mm]6\ea\,\right)
+\left(\,\ba{c}n\pl3\\[1mm]6\ea\,\right).\end{equation}
Finally, note that in addition to the equalities between the three distribution functions of~(\ref{ExN}),~(\ref{NEq}) also
gives $\mathcal{N}(\emptyset,6,\prec_q,i)=\mathcal{N}(\emptyset,6,<,i)=\mathcal{N}(\emptyset,6,>,i)$ for any
$q$-order $\prec_q$, $q\in\Z$ and $i\in[0,5]$.

\subsubsection{Nonintersecting Paths and Semistandard Young Tableaux}
In this section, the special case of tuples of nonintersecting paths in an $a$ by $b$ rectangle, with the paths starting and
ending at points (specified by $(\a,\b)\in\BP(a,b)$) along the lower and right boundaries, and
taking only unit steps rightwards or upwards, will be considered.  It will be seen that,
using the lexicographic ordering pair, there is
a bijection between such path tuples and certain semistandard Young tableaux.
Three different enumeration formulae will then be obtained, the derivations of two of which involve using
the Lindstr\"{o}m-Gessel-Viennot theorem.

A tuple of nonintersecting paths, $P\in\OP(a,b,\a,\b)$, is a path tuple which has no osculations, i.e.,
$X(P)=\emptyset$, $\chi(P)=0$ and $Z(P)=N(P)$, and so, by~(\ref{Pmag}), $|Z(P)|=|\l_{a,b,\a,\b}|$.
Accordingly, the set of tuples of nonintersecting paths in $[a]\t[b]$ with boundary point pair
$(\a,\b)\in\BP(a,b)$ is defined as
\begin{equation}\NP(a,b,\a,\b)\::=\:\OP(a,b,\a,\b,|\l_{a,b,\a,\b}|)\,.\end{equation}
It also follows from~(\ref{DChi}) that for any $S\in\VOS(a,b,\a,\b,|\l_{a,b,\a,\b}|)$ and any progression~$s$ for $S$,
the oscillating tableau $\e(s)$ associated with~$s$ has no deletions.  Therefore, as indicated in Section~11,
$\e(s)$ corresponds to a standard Young tableau.
Of primary interest here will be the progression $s=((i_1,j_1),\ldots,(i_l,j_l))$ for~$S$ induced
by the lexicographic ordering pair $(\L,<)$ (i.e., the lexicographic ordering of~$S$), where $l=|\l_{a,b,\a,\b}|$.
In this case, the generalized oscillating tableaux associated with $S$
via $(\L,<)$ is, using~(\ref{lex}) and~(\ref{the11}),
$\Theta_{(\L,<)}(S)=((\L(s_1),\ldots,\L(s_l)),\,\e(s))$
$=((i_1,\ldots,i_l),\,\e(s))\in\GOT([a],<,\l_{a,b,\a,\b},l)$.

\begin{figure}[h]\centering\bpic(0.3,3.6)\put(0,0.15){\pp{l}{4}}\put(0,1.15){\pp{l}{3}}\put(0,2.15){\pp{l}{2}}\put(0,3.15){\pp{l}{1}}
\put(0.45,3.45){\pp{b}{1}}\put(1.45,3.45){\pp{b}{2}}\put(2.45,3.45){\pp{b}{3}}
\put(3.45,3.45){\pp{b}{4}}\put(4.45,3.45){\pp{b}{5}}\put(5.45,3.45){\pp{b}{6}}\epic
\includegraphics[width=50mm]{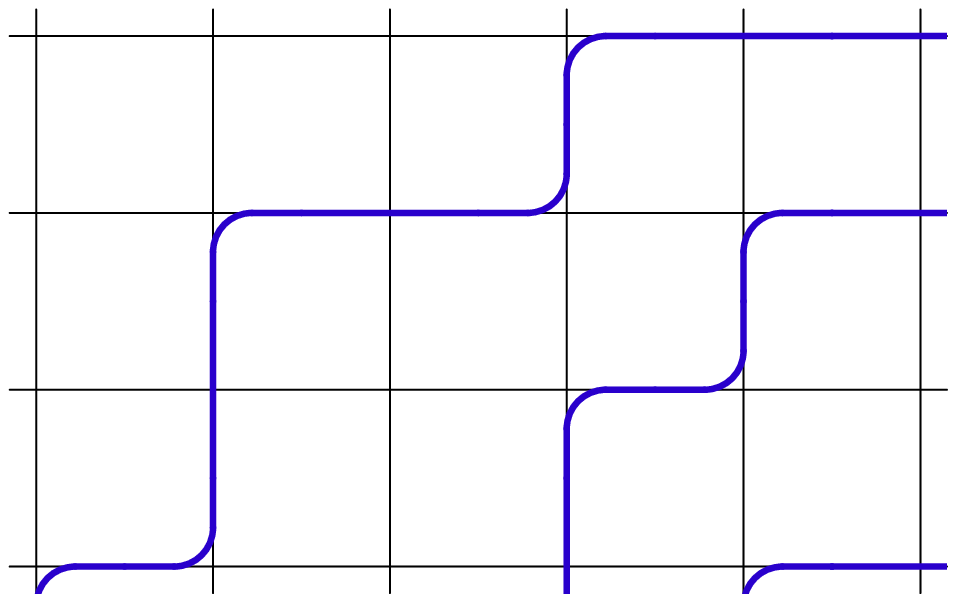}
\caption{Example of a tuple of nonintersecting paths.\label{NP}}\end{figure}

An example of a tuple $P$ of nonintersecting paths in $\NP(4,6,\{1,2,4\},\{1,4,5\})$
is shown diagrammatically in Figure~\ref{NP}.
\setlength{\unitlength}{9.43mm}
The progression for the vacancy-osculation set~$S$ of~$P$ induced by the lexicographic
ordering pair (i.e., the lexicographic ordering of the vacancies of $P$) is
$s=((i_1,j_1),\ldots,(i_8,j_8))=((1,1)$, $(1,2)$, $(1,3)$, $(2,1)$, $(3,1)$, $(3,3)$, $(3,6)$, $(4,3))$.
It can be seen that $(j_1\mi i_1,\ldots,j_8\mi i_8)=(0,1,2,-1,-2,0,3,-1)$ matches
the profile of the oscillating tableau of Table~\ref{SYTOT}, implying that this is the oscillating
tableau $\e(s)$ associated with $s$, and that the corresponding standard Young tableau is~$\s$ of~(\ref{ExSYT}).
The generalized oscillating tableau associated with $S$ via $(\L,<)$, $\Theta_{(\L,<)}(S)=((i_1,\ldots,i_8),\,\e(s))$, is therefore
\begin{equation}\label{ExSSYTGOT}\ba{l}(t,\e)\,=\,\Bigl((1,1,1,2,3,3,3,4),\\[1.2mm]
\hspace{31mm}(\emptyset,(1),(2),(3),(3,1),(3,1,1),(3,2,1),(4,2,1),(4,2,2))\Bigr).\ea\end{equation}

For a subset $T$ of~$\P$, and a partition~$\l$, define the corresponding set of \emph{semistandard Young tableaux} as
\begin{equation}\ba{r@{}l}\SSYT(T,\l)\;:=\;\Bigl\{\tau\,\Big|\,\:&
\bullet \ \tau\mbox{ is an array of entries }\tau_{ij}\in T\mbox{ with }(i,j)\in Y(\l)\\[1.5mm]
&\bullet \ \tau_{ij}\le\tau_{i,j\pl1}\mbox{ \ for each }i\in[\ell(\l)], \ j\in[\l_i\mi1]\\[1.5mm]
&\bullet \ \tau_{ij}<\tau_{i\pl1,j}\mbox{ \ for each }i\in[\ell(\l)\mi1], \ j\in[\l_{i\pl1}]\Bigr\}\,.\ea\end{equation}
As with standard Young tableaux, $\tau\in\SSYT(T,\l)$ is usually shown by writing each~$\tau_{ij}$ within the
unit square centered at $(i,j)$ in the Young diagram of $\l$.  These entries now increase weakly along each row
and increase strictly down each column.

It can be checked that, for any $n\in\P$ and partition $\l$,
each generalized oscillating tableau $(t,\e)\in\GOT([n],<,\l,|\l|)$ can be bijectively associated
with a semistandard Young tableau $\tau\in\SSYT([n],\l)$ where, for each $k\in[|\l|]$,
if $Y(\e_k)=Y(\e_{k\mi1})\cup\{(i,j)\}$, then $\tau_{ij}=t_k$.
For the inverse mapping from $\tau$ to $(t,\e)$, the points of $Y(\l)$ should first be
ordered as $(i_1,j_1),\ldots,(i_{|\l|},j_{|\l|})$ such that
$\tau_{i_kj_k}<\tau_{i_{k\pl1}j_{k\pl1}}$, or $\tau_{i_kj_k}=\tau_{i_{k\pl1}j_{k\pl1}}$ and
$j_k\mi i_k<j_{k\pl1}\mi i_{k\pl1}$, for each $k\in[|\l|\mi1]$. Then $(t,\e)$ is
given by $t_k=\tau_{i_kj_k}$ and $\O(\e)_k=j_k\mi i_k$ for each $k\in[|\l|]$.
For example, the semistandard Young tableau which corresponds to
the generalized oscillating tableau of~(\ref{ExSSYTGOT}) is
\setlength{\unitlength}{4.2mm}
\begin{equation}\label{ExSSYT}\raisebox{-1.27\unitlength}[2\unitlength][1.3\unitlength]{
\bpic(2.8,3)\put(0,1.5){\p{l}{\tau\;=}}\epic
\bpic(4,3)\multiput(0,0)(1,0){3}{\line(0,1){3}}\multiput(0,0)(0,1){2}{\line(1,0){2}}
\multiput(0,2)(0,1){2}{\line(1,0){4}}\multiput(3,2)(1,0){2}{\line(0,1){1}}
\put(0.5,0.5){\pp{}{3}}\put(0.5,1.5){\pp{}{2}}\put(0.5,2.5){\pp{}{1}}
\put(1.5,0.5){\pp{}{4}}\put(1.5,1.5){\pp{}{3}}\put(1.5,2.5){\pp{}{1}}
\put(2.5,2.5){\pp{}{1}}\put(3.5,2.5){\pp{}{3}}\epic}\,.\end{equation}

It follows from this bijection between $\GOT([n],<,\l,|\l|)$ and $\SSYT([n],\l)$, the bijection
between $\OT(\l,|\l|)$ and $\SYT(\l)$, and~(\ref{NGOT}) that
\begin{equation}\label{NSSYT}
|\SSYT([n],\l)|\;=\,\ds\sum_{\s\in\SYT(\l)}\left(\;\ba{c}n+\AA(\s)\\[1.3mm]|\l|\ea\;\right),\end{equation}
where $\AA$ is the number of ascents for standard Young tableaux, as defined in Section~11.
This result, which could have been derived directly rather than by using oscillating tableaux, is contained
in~\cite[Proposition~7.19.12]{Sta99}. Note that there also exists a simpler formula for $|\SSYT([n],\l)|$,
namely the \emph{hook-content formula} (see for example~\cite[Corollary~7.21.4]{Sta99}).

Now define, for any $S\in\VOS(a,b,\a,\b,|\l_{a,b,\a,\b}|)$,
$\Upsilon(S)$ to be the semistandard Young tableau associated with the generalized oscillating tableau $\Theta_{(\L,<)}(S)$.
Thus, for the vacancy-osculation set $S$ corresponding to the path tuple of Figure~\ref{NP}, $\Upsilon(S)$
is given by~(\ref{ExSSYT}).

Also define, for any $a,b\in\P$ and $\l\in\Par$,
\begin{equation}\label{SSSYT}
\SSSYT(a,b,\l)\;:=\;\{\tau\in\SSYT([a],\l)\;|\;\tau_{i,\l_i}\le b\mi\l_i\pl i\mbox{ \ for each }i\in[\ell(\l)]\}\,.
\end{equation}
Note that the condition $\tau_{i,\l_i}\le b\mi\l_i\pl i$\ru{1.5} places an additional
constraint on the largest (i.e., rightmost) entry of row $i$ of $\tau$ whenever $\l_i\mi i>b\mi a$.

By Theorem~11, $\Upsilon$ is an injective mapping from
$\VOS(a,b,\a,\b,|\l_{a,b,\a,\b}|)$ to $\SSYT([a],$
$\l_{a,b,\a,\b})$, and it can be shown
using Lemma~12 and arguments similar to those used in previous sections, that~$\Upsilon$ is a bijection from
$\VOS(a,b,\a,\b,|\l_{a,b,\a,\b}|)$ to $\SSSYT(a,b,\l_{a,b,\a,\b})$.
Also, by Theorem~16, if $\a=[a]$, then
$\SSSYT(a,b,\l_{a,b,[a],\b})=\SSYT([a],\l_{a,b,[a],\b})$.

Note too that in addition to $\Upsilon$, which involves the lexicographic ordering pair and gives a bijection
from $\VOS(a,b,\a,\b,|\l_{a,b,\a,\b}|)$ to a set of semistandard Young tableaux, there is also an analogous bijection
which involves the antilexicographic ordering pair, and gives a bijection from
$\VOS(a,b,\a,\b,|\l_{a,b,\a,\b}|)$ to a set of row-strict Young tableaux (i.e., arrays whose
entries increase strictly along each row and increase weakly down each column).

It now follows using the bijection $\Upsilon$ and Corollary~14 that
\begin{equation}\label{NNP1}\ba{l}|\NP(a,b,\a,\b)|\;=\;|\SSSYT(a,b,\l_{a,b,\a,\b})|\;=\\[1mm]
\hspace{49mm}\ds\sum_{\e\in\OT(\l_{a,b,\a,\b},|\l_{a,b,\a,\b}|)}\left(\;\ba{c}\min(a,b)+
\AA(\prec_{b\mi a},\e)\\[1.3mm]|\l_{a,b,\a,\b}|\ea\;\right),\ea
\end{equation}
for any $(b\mi a)$-order.  The sum over oscillating tableaux in~(\ref{NNP1}) can instead be regarded as a sum over
standard Young tableaux $\SYT(\l_{a,b,\a,\b})$. For the case $\a=[a]$, it follows that $|\NP(a,b,[a],\b)|=|\SSYT([a],\l_{a,b,[a],\b})|$,
so $\NP(a,b,[a],\b)$ can alternatively be enumerated using~(\ref{NSSYT}) or the hook-content formula.

The Lindstr\"{o}m-Gessel-Viennot theorem (see for
example~\cite[Theorem~1]{GesVie85},~\cite[Corollary~2]{GesVie89} or~\cite[Theorem~2.7.1]{Sta86}) states that
the number of tuples of nonintersecting paths with certain fixed endpoints can be expressed as the determinant of a
matrix of binomial coefficients.
For $r$ paths, i.e., $|\a|=|\b|=r$, this theorem then gives
an alternative formula for $|\NP(a,b,\a,\b)|$ as the determinant of an $r\t r$ matrix,
\begin{equation}\label{NNP2}
|\NP(a,b,\a,\b)|\;=\;\det\left(\left(\,\ba{c}i\pl j\\[1.3mm]i\ea\,\right)\right)_{\!\!i=a\mi\a_1,\ldots,a\mi\a_r,\;
j=b\mi\b_1,\ldots,b\mi\b_r}\;,
\end{equation}
where as usual the elements of $\a$ and $\b$ are labeled in increasing order.

Finally, as indicated in Section~1, semistandard Young tableaux can also be associated with tuples of nonintersecting paths
using a method in which each row of a tableau is regarded as a partition, by reading the row's entries from right to left,
and then associated with the path formed by the lower and right
boundary edges of the Young diagram of that partition.  By translating these paths appropriately on the lattice,
this gives a bijection between
$\SSSYT(a,b,\l_{a,b,\a,\b})$ and a certain set of $\ell(\l_{a,b,\a,\b})$-tuples of nonintersecting paths,
which is therefore in bijection with $\NP(a,b,\a,\b)$.
Using~\cite[Theorem~3]{GesVie89}, which combines this bijection with~\cite[Corollary~2]{GesVie89}, then gives
a third formula for $|\NP(a,b,\a,\b)|$ as the determinant of an $\ell(\l_{a,b,\a,\b})\t\ell(\l_{a,b,\a,\b})$ matrix,
\begin{equation}\label{NNP3}\ba{l}|\NP(a,b,\a,\b)|\;=\\[1.5mm]
\qquad\det\left(\left(\,\ba{c}(\l_{a,b,\a,\b})_i\mi i\pl j\pl\min(b\mi(\l_{a,b,\a,\b})_i\pl i,a)\mi1\\[1.3mm]
(\l_{a,b,\a,\b})_i\mi i\pl j\ea\;\right)\right)_{i,j=1,\ldots,\ell(\l_{a,b,\a,\b})}\;.\ea\end{equation}

\subsubsection{Discussion}
All of the results of this paper have now been presented, so in this final section, the scope for additional,
related work will be discussed.

Some possible directions for further research are:
\vspace{-4mm}
\begin{enumerate}
\item Studying the combinatorics of osculating paths with external configurations other than fixed start and end points
on the lower and right boundaries of a rectangle.\vspace{-2mm}
\item Investigating whether the generalized oscillating tableaux defined in this paper can be interpreted
within a representation theoretic context.\vspace{-2mm}
\item Applying further bijections to generalized oscillating tableaux, with the aim of
obtaining bijections between osculating paths and other known combinatorial objects.
\end{enumerate}
\vspace{-4mm}
With regard to (1), some obvious possibilities for other external configurations are fixed start points
on the lower and left boundaries, and end points on the upper and right boundaries, of a rectangle, boundary conditions which
correspond to previously-studied subclasses of alternating sign matrices other than those outlined in Section~4
(see for example \cite{Kup02,Oka06}), boundary conditions which correspond to previously-studied cases of the
six-vertex model of statistical mechanics (see for example \cite{BatBaxOroYun95,Bax82,IzeCokKor92,OwcBax89,YunBat95,ZinJus02}), and the
configuration used in studies of friendly walkers in which
the paths start and end on parallel lines at $45^\circ$ to the rows
or columns of the lattice (see for example \cite{Bou06,Bra97,Ess03,GutVog02,KraGutVie03}).

With regard to (2), oscillating tableaux have often been studied in the context of
the representation theory of the symplectic group or Brauer algebra (see for example~\cite{Ber86,DelDulFav88,Pro91,Sun86,Sun90,Ter01}).
Certain generalizations of oscillating tableaux have also been introduced (see for
example~\cite{ChaDul02,Cho03,DulSag95,Ges93,Kra96,Rob95}), but these all appear to differ from the generalization
used in this paper.  It would therefore seem worthwhile to investigate whether the generalized oscillating tableaux used here
can be related to representation theory.  Further motivation for searching for representation-theoretic interpretations
of aspects of this paper is that it is shown in~\cite{Oka06} that enumeration formulae for standard
alternating sign matrices, and various subclasses thereof,
can be obtained using determinantal formulae which are related to characters of irreducible representations of classical groups.

With regard to (3), it would seem natural to combine the bijections obtained in this paper, between
tuples of osculating paths and generalized oscillating tableaux, with further bijections involving oscillating tableaux.
Certain such bijections, which can be regarded as analogs of the Robinson-Schensted-Knuth correspondence (see for
example~\cite{Ful99},~\cite{Sag01} and~\cite[Ch.~7]{Sta99}), are already known (see for
example~\cite{Ber86,ChaDul02,Cho03,DelDulFav88,DulSag95,Pro91,Rob91b,Rob95,RobTer05,Sun86,Sun90}), although
these need to be generalized in order to be applicable to the generalized oscillating tableaux being used here.
It is expected that the results of such work will be reported in a sequel to this paper.  In particular, it is hoped that
it will be possible to exhibit bijections between certain tuples of osculating paths and certain tuples of nonintersecting paths.

In conclusion, an example of such a result, and its consequences, will now be outlined.
For any $n\in\P$ and $l\in\N$,
there appears to be a bijection between the set of tuples of osculating paths
$\OP(n,n,[n],[n],2l)$, and the set of tuples of nonintersecting paths
$\cup_{\l{\sss\models}\,l,\l_1<n}\,\NP(n,n,\a_{n,n,D(\l)},\b_{n,n,D(\l)})$, where~(\ref{BPInv}) is being used,
$\cup_{\l{\sss\models}\,l,\l_1<n}$ denotes the union over all strict (i.e., with distinct parts) partitions $\l$ of $l$ with
largest part less than $n$,
and $D(\l)$ is the double of~$\l$, defined using Frobenius notation to be the partition
$(\l_1,\ldots,\l_r\,|\,\l_1\mi1,\ldots,\l_r\mi1)$, where $r$ is the length of $\l$.
Note that $\OP(n,n,[n],[n],2l)$ corresponds to the set of $n\times n$ standard alternating sign matrices whose
corresponding path tuple has~$l$ vacancies and $l$ osculations.
Such a bijection seems related to the analog of
the Robinson-Schensted-Knuth correspondence described in~\cite[Sec.~4]{Bur74}.
Applying~(\ref{NNP2}) to the individual sets of nonintersecting paths gives
a formula for $|\OP(n,n,[n],[n],2l)|$ as a sum of determinants, and then summing over all
possible values of $l$, with weight~$x^l$ for a parameter $x$, and combining these determinants into a single determinant, gives
\begin{equation}\ba{l}\ds\sum_{l=0}^{n(n\mi1)/2}\!|\OP(n,n,[n],[n],2l)|\:x^l\;=\\[1.4mm]
\quad\qquad\qquad\qquad\qquad\qquad\qquad x^{n(n\mi1)/2}\:\det\left(x^{-i}\left(\,\ba{c}i\pl j\\[1.3mm]i\ea\,\right)-
\d_{i,j\pl1}\right)_{\!\!i,j=0,\ldots,n\mi1}\,.\ea\end{equation}
Setting $x=1$ gives a determinantal formula for $|\OP(n,n,[n],[n])|$, or equivalently the number of
$n\t n$ standard alternating sign matrices, which was obtained
in~\cite[Remark~5.2]{GesXin06} by transforming the determinant used in~\cite{And79} for the
enumeration of descending plane partitions.  However, it is hoped that the method being used here will clarify the
underlying combinatorics.

\end{document}